\documentclass[final]{siamltex}

\usepackage{epsfig,amssymb,latexsym}
\usepackage{amsfonts,psfrag,amsmath,bbm,color}
\usepackage{cancel}
\usepackage{comment}
\usepackage{mathrsfs}
\usepackage{graphicx}
\usepackage{textcomp}
\usepackage{multirow}
\usepackage{enumerate}
\usepackage{cancel}
\usepackage{algpseudocode}
\usepackage{caption}  
\usepackage{subcaption}
\usepackage{url}
\usepackage{rotating}
\usepackage{bigints}
\usepackage{hyperref}
\usepackage{xcolor}
\usepackage{ulem}
\usepackage{tabularx}

\usepackage{placeins}

\newtheorem{Lemma}{Lemma}
\newtheorem{Th}{Theorem}
\usepackage[ruled,vlined]{algorithm2e}
\usepackage{geometry}
\vbadness=\maxdimen

\newcommand{\footremember}[2]{%
    \footnote{#2}
    \newcounter{#1}
    \setcounter{#1}{\value{footnote}}%
}
\newcommand{\footrecall}[1]{%
    \footnotemark[\value{#1}]%
} 
\renewcommand{\figurename}{Figure}
\usepackage[most]{tcolorbox}

\newcolumntype{L}{>{\raggedright\arraybackslash}X}

\usepackage{tikz}
\usetikzlibrary{shapes.geometric, arrows}
\tikzstyle{startstop} = [rectangle, rounded corners, minimum width=1cm, minimum height=1cm,text centered, draw=black]
\tikzstyle{io} = [trapezium, trapezium left angle=70, trapezium right angle=110, minimum width=1cm, minimum height=1cm, text centered, draw=black, fill=blue!30]
\tikzstyle{method} = [rectangle, rounded corners, minimum width=1cm, minimum height =1cm, text centered, draw=black]
\tikzstyle{process} = [rectangle, minimum width=1cm, minimum height=1cm, text centered, draw=black]
\tikzstyle{decision} = [diamond, minimum width=0.5cm, minimum height=0.5cm, text centered, draw=black, fill=green!30]
\tikzstyle{arrow} = [thick,->,>=stealth]

\usepackage{amscd}
\usepackage{xcolor}
\graphicspath{{./figs/}}

\def\PP{{{\rm l}\kern - .15em {\rm P} }}
\def\PN2{{\PP_{N}-\PP_{N-2}}}

\newcommand{\bn}{\boldsymbol{n}}

\newcommand{\bx}{\boldsymbol{x}}

\newcommand{\deleted}[1]{{}}

\begin{document}
\title{Stocking and Harvesting Effects in Advection-Reaction-Diffusion Model: Exploring Decoupled Algorithms and Analysis}

\author{
Mayesha Sharmim Tisha\footremember{DU}{D\MakeLowercase{epartment of} M\MakeLowercase{athematics}, U\MakeLowercase{niversity of} D\MakeLowercase{haka}, D\MakeLowercase{haka} 1000, B\MakeLowercase{angladesh}.}\hspace{-1ex}
\and Md. Kamrujjaman\footrecall{DU}\hspace{1ex}\footnote{p\MakeLowercase{artially supported by the }U\MakeLowercase{niversity} G\MakeLowercase{rants} C\MakeLowercase{ommission} (UGC), \MakeLowercase{and the}  U\MakeLowercase{niversity of }D\MakeLowercase{haka}, B\MakeLowercase{angladesh}.}%
\and
Muhammad Mohebujjaman\footremember{uabm}{D\MakeLowercase{epartment of} M\MakeLowercase{athematics}, U\MakeLowercase{niversity of} A\MakeLowercase{labama at} B\MakeLowercase{irmingham}, AL 35294, USA; P\MakeLowercase{artially supported by the} N\MakeLowercase{ational} S\MakeLowercase{cience} F\MakeLowercase{oundation} (NSF) \MakeLowercase{grant} DMS-2213274. C\MakeLowercase{orresponding author.} E\MakeLowercase{mail: mmohebuj@uab.edu}.}\hspace{-1ex}
 \and Taufiquar Khan\footremember{UNCC}{D\MakeLowercase{epartment of }M\MakeLowercase{athematics and} S\MakeLowercase{tatistics}, U\MakeLowercase{niversity of} N\MakeLowercase{orth} C\MakeLowercase{arolina at} C\MakeLowercase{harlotte}, NC 28223, USA.}
 }

\maketitle

\pagestyle{myheadings}
\thispagestyle{plain}

\markboth{\MakeUppercase{Stocking and Harvesting Effects in Advection-Reaction-Diffusion Model}}{\MakeUppercase{M. S. Tisha, M. Kamrujjaman, M. Mohebujjaman, and T. Khan}}

\begin{abstract}
    We propose a time-dependent Advection Reaction Diffusion (ARD) $N$-species competition model to investigate the Stocking and Harvesting (SH) effect on population dynamics. For ongoing analysis, we explore the outcomes of a competition between two competing species in a heterogeneous environment under no-flux boundary conditions, meaning no individual can cross the boundaries. We establish results concerning the existence, uniqueness, and positivity of the solution. As a continuation, we propose, analyze, and test two novel fully discrete decoupled linearized algorithms for a nonlinearly coupled ARD $N$-species competition model with SH effort. The time-stepping algorithms are first and second order accurate in time and optimally accurate in space. Stability and optimal convergence theorems of the decoupled schemes are proved rigorously. We verify the predicted convergence rates of our analysis and the efficacy of the algorithms using numerical experiments and synthetic data for analytical test problems. We also study the effect of harvesting or stocking and diffusion parameters on the evolution of species population density numerically and observe the coexistence scenario subject to optimal stocking or harvesting.
\end{abstract}



\section{Introduction} 
	Stocking and harvesting modeling and their strategies for interactions between species have a significant impact on population economics, in particular, on the productivity and profitability of farms \cite{clark1974mathematical,dai1998coexistence,li2011optimal,otunuga2021time,yang2019optimal,zhang2003optimal, schons2021rotatinuous, brissette1996effects, forsberg1996optimal, lorenzen2005population, shaner1996assessment}. Harvesting is a management strategy involving the reduction of the size of the species by hunting, fishing, or other forms of capture, and natural disasters. The release of cultured organisms into the environment to increase the number of wild species is known as stocking (or immigration).
	
	In recent decades, numerous researchers have conducted investigations, mathematical modeling, and analyses of SH effects based on ordinary differential equations \cite{myerscough1992analysis, brauer1982coexistence,brauer1982constant, liu2007impulsive, hening2022effects, liu2010coexistence, dong2006extinction, jiao2008delayed, wang2014permanence, jiao2008permanence, ruan2023imperfect}, difference equations \cite{selgrade1998reversing, alsharawi2009coexistence, selgrade1998using}, delay-differential equations \cite{delay1,delay2,delay3,delay4}, and Partial Differential Equations (PDEs) \cite{braverman2009optimal, korobenko2013persistence,kamrujjaman2022spatio, adan2023interplay, mohebujjaman2024decoupled}. In \cite{ainseba2008reaction}, PDEs Lotka-Volterra interactions model with no-flux boundary conditions in the presence of prey-taxis and spatial diffusion is given, and the existence and uniqueness of the weak solution are examined. In \cite{korobenko2013persistence, braverman2016competitive, roques2007population}, there are few studies about the harvesting of one or two populations, and in practice, they often fail to illustrate the actual circumstances. Additional intriguing situations that show either coexistence or competitive exclusion by others are discovered when harvesting is considered for multiple interacting populations \cite{clayton1997bringing,leung1995optimal,liu2016optimal, stigter2004optimal}. A predator-prey system's total behavior is explored in \cite{brauer1982constant} under constant stocking or harvesting of one or both species.
	
	PDEs Stocking and Harvesting (PDEs-SH) modelings are complex but more realistic. The spatial-temporal impacts of logistic and Gilpin-Ayala growth functions with starving type diffusion in a single species population with harvesting was investigated in \cite{kamrujjaman2022spatio, Zahan2022Mathmatical} using Reaction Diffusion PDEs-SH (RD-PDEs-SH) model. They investigated the properties of species stability in terms of survival and extinction and explored the optimal harvesting attempts in the situation of space-dependent carrying capacity. Two species RD-PDEs-SH model is proposed, analyzed the existence and uniqueness of the model solution, and stated the coexistence conditions in \cite{adan2023interplay}. $N-$species non-linearly coupled system of the RD-PDEs-SH model is presented in \cite{mohebujjaman2024decoupled} together with space- and time-dependent intrinsic growth rate and carrying capacity.
	
	Unidirectional drift \cite{Lutscher1}, often referred to as advection in mathematical ecology, pushes individuals out of the system, leading to a decline in population across various environments \cite{Lutscher2}. In recent years, ARD-PDEs have gained popularity for modeling problems in population ecology \cite{Zhou1,Xu1,Zhou2}. In cases where there is no control and only classical diffusion is at play, the suitable model is a diffusion equation without a reaction term.   When a species senses an environmental gradient and moves towards or away from it, the ARD-PDE model \cite{Cantrell1,Zhang1} is employed. This gradient is more competitive than the classical reaction-diffusion one due to the inclusion of the advection term.   	

 The ARD-PDEs Stocking and Harvesting (ARD-PDEs-SH) model for $N$-species competition  is governed by the following system of nonlinearly coupled time evolutionary equations \cite{adan2023interplay,braverman2019interplay, wong2009analysis}: For $i=1,2,\cdots,N$
\begin{align}
\frac{\partial u_i}{\partial t}=d_i\Delta u_i-\beta_i\nabla\cdot(u_i\nabla K)+r_iu_i\left(1-\gamma_i-\frac{1}{K}\sum\limits_{j=1}^N u_j\right)+f_i, \hspace{2mm}\forall (t, \bx)\in (0,T]\times\Omega,\label{RDE1}
\end{align}
along with known initial and boundary conditions, where $u_i$, $d_i$, $\beta_i$, $r_i(t,\bx)$, and $\gamma_i$ indicate the $i^{th}$ competing species' population density, diffusion rate, advection rate, intrinsic growth rate, and harvesting or stocking coefficient, respectively. The number of species in the competition is denoted by $N$; Single species simple logistic growth model can be retrieved for $N=1$. $K(t,\bx)$ represents the heterogeneous environment's carrying capacity, $f_i$ the forcing, $t$ the time, $\bx$ the spatial variable, $\Omega$ the domain, and $T$ the simulation end time. In the model \eqref{RDE1}, the stocking and harvesting rate is considered to be proportional to the intrinsic growth rate, and  $\gamma_i<1$. When $\gamma_i=1$, the growth function of the  $i^{th}$ species decreases. The stocking is represented as the case when $\gamma_i<0$.

In the continuous analysis, the primary objective of this study is to explore the steady state under various and unequal diffusion and advection rates through both theoretical and numerical analyses. Solving the above system is computationally expensive because we need to solve a system of non-linearly coupled ARD-PDEs at each time-step together with space- and time-dependent growth rate and carrying capacity. It is an open problem how to decouple the above ARD-PDEs-SH system \eqref{RDE1} in a stable way.

Finite Element (FE) analysis and simulations of PDEs are very popular in the scientific community, however, FE analysis and implementation of the PDEs-SH model are scarce in the literature.  FE analysis of a three-species competition-diffusion model with constant intrinsic growth rate in a homogeneous environment ($K\equiv$ constant) without SH effect is studied in \cite{wong2009analysis}. In \cite{adan2023interplay}, the optimal harvesting in regulating species density in a two-species RD-PDEs-SH model with a diverse habitat is examined using a FE fully-discrete backward-Euler decoupled linearized time-stepping algorithm without any discrete analysis. First- and second-order temporal and optimally accurate in space, fully discrete, two decoupled time-stepping algorithms of a system of RD-PDEs-SH model are given in \cite{mohebujjaman2024decoupled} and the stability, and convergence theorems are proven rigorously. The authors also investigated numerically the effect of SH on the coexistence of the species. The main novelties of our work are the following:
\begin{itemize}
\item We propose a non-stationary system of ARD-PDEs-SH model \eqref{RDE1} for the population dynamics of $N-$species competition with space- and time-dependent intrinsic growth rate and heterogeneous carrying capacity.
\item We establish and demonstrate the existence, uniqueness, and positivity of solutions for the ARD-PDEs-SH model \eqref{RDE1} with $N=1$, and 2 in a heterogeneous advective environment under no-flux boundary conditions. A unique solution exists when the harvesting rate is in the interval (0,1]. If a species' ratio of advection rate to diffusion rate is smaller, it will always dominate.
\item We propose two fully discrete decoupled linearized implicit-explicit time-stepping FE schemes for approximating \eqref{RDE1} together with appropriate initial and boundary conditions:
\begin{itemize}
	\item We rigorously prove the stability and convergence theorems of the discrete schemes. We found that the first- and second-order temporal schemes are optimally accurate in space and time.
	\item Several numerical experiments are given to examine the predicted convergence rates with manufactured analytical solutions.
\end{itemize}
\item Several tests are given to examine the effect of advection, diffusion, SH on the coexistence of the species.

\item To the best of our knowledge, the ARD-PDEs-SH model, its continuous analysis, the proposed two discrete schemes, and their analysis for modeling the impact of SH are novel.
\end{itemize}
The rest of the manuscript is organized as follows: To follow a smooth analysis we provide the necessary notations and mathematical preliminaries in Section \ref{notation-preliminaries}. The existence, uniqueness, and positivity of the solution are discussed in Section \ref{Existence-Uniqueness-Positivity-Solution} for single and double species models. Two fully discrete decoupled linearized FE schemes are proposed, and rigorously proven their stability and convergence theorems in Section \ref{fully-discrete-schme}. In Section \ref{numerical-test}, to support the theoretical discrete analysis, we examine the spatial and temporal convergence rates of both of the schemes for a three-species model. Several numerical tests are also given to examine the effect of diffusion, advection, and SH on the population dynamics in Section \ref{numerical-test}. Finally, in Section \ref{conclusion}, we provide the conclusion and discussion on future research directions.

\section{Notation and Preliminaries} \label{notation-preliminaries} Let $\Omega\subset\mathbb{R}^d (d\in\{1,2,3\})$ be a  convex domain with boundary $\partial\Omega$. For a given carrying capacity $K:(0,T]\times\Omega\rightarrow\mathbb{R}$, we define
\begin{align}
K_{\min}:=\inf\limits_{(t,{\bx})\in (0,T]\times\Omega}|K(t,\bx)|,\label{kmin-def}
\end{align}
and make the assumption that $K_{\min}>0$. The standard $L^2(\Omega)$ norm and inner product are denoted by $\|.\|$ and $(.,.)$. Similarly, the $L^p(\Omega)$ norms and the Sobolev $W_p^k(\Omega)$ norms for $k\in\mathbb{N},\hspace{1mm}1\le p\le \infty$ are $\|.\|_{L^p}$ and $\|.\|_{W_p^k}$, respectively. The Sobolev space $W_2^k(\Omega)^d$ is represented by Hilbert spaces $H^k(\Omega)^d$ with norm $\|.\|_k$.
For $X$ being a normed function space in $\Omega$, $L^p(0,T;X)$ is the space of all functions defined on $(0,T]\times\Omega$ for which the following norm 
\begin{align*}
\|u\|_{L^p(0,T;X)}=\left(\int_0^T\|u\|_{X}^pdt\right)^\frac{1}{p},\hspace{2mm}p\in[1,\infty)
\end{align*}
is finite. For $p=\infty$, the usual modification is used in the definition of this space. We denote $$\|u\|_{\infty,\infty}:=\|u\|_{L^\infty\big(0,T;L^\infty(\Omega)^d\big)}.$$ The natural function spaces for our problem are
\begin{align*}
X:&=H_0^1(\Omega)=\big\{v\in L^2(\Omega) :\nabla v\in L^2(\Omega)^{d},\hspace{1mm}  v=0 \hspace{1mm} \mbox{on}\hspace{1mm}   \partial \Omega\big\}.
\end{align*}
In the dual space of $X$ and for an element $f$ , we define the norm as
$$\|f\|_{-1}:=\sup\limits_{v\in X}\frac{(f,v)}{\|\nabla v\|}.$$
Recall the Poincar\'e inequality holds in $X$: There exists $C$ depending only on $\Omega$ satisfying for all $\phi\in X$,
\[
\| \phi \| \le C \| \nabla \phi \|.
\]
We define the initial conditions $u_i^0:=u_i(0,{\bx})$, and $u_i^1:=u_i(\Delta t,{\bx})$.
We obtain the following by multiplying both sides of \eqref{RDE1} by $v\in X$ and integrating over $\Omega$: Regarding $i=1,2,\cdots,N$
\begin{align}
\left(\frac{\partial u_i}{\partial t},v\right)&+d_i\left(\nabla u_i,\nabla v\right)-\beta_i\left(u_i\nabla K,\nabla v\right)=(1-\gamma_i)\big(r_iu_i,v\big)-
\left(\frac{r_i u_i}{K}\sum\limits_{j=1}^Nu_j,v\right)+\left(f_i,v\right).\label{vec-weak-form}
\end{align}
The complying finite element space is denoted by $X_h\subset X$, and the inverse inequality holds if the triangulation is sufficiently regular $\tau_h(\Omega)$, where $h$ is the maximum triangle diameter. We have the following approximation features in $X_h$ that are typical of piecewise polynomials of degree $k$: \cite{BS08,linke2019pressure}
\begin{align}
\|u- P^{L^2}_{X_h}(u) \|&\leq Ch^{k+1}|u|_{k+1},\hspace{2mm}u\in H^{k+1}(\Omega),\label{AppPro3}\\
\| \nabla (u- P^{L^2}_{X_h}(u)  ) \|&\leq Ch^{k}|u|_{k+1},\hspace{2mm}u\in H^{k+1}(\Omega),\label{AppPro4}
\end{align}
where $P^{L^2}_{X_h}(u)$ is the $L^2$ projection of $u$ into $X_h$ and $|\cdot|_r$ denotes the $H^r$ seminorm. Note that $C>0$ is a generic constant and changes in computation. The following lemma for the discrete Gr\"onwall inequality.
\begin{Lemma}\label{dgl} Let $\mathbb{N}$ denote the set of all natural numbers and
$\Delta t$, $\mathcal{E}$, $a_n$, $b_n$, $c_n$, $d_n$ be non-negative numbers for $n=1,\cdots\hspace{-0.35mm},M$ such that
$$a_M+\Delta t \sum_{n=1}^Mb_n\leq \Delta t\sum_{n=1}^{M-1}{d_na_n}+\Delta t\sum_{n=1}^Mc_n+\mathcal{E}\hspace{3mm}\mbox{for}\hspace{2mm}M\in\mathbb{N},$$
then for all $\Delta t> 0,$
$$a_M+\Delta t\sum_{n=1}^Mb_n\leq \left(\Delta t\sum_{n=1}^Mc_n+\mathcal{E}\right)\mbox{exp}\left(\Delta t\sum_{n=1}^{M-1}d_n\right)\hspace{2mm}\mbox{for}\hspace{2mm}M\in\mathbb{N}.$$
\end{Lemma}
In the following section, we want to state the existence, uniqueness, and positivity of the solutions for the problems with $N=1,$ and 2.
\section{Existence, Uniqueness and Positivity of Solution}\label{Existence-Uniqueness-Positivity-Solution} Without forcing consider the following model for two species competition with homogeneous Neumann boundary conditions: For $i=1,2$
\begin{align} 
	\label{harvest_system}
	\begin{cases}
  &\frac{\partial u_i}{\partial t}=d_i\Delta u_i-\beta_i\nabla\cdot(u_i\nabla K)+r_iu_i\left(1-\gamma_i-\frac{u_1+u_2}{K}\right), \hspace{2mm}\forall (t, \bx)\in (0,T]\times\Omega,\\
		\vspace{0.2cm}
		& d_i\frac{\partial u_i}{\partial \bn}-\beta_iu_i\frac{\partial K}{\partial \bn}=0,\  \bx\in\partial \Omega,\\
		\vspace{0.2cm}
		&	u_i(\bx,0)={{u_i}^0}(\bx)\geq,\not\equiv0,\  \bx\in \Omega,
	\end{cases}
\end{align}
where $\bn$ is an outward unit normal vector. To analyze the system \eqref{harvest_system}, we consider the following single species model:
\begin{align} 
\label{harvest_system_decoupled}
\begin{cases}
	\vspace{0.2cm}
	&\frac{\partial u_1}{\partial t}=d_1\Delta u_1-\beta_1\nabla\cdot(u_1\nabla K)+r_1u_1\left(1-\gamma_1-\frac{u_1}{K}\right), \hspace{2mm}\forall (t,{\bx})\in (0,T]\times\Omega,\\
	\vspace{0.2cm}
	& d_1\frac{\partial u_1}{\partial \bn}-\beta_1u_1\frac{\partial K}{\partial \bn}=0,\  \bx\in\partial \Omega,\\
	\vspace{0.2cm}
	&	u_1(\bx,0)={{u_1}^0}(\bx)\geq,\not\equiv0,\ \bx\in \Omega.
\end{cases}
\end{align}
\begin{Lemma}
\label{lm17}
Let  ${u_1}^0(\bx)$ be the non-negative and non-trivial initial continuous function and $ {u_1}^0(\bx) > 0 $	in some non-empty open bounded sub-domain $ \Omega_1 \subset \Omega $. Then for any $ t > 0 $, there exists a unique solution $ u_1(t, \bx) $ of the problem \eqref{harvest_system_decoupled} and it is positive.
\end{Lemma}
\begin{proof}
We define
$$\mathcal{F}(\bx,u_1):=\mathcal{G}(\bx,u_1)u_1=r_1u_1\left(1-\gamma_1-\frac{u_1}{K}\right),$$
where
$$\mathcal{G}(\bx,u_1):=r_1\left(1-\gamma_1-\frac{u_1}{K}\right).$$
Then, the system \eqref{harvest_system_decoupled} becomes
\begin{align} 
	\label{harvest_system_decoupled2}
	\begin{cases}
		\vspace{0.2cm}
		&\frac{\partial u_1}{\partial t}=d_1\Delta u_1-\beta_1\nabla\cdot(u_1\nabla K)+\mathcal{G}u_1,\hspace{2mm}\forall (t,{\bx})\in (0,T]\times\Omega,\\
		&d_1\frac{\partial u_1}{\partial \bn}-\beta_1u_1\frac{\partial K}{\partial \bn}=0,\ \displaystyle \bx\in\partial \Omega,\\
	\vspace{0.2cm}
	&	u_1(\bx,0)={{u_1}^0}(\bx)\geq,\not\equiv0,\ \bx\in \Omega.
	\end{cases}
\end{align}
Here, $\mathcal{F}(\bx,u_1)$ is Lipschitz in $u_1$ and is a measurable function in $\bx$ which is bounded if $u_1$ is restricted to a bounded set, $\Omega$ and $\partial\Omega$
is of class $ C^{2+{\alpha}}(\Omega)$, for some $\alpha>0$. The function $\mathcal{F}(\bx,u_1)$ and $\mathcal{G}(\bx,u_1)$ are in the class $ C^2 $ in $ u_1 $, and there exists $ K> 0 $ such that $\mathcal{G}(\bx,u_1)<0$
for $ u_1 > K $. The corresponding eigenvalue problem of \eqref{harvest_system_decoupled2} is represented as follows:
\begin{align}
	\label{ef}
	\begin{cases}
		& d_1\Delta\xi-\beta_1\nabla\cdot\left(\xi\nabla K\right)+\mathcal{G}(\bx,0)=\lambda\xi,\quad \bx\in \Omega,\\
		& d_1\frac{\partial \xi}{\partial \bn}-\beta_1\xi\frac{\partial K}{\partial \bn}=0,\quad \bx\in\partial\Omega,
	\end{cases}
\end{align}
where $\xi$ is an eigenfunction. The assumptions on $\mathcal{F}(\bx,u_1)$ imply that we can write $ \mathcal{F}(\bx,u_1) = (\mathcal{G}(\bx,0) + \mathcal{H}(\bx,u_1)u_1) u_1 $, where $\mathcal{H}(\bx,u_1)$ is $C^1$ in $u_1$.
If the principal eigenvalue $\lambda$ is positive of this problem. Let $\Phi$ be an eigenfunction for \eqref{ef} with $\Phi>0$ on $\Omega$. For $\epsilon>0$ sufficiently small,
\begin{align*}
	d_1\Delta(\epsilon\Phi)-\beta_1\nabla\cdot\left((\epsilon\Phi)\nabla K\right)+\mathcal{F}(\bx,\epsilon\Phi)&=\epsilon\left[d_1\nabla\Phi-\beta_1\Phi\nabla K+\mathcal{G}(\bx,0)\Phi\right]+\mathcal{H}(\bx,\epsilon\Phi)\epsilon^2\Phi^2\\
	&=\epsilon\lambda\Phi+\mathcal{H}(\bx,\epsilon\Phi)\epsilon^2\Phi^2\\
	&=\epsilon\Phi\left(\lambda+\mathcal{H}(\bx,\epsilon\Phi)\epsilon\Phi\right).
\end{align*}
Hence, for $\epsilon > 0$ small, $\epsilon\Phi$	 is a sub-solution for the elliptic problem
\begin{align}
	\label{ef2}
	\begin{cases}
		&d_1\Delta u_1-\beta_1\nabla\cdot(u_1\nabla K)+\mathcal{F}(\bx,u_1)=\lambda u_1,\quad \bx\in \Omega,\\
		& d_1\frac{\partial u_1}{\partial \bn}-\beta_1u_1\frac{\partial K}{\partial \bn}=0,\ \displaystyle \bx\in\partial \Omega,\\
		&u_1(\bx,0)={u_1}^0(\bx)\geq,\not\equiv0,\  \bx\in \Omega,
	\end{cases}
\end{align}
corresponding to \eqref{harvest_system_decoupled}. If $ \underline{u_1}(\bx,t) $ is a solution to \eqref{harvest_system_decoupled} with $\underline{u_1}(\bx,0) = \epsilon\Phi$, then at
$ t = 0 $, $\frac{\partial\underline{u_1}}{\partial t}>0$ on $\Omega$
and general properties of sub-solutions and super-solutions imply
that $\underline{u_1}(\bx,t)$ is increasing in $ t $. Since $ K > u_1 $ is a super-solution to \eqref{harvest_system_decoupled} we must have
$\underline{u_1}(\bx,t)\uparrow {u_1}^\ast(\bx) $ as $ t\rightarrow\infty $, where ${u_1}^\ast$ is the minimal positive solution of \eqref{harvest_system_decoupled}
(we can be sure that ${u_1}^\ast$ is minimal because $\Phi$ will be a strict sub-solution for all $ \epsilon> 0 $
sufficiently small). If $ {u_1}(\bx,t) $ is a solution to \eqref{harvest_system_decoupled} which is initially nonnegative and is
positive on an open subset of $\Omega$, then the strong maximum principle implies $ {u_1}(\bx,t) > 0 $ on $\overline{\Omega}$ for $ t > 0 $, which completes the proof.
\end{proof}
\begin{Lemma}\label{lm18}
A unique equilibrium solution of \eqref{harvest_system_decoupled} ${u_1}^\ast(\bx)$ exists. Then for any initial solution ${u_1}^0(\bx)\ge 0,\ {u_1}^0(\bx)\not\equiv 0$, the solution ${u_1}(\bx,t)$\ of \eqref{harvest_system_decoupled} satisfies the condition $$\lim_{t\rightarrow\infty}{u_1}(\bx,t)={u_1}^\ast(\bx)$$ uniformly for $\bx\in \overline{\Omega}$.
\end{Lemma}
\begin{proof}
Suppose that the hypotheses of Lemma \ref{lm17} are satisfied and $ \mathcal{F}(\bx,u_1) =
\mathcal{G}(\bx,u_1)u_1 $ with $\mathcal{G}(\bx,u_1)$ strictly decreasing in $ u_1 $ for $ u_1 \ge 0 $. Then, the minimal positive
equilibrium ${u_1}^\ast$ is the only positive equilibrium for \eqref{harvest_system_decoupled}. Let ${u_1}^{\ast\ast}$ be an another positive
equilibrium of \eqref{harvest_system_decoupled} with ${u_1}^{\ast\ast}\ne {u_1}^\ast$, then since $ {u_1}^\ast $ is minimal, we must have $ {u_1}^{\ast\ast}> {u_1}^\ast $
somewhere on $\Omega$. Since ${u_1}^\ast > 0 $ is an equilibrium of \eqref{harvest_system_decoupled}, it is a positive solution to
\begin{align}
	\label{ef3}
	\begin{cases}
		&d_1\Delta\xi-\beta_1\nabla\cdot(\xi\nabla K)+\mathcal{G}(\bx,{u_1}^\ast)=\lambda\xi,\quad \bx\in \Omega,\\
  & d_1\frac{\partial \xi}{\partial \bn}-\beta_1\xi\frac{\partial K}{\partial \bn}=0,\quad \bx\in\partial\Omega.
	\end{cases}
\end{align}
Let $\lambda=0$ is an eigenvalue for any $\xi$, then $\lambda_1=0$ must be the principal eigenvalue of \eqref{ef3}. Similarly $ {u_1}^{\ast\ast}$ satisfies 
\begin{align}
	\label{ef4}
	\begin{cases}
		&d_1\Delta\xi-\beta_1\nabla\cdot(\xi\nabla K)+\mathcal{G}(\bx,u_1^{\ast\ast})=\lambda\xi,\quad \bx\in \Omega,\\
		&d_1\frac{\partial \xi}{\partial \bn}-\beta_1\xi\frac{\partial K}{\partial \bn}=0,\quad \bx\in\partial\Omega,
	\end{cases}
\end{align}with $\lambda=0$ for any $\xi$, and therefore, $\lambda_1=0$ is also valid for  \eqref{ef4} as a principal eigenvalue. However, since $ \mathcal{G}(\bx,u_1) $ is strictly
decreasing in $u_1$ and $ {u_1}^{\ast\ast}> {u_1}^\ast $ on at least part of $\Omega$, the principal eigenvalue in \eqref{ef4}
must be less than the principal eigenvalue in \eqref{ef3}. Thus, we cannot have $\lambda=0$ in
both \eqref{ef3} and \eqref{ef4}, and therefore, \eqref{harvest_system_decoupled} cannot have any equilibrium other than the minimal
equilibrium ${u_1}^\ast$, which completes the proof.
\end{proof}
\begin{Lemma}
\label{lm19}
Let ${u_2}^0(\bx)$ be the non-negative and non-trivial initial continuous function and $ {u_2}^0(\bx) > 0 $	in some non-empty open bounded sub-domain $ \Omega_1 \subset \Omega $. Then for any $ t > 0 $, there exists a unique solution $ u_2(\bx, t) $ of the problem \eqref{harvest_system} and it is positive.
\end{Lemma}
\begin{Lemma}\label{lm20}
A unique equilibrium solution of \eqref{harvest_system} ${u_2}^\ast(\bx)$ exists. Then for any initial solution ${u_2}^0(\bx)\ge 0,\ {u_2}^0(\bx)\not\equiv 0$, the solution ${u_2}(\bx,t)$\ of \eqref{harvest_system} satisfies the condition $$\lim_{t\rightarrow\infty}{u_2}(\bx, t)={u_2}^\ast(\bx),$$ uniformly for $\bx\in \overline{\Omega}$.
\end{Lemma}
\begin{proof}
The Lemmas \ref{lm19} and \ref{lm20} are analogous to Lemma \ref{lm17} and \ref{lm18}, so we can omit the proofs here. 
\end{proof}
\begin{Th}
Let for $\gamma_1, \gamma_2\in [0,1)$ and any ${u_1}^0(\bx)$, and ${u_2}^0(\bx)$ are continuous on $\Omega$, then the system \eqref{harvest_system} has a unique solution $ (u_1, u_2) $. Further, when both initial conditions ${u_1}^0$ and $ {u_2}^0 $ are
non-negative and non-trivial, then $ u_1(\bx, t) > 0 $ and $ u_2(\bx, t) > 0 $ for any $ t > 0 $.
\end{Th}
\begin{proof}
The system \eqref{harvest_system} has a non trivial time dependent solution according to \cite{e3}. We consider \cite{e3}
$$\rho_{u_1}=\max\left\{\sup_{\bx\in \Omega}{u_1}^0(\bx),K(\bx)\right\},~~\rho_{u_2}=\max\left\{\sup_{\bx\in \Omega}{u_2}^0(\bx),K(\bx)\right\}.$$
Suppose that 
\begin{align}
	&g_1(\bx, t, u_1, u_2)=r_1u_1{\left(1-\gamma_1-\frac{u_1+u_2}{K}\right)},\\
	& g_2(\bx, t, u_1,u_2)=r_2{u_2}{\left(1-\gamma_2-\frac{u_1+u_2}{K}\right)},
\end{align}
and define
$$S^\ast:= \{(\check{u_1},\check{u_2})\in C([0,\infty)\times \Omega):0\le \check{u_1}\le \rho_{u_1},0\le \check{u_2}\le \rho_{u_2}\}.$$
Now, the functions $g_1$\ and\ $g_2$\ are quasi-monotone non-increasing Lipshitz functions in $S^\ast$, then the following conditions are satisfied according to \cite{book} 
\begin{align}
	&g_1(\bx,t,\rho_{u_1},0)\le 0\le g_1(\bx,t,0,\rho_{u_2}),\\
	&g_2(\bx,t,0,\rho_{u_2})\le 0\le g_2(\bx,t,\rho_{u_1},0).
\end{align}
Then, for the initial solution ${u_1}^0(\bx),{u_2}^0(\bx)$ and the class of continuous functions $C([0,\infty)\times \Omega)$ on $[0,\infty)\times \Omega$\ such that
$$({u_1}^0,{u_2}^0)\in S^\ast.$$
Then, a unique solution $(u_1(\bx, t),u_2(\bx,t))$\ of \eqref{harvest_system} exists in $S^\ast$ for all $(\bx,t)\in [0,\infty)\times \Omega$. Therefore, $ (u_1(\bx, t), u_2(\bx, t))$ is unique and positive solution.
\end{proof}
\section{Fully Discrete Schemes} \label{fully-discrete-schme}
For the purpose of approximating solutions of \eqref{RDE1}, we propose and analyze two fully discrete, decoupled, and linearized time-stepping algorithms. Decoupled and linearized algorithms are computationally efficient and popular and have been implemented in many areas of scientific fields including in population dynamics \cite{adan2023interplay,  mohebujjaman2024decoupled, wong2009analysis}, Navier-Stokes \cite{linke2017connection} and magnetohydrodynamcs  \cite{AKMR15,HMR17, mohebujjaman2022efficient} problems. The first-order backward-Euler formula is used in Algorithm \ref{Algn1} to approximate the temporal derivative, and the immediately preceding time-step solution is used to linearize the non-linear component. In Algorithm \ref{Algn2}, we provide the Decoupled Backward Difference Formula 2 (DBDF-2) scheme, which consists of a second-order accurate time derivative approximation formula that linearizes the non-linear term by approximating the unknown solution at the previous time-step.  For the simplicity of our analysis, we define for $i=1,2,\cdots,N$
\begin{align}
\alpha_i:=d_i-C\beta_i\|K\|_{\infty,2}-C\|r_i\|_{\infty,\infty}\left(|1-\gamma_i|+\frac{1}{K_{\min}}\right).\label{alpha-def}
\end{align}
\begin{algorithm}	
\caption{DBE scheme}\label{Algn1}
Given time-step $\Delta t>0$, end time $T>0$, for $i=1,2,\cdots,N$, initial conditions $u_i^0\in L^2(\Omega)^d$, $f_i\in L^2\left(0,T;H^{-1}(\Omega)^d\right)$, $K_{\min}>0$, $\alpha_i>0$, and $r_i\in L^\infty(0,T;L^\infty(\Omega)^d)$. Set $M=T/\Delta t$ and for $n=0,1,\cdots\hspace{-0.35mm},M-1$, compute:
Find $u_{i,h}^{n+1}\in X_h$ satisfying,  $\forall v_{h}\in X_h$:
\begin{align}
	&\left(\frac{u_{i,h}^{n+1}-u_{i,h}^{n}}{\Delta  t},v_{h}\right)+d_i\left(\nabla u_{i,h}^{n+1},\nabla v_{h}\right)-\beta_i\left(u_{i,h}^{n+1}\nabla K(t^{n+1}),\nabla v_{h}\right)\nonumber\\
	&=(1-\gamma_i)\left(r_i(t^{n+1})u_{i,h}^{n+1},v_{h}\right)-\left(\frac{r_i(t^{n+1}) u_{i,h}^{n+1}}{K(t^{n+1})}\sum\limits_{j=1}^Nu_{j,h}^n,v_{h}\right)+\left(f_i(t^{n+1}),v_{h}\right).\label{disc-weak-form}
\end{align}
\end{algorithm}
\begin{algorithm}
\caption{DBDF-2 scheme}\label{Algn2}
Given time-step $\Delta t>0$, end time $T>0$, for $i=1,2,\cdots,N$, initial conditions $u_i^0,\;u_i^1\in L^2(\Omega)^d$, $f_i\in L^2\left(0,T;H^{-1}(\Omega)^d\right)$, $\alpha_i>0$, $K_{\min}>0$, and $r_i\in L^\infty(0,T;L^\infty(\Omega)^d)$. Set $M=T/\Delta t$ and for $n=1,\cdots\hspace{-0.35mm},M-1$, compute:
Find $u_{i,h}^{n+1}\in X_h$ satisfying,  $\forall v_{h}\in X_h$:
\begin{align}
	& \Bigg(\frac{3u_{i,h}^{n+1}-4u_{i,h}^{n}+u_{h}^{n-1}}{2\Delta  t},v_{h}\Bigg)+d_i\left(\nabla u_{i,h}^{n+1},\nabla v_{h}\right)-\beta_i\left(u_{i,h}^{n+1}\nabla K(t^{n+1}),\nabla v_{h}\right)\nonumber\\
	&=(1-\gamma_i)\left(r_i(t^{n+1})u_{i,h}^{n+1},v_{h}\right)-\left(\frac{r_i(t^{n+1}) u_{i,h}^{n+1}}{K(t^{n+1})}\sum\limits_{j=1}^N(2u_{j,h}^n-u_{j,h}^{n-1}),v_{h}\right)+\left(f_i(t^{n+1}),v_{h}\right).\label{disc-weak-form2}
\end{align}
\end{algorithm}
%
\subsection{Stability Analysis}\label{stability-analysis}

The stability theorems and well-posedness of DBE and DBDF-2 schemes are investigated in this section. 

\begin{Lemma}
There exists a constant $C_*>0$ such that $\|u_{i,h}^n\|_{\infty}\le C_*$ for $i=1,2,\cdots,N$, where $u_{i,h}^n$ is a solution of the DBE/DBDF-2 scheme. \label{assumption-1}
\end{Lemma}
\begin{proof}
The proof is straightforward and similar as in the Lemma 3.5 in \cite{mohebujjaman2024decoupled}.
\end{proof}

\begin{Th}(Stability Analysis of DBE)\label{stability-theorm}
For $i=1,2,\cdots,N$, 
assume $u_{i,h}^0\in L^2(\Omega)^d$, $f_i\in L^2\left(0,T;H^{-1}(\Omega)^d\right)$, $r_i\in L^\infty(0,T;L^\infty(\Omega)^d)$, $K\in L^\infty(0,T;H^2(\Omega)^d),$ $K_{\min}>0$ and using the Lemma \ref{assumption-1}, if $\alpha_i> 0$, then for any $\Delta t>0$:
\begin{align}
	\|u_{i,h}^M\|^2+\alpha_i\Delta t\sum\limits_{n=1}^M\|\nabla u_{i,h}^n\|^2\le\|u_{i,h}^0\|^2+\frac{\Delta t}{\alpha_i}\sum\limits_{n=1}^{M} \|f_i(t^{n})\|_{-1}^2.
\end{align}
\end{Th}

\begin{proof}
Taking $v_{h}=u_{i,h}^{n+1}$ in \eqref{disc-weak-form}, and using the polarization identity
$$(q-p,q)=\frac12\left(\|q-p\|^2+\|q\|^2-\|p\|^2\right),$$
we have
\begin{align}
	&\frac{1}{2\Delta t}\Big(\|u_{i,h}^{n+1}-u_{i,h}^{n}\|^2+\|u_{i,h}^{n+1}\|^2-\|u_{i,h}^{n}\|^2\Big)+d_i\|\nabla u_{i,h}^{n+1}\|^2=\beta_i\left(u_{i,h}^{n+1}\nabla K(t^{n+1}),\nabla u_{i,h}^{n+1}\right)\nonumber\\
	&+(1-\gamma_i)\left(r_i(t^{n+1})u_{i,h}^{n+1},u_{i,h}^{n+1}\right)-\left(\frac{r_i(t^{n+1}) u_{i,h}^{n+1}}{K(t^{n+1})}\sum\limits_{j=1}^Nu_{j,h}^n,u_{i,h}^{n+1}\right)+\left(f_i(t^{n+1}),u_{i,h}^{n+1}\right).\label{pol-new}
\end{align}
For the first term on the right-hand-side of \eqref{pol-new}, we use H\"{o}lder's inequality, Sobolev embedding theorem, and Poincar\'e inequalities, to find the following estimate
\begin{align}
	\left(u_{i,h}^{n+1}\nabla K(t^{n+1}),\nabla u_{i,h}^{n+1}\right)&\le\|u_{i,h}^{n+1}\|_{L^3}\|\nabla K(t^{n+1})\|_{L^6}\|\nabla u_{i,h}^{n+1}\|\nonumber\\&\le C\|u_{i,h}^{n+1}\|^{\frac12}\|\nabla u_{i,h}^{n+1}\|^{\frac12} \| K(t^{n+1})\|_{H^2(\Omega)}\|\nabla u_{i,h}^{n+1}\|\nonumber\\&\le C\|\nabla u_{i,h}^{n+1}\|^2\|K(t^{n+1})\|_{H^2(\Omega)^d}.\label{first-term-right-up}
\end{align}
With the above estimate, we apply H\"{o}lder's inequality on the second and third terms and Cauchy Schwarz's inequality on the forcing term on the right-hand-side of \eqref{pol-new}, we have
\begin{align}
	\frac{1}{2\Delta t}\Big(\|u_{i,h}^{n+1}-u_{i,h}^{n}\|^2+\|u_{i,h}^{n+1}\|^2-\|u_{i,h}^{n}\|^2\Big)+d_i\|\nabla u_{i,h}^{n+1}\|^2 \nonumber\\\le C\beta_i\|\nabla u_{i,h}^{n+1}\|^2\|K(t^{n+1})\|_{H^2(\Omega)}+|1-\gamma_i|\|r_i(t^{n+1})\|_{\infty}\|u_{i,h}^{n+1}\|^2\nonumber\\+\Big\|\frac{r_i(t^{n+1})}{K(t^{n+1})}\Big\|_\infty\sum\limits_{j=1}^N\|u_{i,h}^{n+1}\|^2\|u_{j,h}^n\|_{\infty}+\|f_i(t^{n+1})\|_{-1}\|\nabla u_{i,h}^{n+1}\|.\label{before-small-data-assumption}
\end{align}
Use of Poincar\'e inequality and the Lemma \ref{assumption-1}, gives
\begin{align}
	\frac{1}{2\Delta t}\big(\|u_{i,h}^{n+1}-u_{i,h}^{n}\|^2+\|u_{i,h}^{n+1}\|^2-\|u_{i,h}^{n}\|^2\big)+d_i\|\nabla u_{i,h}^{n+1}\|^2\nonumber\\\le C\beta_i\|\nabla u_{i,h}^{n+1}\|^2\|K\|_{L^\infty\big(0,T;H^2(\Omega)^d\big)}+ C|1-\gamma_i|\|r_i\|_{L^\infty\big(0,T;L^\infty(\Omega)^d\big)}\|\nabla u_{i,h}^{n+1}\|^2\nonumber\\+\frac{C\|r_i\|_{L^\infty\big(0,T;L^\infty(\Omega)^d\big)}}{\inf\limits_{(t,\bx)\in (0,T]\times\Omega}|K|}\|\nabla u_{i,h}^{n+1}\|^2+\|f_i(t^{n+1})\|_{-1}\|\nabla u_{i,h}^{n+1}\|.
\end{align}
Grouping terms on the left-hand-side and using \eqref{kmin-def}, and\eqref{alpha-def}, yields
\begin{align}
	\frac{1}{2\Delta t}&\left(\|u_{i,h}^{n+1}-u_{i,h}^{n}\|^2+\|u_{i,h}^{n+1}\|^2-\|u_{i,h}^{n}\|^2\right)+\alpha_i\|\nabla u_{i,h}^{n+1}\|^2 \le \|f_i(t^{n+1})\|_{-1}\|\nabla u_{i,h}^{n+1}\|.
\end{align}

Assume $\alpha_i>0$ for $i=1,2,\cdots,N$, apply Young's inequality, and rearrange, to have
\begin{align}
	\frac{1}{2\Delta t}\left(\|u_{i,h}^{n+1}-u_{i,h}^{n}\|^2+\|u_{i,h}^{n+1}\|^2-\|u_{i,h}^{n}\|^2\right)+\frac{\alpha_i}{2}\|\nabla u_{i,h}^{n+1}\|^2 \le \frac{1}{2\alpha_i}\|f_i(t^{n+1})\|_{-1}^2.
\end{align}
Now, multiplying both sides by $2\Delta t$, summing over time steps  $n=0,1,\cdots, M-1$, and dropping non-negative terms from the left-hand-side finishes the proof. 
\end{proof}
\begin{Th}(Stability analysis of DBDF-2)\label{stability-theorm2}
For $i=1,2,\cdots,N$, 
suppose $u_{i,h}^0,u_{i,h}^1\in L^2(\Omega)^d$, $f_i\in L^2\left(0,T;H^{-1}(\Omega)^d\right)$, $K\in L^\infty(0,T;H^2(\Omega)^d),$ $K_{\min}>0$, $r_i\in L^\infty(0,T;L^\infty(\Omega)^d)$, and using the result of the Lemma \ref{assumption-1}, if $\alpha_i> 0$, then for any $\Delta t>0$:
\begin{align}
	\|u_{i,h}^{M}\|^2&+\|2u_{i,h}^{M}-u_{i,h}^{M-1}\|^2+2\alpha_i\Delta t\sum\limits_{n=2}^{M}\|\nabla u_{i,h}^{n}\|^2\nonumber\\&\le\|u_{i,h}^{1}\|^2+\|2u_{i,h}^{1}-u_{i,h}^{0}\|^2
	+\frac{2\Delta t}{\alpha_i}\sum\limits_{n=2}^{M}\|f_i(t^{n})\|_{-1}^2.
\end{align}
\end{Th}
\begin{proof}
Set $v_{h}=u_{i,h}^{n+1}$ in \eqref{disc-weak-form2}, to obtain
\begin{align}\label{pol-new2}
	&  \Bigg(\frac{3u_{i,h}^{n+1}-4u_{i,h}^{n}+u_{i,h}^{n-1}}{2\Delta  t},u_{i,h}^{n+1}\Bigg)+d_i\|\nabla u_{i,h}^{n+1}\|^2\nonumber\\&=\beta_i\left(u_{i,h}^{n+1}\nabla K(t^{n+1}),\nabla u_{i,h}^{n+1}\right)+(1-\gamma_i)\left(r_i(t^{n+1})u_{i,h}^{n+1},u_{i,h}^{n+1}\right) \nonumber\\
	&-\left(\frac{r_i(t^{n+1}) u_{i,h}^{n+1}}{K(t^{n+1})}\sum\limits_{j=1}^N(2u_{j,h}^n-u_{j,h}^{n-1}),u_{i,h}^{n+1}\right)+\left(f_i(t^{n+1}),u_{i,h}^{n+1}\right).
\end{align}
Using the upper bound of the first term on the right-hand-side of \eqref{pol-new2} given in  \eqref{first-term-right-up}, and using the following identity
\begin{eqnarray}
	(3a-4b+c,a)=\frac{a^2+(2a-b)^2}{2}-\frac{b^2+(2b-c)^2}{2}+\frac{(a-2b+c)^2}{2},\label{ident}
\end{eqnarray}
we obtain
\begin{align}
	\frac{1}{4\Delta t}\Big(&\|u_{i,h}^{n+1}\|^2-\|u_{i,h}^{n}\|^2+\|2u_{i,h}^{n+1}-u_{i,h}^{n}\|^2-\|2u_{i,h}^{n}-u_{i,h}^{n-1}\|^2+\|u_{i,h}^{n+1}-2u_{i,h}^{n}+u_{i,h}^{n-1}\|^2\Big)\nonumber\\&+d_i\|\nabla u_{i,h}^{n+1}\|^2\le C\|\nabla u_{i,h}^{n+1}\|^2\|K(t^{n+1})\|_{H^2(\Omega)^d}+(1-\gamma_i)\left(r_i(t^{n+1})u_{i,h}^{n+1},u_{i,h}^{n+1}\right)\nonumber\\
	&-\left(\frac{r_i(t^{n+1}) u_{i,h}^{n+1}}{K(t^{n+1})}\sum\limits_{j=1}^N(2u_{j,h}^n-u_{j,h}^{n-1}),u_{i,h}^{n+1}\right)+
	\left(f_i(t^{n+1}),u_{i,h}^{n+1}\right).\label{BDF2-identity}
\end{align}
Applying H\"{o}lder's inequality on the second term, H\"{o}lder's and triangle inequalities on the  third term, and Cauchy Schwarz's inequality on the forcing term on the right-hand-side of \eqref{BDF2-identity}, yields
\begin{align*}
	\frac{1}{4\Delta t}\Big(\|u_{i,h}^{n+1}\|^2-\|u_{i,h}^{n}\|^2+\|2u_{i,h}^{n+1}-u_{i,h}^{n}\|^2-\|2u_{i,h}^{n}-u_{i,h}^{n-1}\|^2+\|u_{i,h}^{n+1}-2u_{i,h}^{n}+u_{i,h}^{n-1}\|^2\Big)\nonumber\\+d_i\|\nabla u_{i,h}^{n+1}\|^2\le C\|\nabla u_{i,h}^{n+1}\|^2\|K(t^{n+1})\|_{H^2(\Omega)^d}+|1-\gamma_i|\|r_i(t^{n+1})\|_{\infty}\|u_{i,h}^{n+1}\|^2\nonumber\\+\Big\|\frac{r_i(t^{n+1})}{K(t^{n+1})}\Big\|_\infty\sum\limits_{j=1}^N\|u_{i,h}^{n+1}\|^2\big(2\|u_{j,h}^n\|_{\infty}+\|u_{j,h}^{n-1}\|_{\infty}\big)+\|f_i(t^{n+1})\|_{-1}\|\nabla u_{i,h}^{n+1}\|.
\end{align*}
Using the Poincar\'e inequality and the Lemma \ref{assumption-1}, we obtain
\begin{multline}\label{eq3}
	\frac{1}{4\Delta t}\Big(\|u_{i,h}^{n+1}\|^2-\|u_{i,h}^{n}\|^2+\|2u_{i,h}^{n+1}-u_{i,h}^{n}\|^2-\|2u_{i,h}^{n}-u_{i,h}^{n-1}\|^2+\|u_{i,h}^{n+1}-2u_{i,h}^{n}+u_{i,h}^{n-1}\|^2\Big)\\+d_i\|\nabla u_{i,h}^{n+1}\|^2\le C\|\nabla u_{i,h}^{n+1}\|^2\|K\|_{L^\infty\big(0,T;H^2(\Omega)^d\big)}+ C|1-\gamma_i|\|r_i\|_{L^\infty\big(0,T;L^\infty(\Omega)^d\big)}\|\nabla u_{i,h}^{n+1}\|^2\\
	+\frac{C\|r_i\|_{L^\infty\big(0,T;L^\infty(\Omega)^d\big)}}{\inf\limits_{(t,{\bx})\in (0,T]\times\Omega}|K|}\|\nabla u_{i,h}^{n+1}\|^2+\|f_i(t^{n+1})\|_{-1}\|\nabla u_{i,h}^{n+1}\|.
\end{multline}
Grouping the terms on the left-hand-side and using \eqref{alpha-def}, and \eqref{kmin-def}, the inequality \eqref{eq3} yields
\begin{align}
	\frac{1}{4\Delta t}\Big(\|u_{i,h}^{n+1}\|^2-\|u_{i,h}^{n}\|^2&+\|2u_{i,h}^{n+1}-u_{i,h}^{n}\|^2-\|2u_{i,h}^{n}-u_{i,h}^{n-1}\|^2+\|u_{i,h}^{n+1}-2u_{i,h}^{n}+u_{i,h}^{n-1}\|^2\Big)\nonumber\\
	&+\alpha_i\|\nabla u_{i,h}^{n+1}\|^2 \le \|f_i(t^{n+1})\|_{-1}\|\nabla u_{i,h}^{n+1}\|.
\end{align}
Assume $\alpha_i>0$ for $i=1,2,\cdots,N$ and use Young's inequality, to get
\begin{align}
	\frac{1}{4\Delta t}\Big(\|u_{i,h}^{n+1}\|^2-\|u_{i,h}^{n}\|^2+\|2u_{i,h}^{n+1}-u_{i,h}^{n}\|^2-&\|2u_{i,h}^{n}-u_{i,h}^{n-1}\|^2+\|u_{i,h}^{n+1}-2u_{i,h}^{n}+u_{i,h}^{n-1}\|^2\Big)\nonumber\\
	&+\frac{\alpha_i}{2}\|\nabla u_{i,h}^{n+1}\|^2 \le \frac{1}{2\alpha_i}\|f_i(t^{n+1})\|_{-1}^2.
\end{align}
Now, dropping non-negative term from the left-hand-side, multiplying both sides by $4\Delta t$, and summing over time-steps $n=1,\cdots, M-1$, completes the proof.
\end{proof}
\subsection{Convergence Analysis}\label{Convergence-sec}
In this section, we provide finite element error analyses and apriori estimates of the errors in the DBE and DBDF-2 schemes.

\begin{Th}(Error estimation of DBE) Consider $m=\max\{2,k+1\}$, and $i=1,2,\cdots,N$, assume $u_i$ solves \eqref{RDE1} and satisfies
\begin{align*}
	u_i\in &L^\infty\big(0,T;H^m(\Omega)^d\big),\hspace{1mm}u_{i,t}\in L^\infty\big(0,T;L^2(\Omega)^d\big),\hspace{1mm}u_{i,tt}\in L^\infty\big(0,T;L^2(\Omega)^d\big),\\&r_i\in L^\infty\left(0,T;L^\infty(\Omega)^d\right), K\in L^\infty(0,T;H^2(\Omega)^d),\text{ and } K_{\min}>0,
\end{align*}

if $\alpha_i>0$ then for $\Delta t>0$ the solution $u_{i,h}$ to the Algorithm  \ref{Algn1} converges to the true solution with
\begin{align}
	\sum\limits_{i=1}^N\|u_i(T)-u_{i,h}^M\|+\sum\limits_{i=1}^N\left\{\alpha_i\Delta t\sum_{n=1}^M\|\nabla \big(u_i(t^n)-u_{i,h}^n\big)\|^2\right\}^{\frac{1}{2}}\le C \big(h^{k}+\Delta t\big).
\end{align}\label{Convergence-analysis-BE}
\end{Th}

\begin{proof}
Firstly, we construct an error equation at the time step $t^{n+1}$, the continuous variational formulations can be written as $\forall v_{h}\in X_h$
\begin{align}
	\Bigg(&\frac{u_{i}(t^{n+1})-u_{i}(t^{n})}{\Delta  t},v_{h}\Bigg)+d_i\left(\nabla u_{i}(t^{n+1}),\nabla v_{h}\right)-\beta_i\left(u_{i}(t^{n+1})\nabla K(t^{n+1}),\nabla v_{h}\right)\nonumber\\&=
	(1-\gamma_i)\left(r_i(t^{n+1})u_{i}(t^{n+1}),v_{h}\right)-\left(\frac{r_i(t^{n+1}) u_{i}(t^{n+1})}{K(t^{n+1})}\sum\limits_{j=1}^Nu_{j}(t^{n+1}),v_{h}\right)\nonumber\\
	&+\left(f_i(t^{n+1}),v_{h}\right)+\left(\frac{u_{i}(t^{n+1})-u_{i}(t^{n})}{\Delta  t}- u_{i,t}(t^{n+1}),v_{h}\right).\label{cont-weak-form}
\end{align}
Denote $e_i^n:=u_i(t^{n+1})-u_{i,h}^n$. Subtract \eqref{disc-weak-form} from \eqref{cont-weak-form} and then rearranging yields
\begin{align}
	\left(\frac{e_{i}^{n+1}-e_{i}^{n}}{\Delta  t},v_{h}\right)+d_i\left(\nabla e_{i}^{n+1},\nabla v_{h}\right)-\beta_i\left( e_{i}^{n+1}\nabla K(t^{n+1}),\nabla v_{h}\right)\nonumber\\-(1-\gamma_i)\left(r_i(t^{n+1})e_{i}^{n+1},v_{h}\right)+\sum\limits_{j=1}^N\left(\frac{r_i(t^{n+1})}{K(t^{n+1})}e_i^{n+1}u_{j,h}^n,v_{h}
	\right)=G(t,u_i,v_{h}),\label{error-equation}
\end{align}
where \begin{align*}
	G(t,u_i,v_{h}):=&\left(\frac{u_{i}(t^{n+1})-u_{i}(t^{n})}{\Delta  t}- u_{i,t}(t^{n+1}),v_{h}\right)
	+\left(\frac{r_i(t^{n+1}) u_{i}(t^{n+1})}{K(t^{n+1})}\sum\limits_{j=1}^N\big\{u_{j,h}^{n}-u_{j}(t^{n+1})\big\},v_{h}\right).
\end{align*}
We now decompose the errors as follows:
\begin{align*}
	e_{i}^n:& = u_i(t^n)-u_{i,h}^n=(u_i(t^n)-\tilde{u}_i^n)-(u_{i,h}^n-\tilde{u}_i^n):=\eta_{i}^n-\phi_{i,h}^n,
\end{align*}
where $\tilde{u}_i^n: =P_{X_h}^{L^2}(u_i(t^n))\in X_h$ is the $L^2$ projections of $u_j(t^n)$ into $X_h$. Note that $(\eta_{i}^n,v_{h})=0\hspace{2mm} \forall v_{h}\in X_h$.  Rewriting, we have for $v_{h}\in X_h$
\begin{align}
	\Bigg(&\frac{\phi_{i,h}^{n+1}-\phi_{i,h}^{n}}{\Delta  t},v_{h}\Bigg)+d_i\left(\nabla \phi_{i,h}^{n+1},\nabla v_{h}\right)-\beta_i\left( \phi_{i,h}^{n+1}\nabla K(t^{n+1}),\nabla v_{h}\right)-(1-\gamma_i)\left(r_i(t^{n+1})\phi_{i,h}^{n+1},v_{h}\right)\nonumber\\
	&+\sum\limits_{j=1}^N\left(\frac{r_i(t^{n+1})}{K(t^{n+1})}\phi_{i,h}^{n+1}u_{j,h}^n,v_{h}
	\right)=d_i\left(\nabla \eta_{i}^{n+1},\nabla v_{h}\right)-\beta_i\left(\eta_{i}^{n+1}\nabla K(t^{n+1}),\nabla v_{h}\right)\nonumber\\
	&-(1-\gamma_i)\left(r_i(t^{n+1})\eta_{i}^{n+1},v_{h}\right)+\sum\limits_{j=1}^N\left(\frac{r_i(t^{n+1})}{K(t^{n+1})}\eta_i^{n+1}u_{j,h}^n,v_{h}
	\right)-G(t,u_i,v_{h}).\label{phi-equation}
\end{align}
Choose $v_{h}=\phi_{i,h}^{n+1}$, and use the polarization identity in \eqref{phi-equation}, to obtain
\begin{align}\label{all-phi}
	&	\frac{1}{2\Delta t}\left(\|\phi_{i,h}^{n+1}-\phi_{i,h}^{n}\|^2+\|\phi_{i,h}^{n+1}\|^2-\|\phi_{i,h}^{n}\|^2\right)+d_i\|\nabla\phi_{i,h}^{n+1}\|^2-\beta_i\left( \phi_{i,h}^{n+1}\nabla K(t^{n+1}),\nabla \phi_{i,h}^{n+1}\right)\nonumber\\
	&~~~~~~	-(1-\gamma_i)\left(r_i(t^{n+1})\phi_{i,h}^{n+1},\phi_{i,h}^{n+1}\right)	+\sum\limits_{j=1}^N\left( \frac{r_i(t^{n+1})}{K(t^{n+1})}\phi_{i,h}^{n+1}u_{j,h}^n,\phi_{i,h}^{n+1}\right)\nonumber\\
	&\le d_i\left(\nabla \eta_{i}^{n+1},\nabla \phi_{i,h}^{n+1}\right)	-\beta_i\left( \eta_{i,h}^{n+1}\nabla K(t^{n+1}),\nabla \phi_{i,h}^{n+1}\right)-(1-\gamma_i)\left(r_i(t^{n+1})\eta_{i}^{n+1},\phi_{i,h}^{n+1}\right)\nonumber\\
	&+\sum\limits_{j=1}^N\left(\frac{r_i(t^{n+1})}{K(t^{n+1})}\eta_{i}^{n+1} u_{j,h}^n,\phi_{i,h}^{n+1}\right)-G\left(t,u_i,\phi_{i,h}^{n+1}\right).
\end{align}
Now, we find the upper-bounds of terms in the above equation. Similar as \eqref{first-term-right-up}, we use H\"{o}lder's inequality, Sobolev embedding theorem, and Poincar\'e inequalities, to find 
\begin{align}
	\left(\phi_{i,h}^{n+1}\nabla K(t^{n+1}),\nabla\phi_{i,h}^{n+1}\right)&\le C\|\nabla \phi_{i,h}^{n+1}\|^2\|K(t^{n+1})\|_{L^\infty\big(0,T;H^2(\Omega)^d\big)},\label{bound-phi-1}\\
	\left(\eta_{i}^{n+1}\nabla K(t^{n+1}),\nabla \phi_{i,h}^{n+1}\right)&\le C\|\nabla\eta_{i}^{n+1}\|\|K(t^{n+1})\|_{L^\infty\big(0,T;H^2(\Omega)^d\big)}\|\nabla \phi_{i,h}^{n+1}\|.
\end{align}
Apply Young's inequality, to have
\begin{align}
	\beta_i\left(\eta_{i}^{n+1}\nabla K(t^{n+1}),\nabla \phi_{i,h}^{n+1}\right)\le\frac{\alpha_i}{12}\|\nabla \phi_{i,h}^{n+1}\|^2+\frac{C\beta_i^2}{\alpha_i}\|\nabla\eta_{i}^{n+1}\|^2\|K\|_{\infty,2}^2.
\end{align}
Using H\"older's, and Poincar\'e inequalities, we have
\begin{align}
	(1-\gamma_i)\left(r_i(t^{n+1})\phi_{i,h}^{n+1},\phi_{i,h}^{n+1}\right)&\le C |1-\gamma_i|\|r_i\|_{\infty,\infty}\|\nabla \phi_{i,h}^{n+1}\|^2.\label{bound-gamma-1}
\end{align}
With the assumption $\alpha_i>0$, use Cauchy-Schwarz, and Young's inequalities, to obtain
\begin{align}
	d_i\left(\nabla \eta_{i}^{n+1},\nabla \phi_{i,h}^{n+1}\right)\le d_i\|\nabla \eta_{i}^{n+1}\|\|\nabla \phi_{i,h}^{n+1}\|\le\frac{\alpha_i}{12}\|\nabla \phi_{i,h}^{n+1}\|^2+\frac{3d_i^2}{\alpha_i}\|\nabla \eta_{i}^{n+1}\|^2.\label{bound-d-1}
\end{align}
Next, using triangle, H\"older's, and Poincar\'e inequalities together with the Lemma \ref{assumption-1}, we have
\begin{align*}
	-\sum\limits_{j=1}^N\left( \frac{r_i(t^{n+1}) }{K(t^{n+1})}\phi_{i,h}^{n+1}u_{j,h}^n,\phi_{i,h}^{n+1}\right)&\le\sum\limits_{j=1}^N\frac{\|r_i(t^{n+1})\|_{\infty}}{\inf\limits_\Omega\|K(t^{n+1})\|}\Big|\left(\phi_{i,h}^{n+1}u_{j,h}^n,\phi_{i,h}^{n+1}\right)\Big|\nonumber\\&\le\sum\limits_{j=1}^N\frac{\|r_i\|_{\infty,\infty}}{K_{\min}}\|u_{j,h}^{n}\|_\infty\|\phi_{i,h}^{n+1}\|^2\nonumber\\&\le\frac{C\|r_i\|_{\infty,\infty}}{K_{\min}}
	\|\nabla\phi_{i,h}^{n+1}\|^2.
\end{align*} 
Using H\"older's, and Poincar\'e inequalities, Lemma \ref{assumption-1}, and Young's inequality, we have
\begin{align}
	(1-\gamma_i)\left(r_i(t^{n+1})\eta_{i}^{n+1},\phi_{i,h}^{n+1}\right)&\le C|1-\gamma_i|\|r_i\|_{\infty,\infty}\|\eta_{i}^{n+1}\|\|\nabla\phi_{i,h}^{n+1}\|\nonumber\\&\le\frac{\alpha_i}{12}\|\nabla\phi_{i,h}^{n+1}\|^2+\frac{C(1-\gamma_i)^2\|r_i\|_{\infty,\infty}^2}{\alpha_i}\|\eta_{i}^{n+1}\|^2,\label{bound-eta-1}\\
	\sum\limits_{j=1}^N\left(\frac{r_i(t^{n+1})}{K(t^{n+1})}\eta_{i}^{n+1} u_{j,h}^n,\phi_{i,h}^{n+1}\right)&\le\sum\limits_{j=1}^N\frac{\|r_i(t^{n+1})\|_{\infty}}{\inf\limits_\Omega\|K(t^{n+1})\|}\Big|\left(\eta_{i}^{n+1} u_{j,h}^n,\phi_{i,h}^{n+1}\right)\Big|\nonumber\\&\le \sum\limits_{j=1}^N\frac{\|r_i\|_{\infty,\infty}}{K_{\min}}\|\eta_{i}^{n+1}\|\| u_{j,h}^n\|_\infty\|\|\phi_{i,h}^{n+1}\|\nonumber\\&\le \frac{C\|r_i\|_{\infty,\infty}}{K_{\min}}\|\eta_{i}^{n+1}\|\|\nabla\phi_{i,h}^{n+1}\|\nonumber\\&\le \frac{\alpha_i}{12}\|\nabla\phi_{i,h}^{n+1}\|^2+\frac{C\|r_i\|_{\infty,\infty}^2}{\alpha_iK_{\min}^2}\|\eta_{i}^{n+1}\|^2.
\end{align}
Now, we want to find the upper-bound of \begin{align}    G(t,u_i,\phi_{i,h}^{n+1})=\left(\frac{u_{i}(t^{n+1})-u_{i}(t^{n})}{\Delta  t}- u_{i,t}(t^{n+1}),\phi_{i,h}^{n+1}\right)\nonumber\\+\left(\frac{r_i(t^{n+1}) u_{i}(t^{n+1})}{K(t^{n+1})}\sum\limits_{j=1}^N\big\{u_{j,h}^{n}-u_{j}(t^{n+1})\big\},\phi_{i,h}^{n+1}\right).\label{G-phi}
\end{align}
We can find the upper-bound of the last term on the right-hand-side of \eqref{G-phi} using H\"older's, triangle, Poincar\'e, Young's inequalities, and  Taylor's series expansion, together with the regularity assumption as
\begin{align*}
	\Bigg(\frac{r_i(t^{n+1}) u_{i}(t^{n+1})}{K(t^{n+1})}&\sum\limits_{j=1}^N\big\{u_{j,h}^{n}-u_{j}(t^{n+1})\big\},\phi_{i,h}^{n+1}\Bigg)\\&\le\frac{C\|r_i\|_{\infty,\infty}}{K_{\min}}\|u_{i}(t^{n+1})\|_\infty\sum\limits_{j=1}^N\|u_{j,h}^{n}-u_{j}(t^{n+1})\|\|\nabla\phi_{i,h}^{n+1}\|\\&\le\sum\limits_{j=1}^N\frac{C\|r_i\|_{\infty,\infty}}{K_{\min}}\|u_{j,h}^{n}-u_{j}(t^{n+1})\|\|\nabla\phi_{i,h}^{n+1}\|\\&\le\sum\limits_{j=1}^N\left(\frac{\alpha_i}{12N}\|\nabla\phi_{i,h}^{n+1}\|^2+\frac{C\|r_i\|_{\infty,\infty}^2}{\alpha_i K_{\min}^2}\|u_{j,h}^{n}-u_{j}(t^{n+1})\|^2\right)\\&\le \frac{\alpha_i}{12}\|\nabla\phi_{i,h}^{n+1}\|^2+\frac{C\|r_i\|_{\infty,\infty}^2}{\alpha_i K_{\min}^2}\sum_{j=1}^N\big(\|u_{j,h}^{n}-u_j(t^n)\|^2+\|u_j(t^n)-u_j(t^{n+1})\|^2\big)\\&\le\frac{\alpha_i}{12}\|\nabla\phi_{i,h}^{n+1}\|^2+\frac{C\|r_i\|_{\infty,\infty}^2}{\alpha_i K_{\min}^2}\sum_{j=1}^N\big(2\|\eta_j^n\|^2+2\|\phi_{j,h}^n\|^2+(\Delta t)^2\|u_{j,t}(s^*)\|^2\big)\\&\le\frac{\alpha_i}{12}\|\nabla\phi_{i,h}^{n+1}\|^2+\frac{C\|r_i\|_{\infty,\infty}^2}{\alpha_i K_{\min}^2}\big(h^{2k+2}+(\Delta t)^2\big)+\frac{C\|r_i\|_{\infty,\infty}^2}{\alpha_i K_{\min}^2}\sum_{j=1}^N\|\phi_{j,h}^n\|^2,
\end{align*}
with $s^*\in[t^n,t^{n+1}]$.
For some $t^*\in[t^n,t^{n+1}]$, we use Taylor's series expansion, Cauchy-Schwarz, Poincar\'e, and Young's inequalities to obtain the following bound for the first term on the right-hand-side of \eqref{G-phi}
\begin{align*}
	\left(\frac{u_{i}(t^{n+1})-u_{i}(t^{n})}{\Delta  t}- u_{i,t}(t^{n+1}),\phi_{i,h}^{n+1}\right)&=\frac{\Delta t}{2}\left(u_{i,tt}(t^*),\phi_{i,h}^{n+1}\right)\\&\le\frac{C\Delta t}{2}\|u_{i,tt}(t^*)\|\|\nabla\phi_{i,h}^{n+1}\|\\&\le\frac{\alpha_i}{12}\|\nabla\phi_{i,h}^{n+1}\|^2+\frac{C(\Delta t)^2}{\alpha_i}\|u_{i,tt}(t^*)\|^2.
\end{align*}
Thus, we have
\begin{align*}   G(t,u_i,\phi_{i,h}^{n+1})\le\frac{\alpha_i}{6}\|\nabla\phi_{i,h}^{n+1}\|^2+C\big(h^{2k+2}+(\Delta t)^2\big)+C\sum_{j=1}^N\|\phi_{j,h}^n\|^2.
\end{align*}
Now, using the above bounds, we can rewrite \eqref{all-phi} as\begin{align}
	&\frac{1}{2\Delta t}\left(\|\phi_{i,h}^{n+1}-\phi_{i,h}^{n}\|^2+\|\phi_{i,h}^{n+1}\|^2-\|\phi_{i,h}^{n}\|^2\right)+\frac{\alpha_i}{2}\|\nabla\phi_{i,h}^{n+1}\|^2\nonumber\\&\le \frac{C\beta_i^2}{\alpha_i}\|\nabla\eta_{i}^{n+1}\|^2\|K\|_{\infty,2}^2+\frac{3d_i^2}{\alpha_i}\|\nabla \eta_{i}^{n+1}\|^2+C\sum_{j=1}^N\|\phi_{j,h}^n\|^2\nonumber\\&+\frac{C(1-\gamma_i)^2\|r_i(t^{n+1})\|_\infty^2}{\alpha_i}\|\eta_{i}^{n+1}\|^2+\frac{C\|r_i\|_{\infty,\infty}^2}{\alpha_i K_{\min}^2}\|\eta_{i}^{n+1}\|^2+C\big(h^{2k+2}+(\Delta t)^2\big).
\end{align}
Using the regularity assumption again, we obtain
\begin{align}
	\frac{1}{2\Delta t}\left(\|\phi_{i,h}^{n+1}-\phi_{i,h}^{n}\|^2+\|\phi_{i,h}^{n+1}\|^2-\|\phi_{i,h}^{n}\|^2\right)+\frac{\alpha_i}{2}\|\nabla\phi_{i,h}^{n+1}\|^2\le C\big(h^{2k}+(\Delta t)^2\big)+C\sum_{j=1}^N\|\phi_{j,h}^n\|^2.
\end{align}
Dropping non-negative term from left-hand-side, multiplying both sides by $2\Delta t$, use $\|\phi_{i,h}^0\|=0$, $\Delta{t}M=T$, and summing over time-steps $n=0,1,\cdots,M-1$, to find
\begin{align}
	\|\phi_{i,h}^{M}\|^2+\alpha_i\Delta t\sum\limits_{n=1}^{M}\|\nabla\phi_{i,h}^{n}\|^2\le \Delta t\sum\limits_{n=1}^{M-1}C\left(\sum_{j=1}^N\|\phi_{j,h}^n\|^2\right)+C\big(h^{2k}+(\Delta t)^2\big).
\end{align}
Sum over $i=1,2,\cdots,N$, we have
\begin{align}
\sum\limits_{i=1}^N\|\phi_{i,h}^{M}\|^2+\Delta t\sum\limits_{n=1}^{M}\left(\sum\limits_{i=1}^N\alpha_i\|\nabla\phi_{i,h}^{n}\|^2\right)\le \Delta t\sum\limits_{n=1}^{M-1}C\left(\sum_{i=1}^N\|\phi_{i,h}^n\|^2\right)+C\big(h^{2k}+(\Delta t)^2\big).
\end{align}
Applying the discrete Gr\"onwall Lemma \ref{dgl}, we get
\begin{align}
	\sum\limits_{i=1}^N\|\phi_{i,h}^{M}\|^2+\Delta t\sum\limits_{n=1}^{M}\left(\sum\limits_{i=1}^N\alpha_i\|\nabla\phi_{i,h}^{n}\|^2\right)\le C \big(h^{2k}+(\Delta t)^2\big),
\end{align}
which gives
\begin{align}
	\|\phi_{i,h}^{M}\|^2+\alpha_i\Delta t\sum\limits_{n=1}^{M}\|\nabla\phi_{i,h}^{n}\|^2\le C \big(h^{2k}+(\Delta t)^2\big)\hspace{4mm}\text{for}\hspace{1mm}i=1,2,\cdots,N.\label{bound-on-phi}
\end{align}
Use of triangle and Young's inequalities allows us to write
\begin{align}
	&\|e_i^M\|^2+\alpha_i\Delta t\sum_{n=1}^M\|\nabla e_i^n\|^2\le 2\left(\|\phi_{i,h}^{M}\|^2+\alpha_i\Delta t\sum\limits_{n=1}^{M}\|\nabla\phi_{i,h}^{n}\|^2+\|\eta_i^M\|^2+\alpha_i\Delta t\sum\limits_{n=1}^{M}\|\nabla\eta_{i}^{n}\|^2\right).
\end{align}
Using regularity assumptions and bound in \eqref{bound-on-phi}, we have, for $i=1,2,\cdots,N$
\begin{align}
	\|u_i(T)-u_{i,h}^M\|^2+\alpha_i\Delta t\sum_{n=1}^M\|\nabla \big(u_i(t^n)-u_{i,h}^n\big)\|^2\le C\big(h^{2k}+(\Delta t)^2\big).
\end{align}
Now, summing over $i=1,2,\cdots,N$ completes the proof.
\end{proof}

\begin{Th} (Error estimation of DBDF-2) For $i=1,2,\cdots,N$, assume $u_i$ solves \eqref{RDE1} and satisfies
\begin{align*}
	u_i\in &L^\infty\left(0,T;H^{k+1}(\Omega)^d\right), u_{i,tt}\in L^\infty\left(0,T;L^2(\Omega)^d\right),u_{i,ttt}\in L^\infty\left(0,T;L^2(\Omega)^d\right),\\
	&r_i\in L^\infty\left(0,T;L^\infty(\Omega)^d\right), K\in L^\infty(0,T;H^2(\Omega)^d),\text{ and } K_{\min}>0,
\end{align*}
if $\alpha_i>0$ then for $\Delta t>0$ the solution $u_{i,h}$ to the Algorithm \ref{Algn2} converges to the true solution with
\begin{align}
	\|u_i(T)-u_{i,h}^M\|+\Big\{\alpha_i\Delta t\sum\limits_{n=2}^M\|\nabla(u_i(t^n)-u_{i,h}^n)\|^2\Big\}^{1/2}\le C(h^k+\Delta t^2).
\end{align}\label{Convergence-analysis-BDF-2}
\end{Th}

\begin{proof}
At first we construct an error equation at the time step $t^{n+1}$, the continuous variational formulations can be written as $\forall v_{h}\in X_h$
\begin{align}
	\Bigg(&\frac{3u_{i}(t^{n+1})-4u_{i}(t^{n})+u_{i}(t^{n-1})}{2\Delta  t},v_{h}\Bigg)+d_i\left(\nabla u_{i}(t^{n+1}),\nabla v_{h}\right)-\beta_i\left( u_{i}(t^{n+1})\nabla K(t^{n+1}),\nabla v_{h}\right)\nonumber\\
	&=(1-\gamma_i)\left(r_i(t^{n+1})u_{i}(t^{n+1}),v_{h}\right)-\left(\frac{r_i(t^{n+1}) u_{i}(t^{n+1})}{K(t^{n+1})}\sum\limits_{j=1}^Nu_{j}(t^{n+1}),v_{h}\right)\nonumber\\
	&+\left(\frac{3u_{i}(t^{n+1})-4u_{i}(t^{n})+u_{i}(t^{n-1})}{2\Delta  t}- u_{i,t}(t^{n+1}),v_{h}\right)+\left(f_i(t^{n+1}),v_{h}\right). \label{cont-weak-form2}
\end{align}
Subtract \eqref{disc-weak-form2} from \eqref{cont-weak-form2} and then rearranging, yields
\begin{align}
	&\left(\frac{3e_{i}^{n+1}-4e_{i}^{n}+e_{i}^{n-1}}{2\Delta  t},v_{h}\right)+d_i\left(\nabla e_{i}^{n+1},\nabla v_{h}\right)-\beta_i\left(e_{i}^{n+1}\nabla K(t^{n+1}),\nabla v_{h}\right)\nonumber\\
	&-(1-\gamma_i)\left(r_i(t^{n+1})e_{i}^{n+1},v_{h}\right)+\sum\limits_{j=1}^N\left(\frac{r_i(t^{n+1})}{K(t^{n+1})}e_i^{n+1}(2u_{j,h}^n-u_{j,h}^{n-1}),v_{h}
	\right)\nonumber\\
	&+\sum\limits_{j=1}^N\left(\frac{r_i(t^{n+1})}{K(t^{n+1})}u_i(t^{n+1})(2e_{j}^n-e_{j}^{n-1}),v_{h}
	\right)=G(t,u_i,v_{h}),
	\label{error-equation2}
\end{align}
where\begin{align*}
	G(t,u_i,v_{h}):=\left(\frac{3u_{i}(t^{n+1})-4u_{i}(t^{n})+u_{i}(t^{n-1})}{2\Delta  t}- u_{i,t}(t^{n+1}),v_{h}\right)\\-\left(\frac{r_i(t^{n+1}) u_{i}(t^{n+1})}{K(t^{n+1})}\sum\limits_{j=1}^N\big\{u_{j}(t^{n+1})-2u_{j}(t^{n})+u_j(t^{n-1})\big\},v_{h}\right).
\end{align*}
Now we decompose the errors and rewrite, then for $\forall v_{h}\in X_h$
\begin{align}
	\Bigg(&\frac{3\phi_{i,h}^{n+1}-4\phi_{i,h}^{n}+\phi_{i,h}^{n-1}}{2\Delta  t},v_{h}\Bigg)+d_i\left(\nabla \phi_{i,h}^{n+1},\nabla v_{h}\right)-\beta_i\left( \phi_{i,h}^{n+1}\nabla K(t^{n+1}),\nabla v_{h}\right)\nonumber\\
	&-(1-\gamma_i)\left(r_i(t^{n+1})\phi_{i,h}^{n+1},v_{h}\right)
	+\sum\limits_{j=1}^N\left(\frac{r_i(t^{n+1})}{K(t^{n+1})}\phi_{i,h}^{n+1}(2u_{j,h}^n-u_{j,h}^{n-1}),v_{h}
	\right)\nonumber\\
	&+\sum\limits_{j=1}^N\left(\frac{r_i(t^{n+1})}{K(t^{n+1})}u_i(t^{n+1})(2\phi_{j,h}^{n}-\phi_{j,h}^{n-1}),v_{h}
	\right)=d_i\left(\nabla \eta_{i}^{n+1},\nabla v_{h}\right)-\beta_i\left(\eta_{i}^{n+1}\nabla K(t^{n+1}),\nabla v_{h}\right)\nonumber\\
	&-(1-\gamma_i)\left(r_i(t^{n+1})\eta_{i}^{n+1},v_{h}\right)+\sum\limits_{j=1}^N\left(\frac{r_i(t^{n+1})}{K(t^{n+1})}\eta_i^{n+1}(2u_{j,h}^n-u_{j,h}^{n-1}),v_{h}
	\right)\nonumber\\ 
	&+\sum\limits_{j=1}^N\left(\frac{r_i(t^{n+1})}{K(t^{n+1})}u_i(t^{n+1})(2\eta_{j}^{n}-\eta_{j}^{n-1}),v_{h}
	\right)-G(t,u_i,v_{h}).\label{phi-equation2}
\end{align}
Choose $v_{h}=\phi_{i,h}^{n+1}$, and use the polarization identity in \eqref{ident}, to obtain
\begin{align}
	&\frac{1}{4\Delta t}\Big(\|\phi_{i,h}^{n+1}\|^2-\|\phi_{i,h}^{n}\|^2+\|2\phi_{i,h}^{n+1}-\phi_{i,h}^{n}\|^2-\|2\phi_{i,h}^{n}-\phi_{i,h}^{n-1}\|^2+\|\phi_{i,h}^{n+1}-2\phi_{i,h}^{n}+\phi_{i,h}^{n-1}\|^2\Big)\nonumber\\
	&+d_i\|\nabla \phi_{i,h}^{n+1}\|^2-\beta_i\left(\phi_{i,h}^{n+1}\nabla K(t^{n+1}),\nabla \phi_{i,h}^{n+1}\right)-(1-\gamma_i)\left(r_i(t^{n+1})\phi_{i,h}^{n+1},\phi_{i,h}^{n+1}\right)\nonumber\\
	&+\sum\limits_{j=1}^N\left(\frac{r_i(t^{n+1})}{K(t^{n+1})}\phi_{i,h}^{n+1}(2u_{j,h}^n-u_{j,h}^{n-1}),\phi_{i,h}^{n+1}
	\right)+\sum\limits_{j=1}^N\left(\frac{r_i(t^{n+1})}{K(t^{n+1})}u_i(t^{n+1})(2\phi_{j,h}^{n}-\phi_{j,h}^{n-1}),\phi_{i,h}^{n+1}
	\right)\nonumber\\
	& \le d_i\left(\nabla \eta_{i}^{n+1},\nabla \phi_{i,h}^{n+1}\right)-\beta_i\left(\eta_{i}^{n+1}\nabla K(t^{n+1}),\nabla \phi_{i,h}^{n+1}\right)-(1-\gamma_i)\left(r_i(t^{n+1})\eta_{i}^{n+1},\phi_{i,h}^{n+1}\right)\nonumber\\
	&+\sum\limits_{j=1}^N\left(\frac{r_i(t^{n+1})}{K(t^{n+1})}\eta_i^{n+1}(2u_{j,h}^n-u_{j,h}^{n-1}),\phi_{i,h}^{n+1}\right)+\sum\limits_{j=1}^N\left(\frac{r_i(t^{n+1})}{K(t^{n+1})}u_i(t^{n+1})(2\eta_{j}^{n}-\eta_{j}^{n-1}),\phi_{i,h}^{n+1}
	\right)\nonumber\\
	&-G\left(t,u_i,\phi_{i,h}^{n+1}\right).\label{all-phi-2}
\end{align}
Now, we find the upper-bounds of terms in the above equation. 
Using triangle, H\"older's, and Poincar\'e inequalities together with the Lemma \ref{assumption-1}, we have
\begin{align*}
	-\sum\limits_{j=1}^N\left(\frac{r_i(t^{n+1})}{K(t^{n+1})}\phi_{i,h}^{n+1}(2u_{j,h}^n-u_{j,h}^{n-1}),\phi_{i,h}^{n+1}
	\right)&\le\sum\limits_{j=1}^N\Big\|\frac{r_i(t^{n+1})}{K(t^{n+1})}\Big\|_{\infty}\|2u_{j,h}^n-u_{j,h}^{n-1}\|_{\infty}\|\phi_{i,h}^{n+1}\|^2\\
	&\le\frac{C\|r_i\|_{\infty,\infty}}{K_{\min}}
	\|\nabla\phi_{i,h}^{n+1}\|^2.
\end{align*}
Using H\"older's, and Poincar\'e inequalities, Lemma \ref{assumption-1}, and Young's inequality, we have
\begin{align*}
\sum\limits_{j=1}^N\left(\frac{r_i(t^{n+1})}{K(t^{n+1})}\eta_{i}^{n+1} \big(2u_{j,h}^n-u_{j,h}^{n-1}\big),\phi_{i,h}^{n+1}\right)&\le\sum\limits_{j=1}^N\frac{\|r_i\|_{\infty,\infty}}{K_{\min}}\|\eta_{i}^{n+1}\|\|2u_{j,h}^n-u_{j,h}^{n-1}\|_{\infty}\|\phi_{i,h}^{n+1}\|\\&\le \sum\limits_{j=1}^N\frac{\|r_i\|_{\infty,\infty}}{K_{\min}}\|\eta_{i}^{n+1}\|\left(2\| u_{j,h}^n\|_\infty+\| u_{j,h}^{n-1}\|_\infty\right)\|\nabla\phi_{i,h}^{n+1}\|\\&\le \frac{C\|r_i\|_{\infty,\infty}}{K_{\min}}\|\eta_{i}^{n+1}\|\|\nabla\phi_{i,h}^{n+1}\|\\&\le \frac{\alpha_i}{12}\|\nabla\phi_{i,h}^{n+1}\|^2+\frac{C\|r_i\|^2_{\infty,\infty}}{\alpha_iK_{\min}^2}\|\eta_{i}^{n+1}\|^2.
\end{align*}
We have used H\"older's inequality, Sobolev embedding theorem,  Poincar\'e inequality, regularity assumption, and  Young's inequality, we have
\begin{align*}
	\sum\limits_{j=1}^N\bigg(&\frac{r_i(t^{n+1})}{K(t^{n+1})}u_i(t^{n+1})(2\eta_{j}^{n}-\eta_{j}^{n-1}),\phi_{i,h}^{n+1}\bigg)\\
	&\le\sum\limits_{j=1}^N\frac{\|r_i\|_{\infty,\infty}}{K_{\min}}\|u_i(t^{n+1})\|_{L^6}\|2\eta_{j}^{n}-\eta_{j}^{n-1}\|\|\phi_{i,h}^{n+1}\|_{L^3}\\
	&\le\sum\limits_{j=1}^N\frac{\|r_i\|_{\infty,\infty}}{K_{\min}}\|u_i(t^{n+1})\|_{H^1}\|2\eta_{j}^{n}-\eta_{j}^{n-1}\|\|\phi_{i,h}^{n+1}\|^{1/2}\|\nabla\phi_{i,h}^{n+1}\|^{1/2}\\
	&\le\sum\limits_{j=1}^N\frac{C\|r_i\|_{\infty,\infty}}{K_{\min}}\|u_i\|_{L^\infty\big(0,T;H^1(\Omega)^d\big)}\|2\eta_{j}^{n}-\eta_{j}^{n-1}\|\|\nabla\phi_{i,h}^{n+1}\|\\
	&\le\frac{\alpha_i}{12}\|\nabla\phi_{i,h}^{n+1}\|^2+\frac{C\|r_i\|_{\infty,\infty}^2}{K_{\min}^2}\sum\limits_{j=1}^N\|2\eta_{j}^{n}-\eta_{j}^{n-1}\|^2.
\end{align*}
Using Taylor’s series, Cauchy-Schwarz and Young’s inequalities the last term is assessed as
\begin{align*}
	-G&\left(t,u_i,\phi_{i,h}^{n+1}\right)\le \Delta t^2\frac{C\|r_i\|_{\infty,\infty}}{K_{\min}}\|u_i\|_{L^\infty\big(0,T;H^1(\Omega)^d\big)}\sum\limits_{j=1}^N\|u_{j,tt}\|_{L^\infty\big(0,T;L^2(\Omega)^d\big)}\|\nabla\phi_{i,h}^{n+1}\|\\&+C\Delta t^2\|u_{i,ttt}\|_{L^\infty\big(0,T;L^2(\Omega)^d\big)}\|\nabla\phi_{i,h}^{n+1}\|\le \frac{\alpha_i}{12}\|\nabla\phi_{i,h}^{n+1}\|^2+\frac{C\Delta t^4\|r_i\|_{\infty,\infty}^2}{\alpha_iK_{\min}^2}.
\end{align*}
	Using the estimates in \eqref{bound-phi-1}-\eqref{bound-eta-1} together with the above estimates in \eqref{all-phi-2} and reducing, yields
	\begin{align}
		&\frac{1}{4\Delta t}\Big(\|\phi_{i,h}^{n+1}\|^2-\|\phi_{i,h}^{n}\|^2+\|2\phi_{i,h}^{n+1}-\phi_{i,h}^{n}\|^2-\|2\phi_{i,h}^{n}-\phi_{i,h}^{n-1}\|^2+\|\phi_{i,h}^{n+1}-2\phi_{i,h}^{n}+\phi_{i,h}^{n-1}\|^2\Big)\nonumber\\
		&+\frac{\alpha_i}{2}\|\nabla\phi_{i,h}^{n+1}\|^2\le \frac{C\|r_i\|_{\infty,\infty}^2}{\alpha_iK_{\min}^2}\sum\limits_{j=1}^N\|2\phi_{j,h}^{n}-\phi_{j,h}^{n-1}\|^2+\frac{3d_i^2}{\alpha_i}\|\nabla \eta_{i}^{n+1}\|^2-\frac{C\beta_i^2}{\alpha_i}\|\nabla\eta_{i}^{n+1}\|^2\|K\|_{\infty,2}^2\nonumber\\
		&+\frac{C(1-\gamma_i)^2\|r_i\|^2_{\infty,\infty}}{\alpha_i}\|\eta_{i}^{n+1}\|^2+\frac{C\|r_i\|^2_{\infty,\infty}}{\alpha_iK_{\min}^2}\|\eta_{i}^{n+1}\|^2+\frac{C\|r_i\|_{\infty,\infty}^2}{K_{\min}^2}\sum\limits_{j=1}^N\|2\eta_{j}^{n}-\eta_{j}^{n-1}\|^2\nonumber\\
		&+\frac{C\Delta t^4\|r_i\|_{\infty,\infty}^2}{\alpha_iK_{\min}^2}\le C\sum\limits_{j=1}^N\|2\phi_{j,h}^{n}-\phi_{j,h}^{n-1}\|^2+Ch^{2k}+Ch^{2k+2}+C\Delta t^4.
	\end{align}
	Dropping non-negative term from left-hand-side, multiply both sides by $4\Delta t$, using $\|\phi_{i,h}^0\|=\|\phi_{i,h}^1\|=0$, and sum over the time-steps $n=1,2,\cdots,M-1$, we have
	\begin{align}
		\|\phi_{i,h}^{M}\|^2+\|2\phi_{i,h}^{M}-\phi_{i,h}^{M-1}\|^2+2\alpha_i\Delta t\sum\limits_{n=2}^M\|\nabla\phi_{i,h}^n\|^2\nonumber\\\le C\Delta t\sum\limits_{n=1}^{M-1}\sum\limits_{j=1}^N\|2\phi_{j,h}^{n}-\phi_{j,h}^{n-1}\|^2+C(h^{2k}+\Delta t^4).
	\end{align}
	Summing over $i=1,2,\cdots,N$, remove non-negative terms from left-hand-side, and reducing, gives
	\begin{align*}
		\sum\limits_{i=1}^N\|\phi_{i,h}^{M}\|^2+2\alpha_i\Delta t\sum\limits_{n=2}^M\sum\limits_{i=1}^N\|\nabla\phi_{i,h}^n\|^2\le C\Delta t\sum\limits_{n=2}^{M-1}\sum\limits_{i=1}^N\|\phi_{i,h}^n\|^2+C(h^{2k}+\Delta t^4).
	\end{align*}
	Applying the discrete Gr\"onwall Lemma \ref{dgl}, we have
	\begin{align}
		\sum\limits_{i=1}^N\|\phi_{i,h}^{M}\|^2+2\alpha_i\Delta t\sum\limits_{n=2}^M\sum\limits_{i=1}^N\|\nabla\phi_{i,h}^n\|^2\le C(h^{2k}+\Delta t^4),
	\end{align}
	for $i=1,2,\cdots,N$, which gives
	\begin{align}
		\|\phi_{i,h}^{M}\|^2+2\alpha_i\Delta t\sum\limits_{n=2}^M\|\nabla\phi_{i,h}^n\|^2\le C(h^{2k}+\Delta t^4).
	\end{align}
	The proof is completed by the use of the triangle and Young's inequalities, as well as the regularity assumption.
\end{proof}

\section{Numerical tests}\label{numerical-test}
Several numerical tests are conducted in this section to corroborate theoretical findings and provide an explanation of the effect of various parameters observed in the simulated results. We have utilised structured triangular meshes and the $\mathbb{P}_2$ element for the finite element simulation in all of the experiments. The codes were written in the finite element plateform Freefem++ \cite{hecht2012new} and direct solver UMFPACK \cite{davis2004algorithm} was used.
\subsection{Convergence rate verification} Define the global error $e_i:=u_i-u_{i,h}$ and its $L^2$-$H^1$ norm as $\|e_i\|_{2,1}:=\|e_i\|_{L^2\big(0,T;H^1(\Omega)^d\big)}$, and use the following approximation $$\|e_i\|_{2,1}\approx\sqrt{\sum_{n=1}^M\Delta t\|e_i^n\|^2_{H^1(\Omega)^d}}.$$  We have seen from the convergence analysis that the predicted error of the Algorithm \ref{Algn1} and \ref{Algn2} and  for $\mathbb{P}_2$ finite element are
\begin{align}
	\|u_i-u_{i,h}\|_{2,1}&\le C(h^2+\Delta t),\hspace{2mm}\text{and}\label{error-estimate-BE}\\
	\|u_i-u_{i,h}\|_{2,1}&\le C(h^2+\Delta t^2),\hspace{1mm}i=1,2,\cdots,N,\label{error-estimate-BDF-2}
\end{align}
respectively. To verify the above convergence rates, we consider a domain $\Omega=(0,1)\times (0,1)$, where species compete in a heterogeneous environment. We plugin the following 3-species manufactured analytical solution
\begin{align*}
u_1(t,\bx) &= \big(1.1 + \sin t\big)(2.0 + \sin(y)),\\
u_2(t,\bx) &=  \big(2.0 + \cos(t)\big)\big(1.1+\cos(x)\big),\\
u_3(t,\bx)&=\big(1.1+\sin(t)\big)\big(1.1+\cos(y)\big),
\end{align*}
in \eqref{RDE1}, together with  $$K(t,\bx)\equiv(2.1+\cos x\cos y)(1.1+\cos t),\;\;\text{and}\;\;r_i(t,\bx)\equiv(1.5+\sin x\sin y)(1.2+\sin t),$$ as the system carrying capacity, and intrinsic growth rate, respectively, to obtain $f_i(t,\bx)$, for $i=\overline{1,3}$.
The produced solution serves as the Dirichlet boundary condition for this experiment, with the following coefficients: For $i=\overline{1,3}$, $d_i=1$ for diffusion, $\beta_i=1.0$ for advection, and $\gamma_1=0.001$, $\gamma_2=0.0006$, and $\gamma_3=0.0$ for harvesting. Next, we compare the computed solutions from the Algorithms \ref{Algn1} and \ref{Algn2} with the above manufactured analytical solution. 

In order to quantify the efficiency of the algorithms, we compute the rates of convergence. For spatial convergence, we keep a fixed but small simulation end time $T=0.001$, run the simulations with progressively decreasing the mesh size $h$ by a factor of 2, and record the errors. However, in order to demonstrate the temporal convergence, we run the simulations with a fixed tiny mesh size of $h$, a progressively refined time-step size $\Delta t$ by a factor of 2, and record the errors.  
Table \ref{spatial-convergence-N-3-BE}-\ref{spatial-convergence-N-3-BDF2} represent the spatial errors and convergence rates, revealing second-order spatial convergence across all three species densities for both algorithms, which is consistent with the theoretical results in \eqref{error-estimate-BE}-\eqref{error-estimate-BDF-2}. Additionally, Table \ref{temporal-convergence-N3-BE}-\ref{temporal-convergence-N3-BDF2} present the three species temporal errors and convergence rates, indicating first-order and secon-order convergence of the DBE, and DBDF-2 scheme, respectively, which are consistent with the priori error estimate given in \eqref{error-estimate-BE}-\eqref{error-estimate-BDF-2}.
\begin{table}
	\begin{center}
		\small\begin{tabular}{|c|c|c|c|c|c|c|}\hline
			\multicolumn{2}{|c}{DBE scheme}&\multicolumn{5}{|c|}{Errors and convergence rates}\\ \hline
			$h$ & $\|e_1\|_{2,1}$ & rate   &$\|e_2\|_{2,1}$ & rate &$\|e_3\|_{2,1}$ & rate \\ \hline
			$1/4$ & 2.1871e-5 & & 3.6321e-5 &  & 1.3318e-5 &\\ \hline
			$1/8$ & 5.4646e-6    &2.00& 9.1122e-6    &1.99    &3.3413e-6    &1.99\\ \hline
			$1/16$ & 1.3660e-6    &2.00& 2.2800e-6    &2.00    &8.3606e-7    &2.00\\ \hline
			$1/32$ & 3.4147e-7    &2.00& 5.7012e-7    &2.00    &2.0906e-7    &2.00\\\hline
			$1/64$ & 8.5366e-8    &2.00& 1.4254e-7    &2.00    &5.2269e-8    &2.00\\ \hline
		\end{tabular}
	\end{center}
	\caption{\footnotesize  Spatial errors and convergence rates with $\gamma_1=0.001$, $\gamma_2=0.0006$, $\gamma_3=0.0$, $T=0.001$, and $\Delta t=T/8$.}\vspace{-2ex} \label{spatial-convergence-N-3-BE}
\end{table}

\begin{table}
	\begin{center}
		\small\begin{tabular}{|c|c|c|c|c|c|c|}\hline
			\multicolumn{2}{|c}{DBDF-2 scheme}&\multicolumn{5}{|c|}{Errors and convergence rates}\\\hline
			$h$ & $\|e_1\|_{2,1}$ & rate   &$\|e_2\|_{2,1}$ & rate &$\|e_3\|_{2,1}$ & rate \\ \hline
			$1/4$ & 2.0459e-5 & & 3.3976e-5 &  & 1.2459e-5 & \\ \hline
			$1/8$ & 5.1118e-6    &2.00& 8.5239e-6    &1.99    &3.1256e-6    &1.99\\ \hline
			$1/16$ & 1.2778e-6    &2.00& 2.1328e-6    &2.00    &7.8209e-7    &2.00\\ \hline
			$1/32$ & 3.1942e-7    &2.00& 5.3331e-7    &2.00    &1.9556e-7    &2.00\\ \hline
			$1/64$ & 7.9850e-8    &2.00&  1.3333e-7    &2.00    &4.8890e-8    &2.00\\ \hline
		\end{tabular}
	\end{center}
	\caption{\footnotesize Spatial errors and convergence rates with $\gamma_1=0.001$, $\gamma_2=0.0006$, $\gamma_3=0.0$, $T=0.001$, and $\Delta t=T/16$.} \vspace{-4ex}\label{spatial-convergence-N-3-BDF2}
\end{table}

\begin{table}
	\begin{center}
		\small\begin{tabular}{|c|c|c|c|c|c|c|}\hline
			\multicolumn{2}{|c}{DBE scheme}&\multicolumn{5}{|c|}{Error and convergence rates}\\\hline
			$\Delta t$ & $\|e_1\|_{2,1}$ & rate   &$\|e_2\|_{2,1}$ & rate &$\|e_3\|_{2,1}$ & rate  \\ \hline
			$T/4$ & 5.5489e-1 & & 7.5417e-1 &  & 4.3128e-1 &\\ \hline
			$T/8$ & 2.2273e-1    &1.32& 3.2590e-1    &1.21    &1.7299e-1    &1.32\\ \hline
			$T/16$ & 1.0245e-1    &1.12& 1.5380e-1    &1.08    &7.9560e-2    &1.12\\ \hline
			$T/32$ & 4.9472e-2    &1.05& 7.5182e-2    &1.03    &3.8416e-2    &1.05\\ \hline
			$T/64$ & 2.4367e-2    &1.02& 3.7276e-2    &1.01    &1.8922e-2    &1.02\\ \hline
			$T/128$ & 1.2103e-2    &1.01& 1.8582e-2    &1.00    &9.3983e-3    &1.01\\ \hline
		\end{tabular}
	\end{center}
	\caption{\footnotesize Temporal convergence with $\gamma_1=0.001$, $\gamma_2=0.0006$, $\gamma_3=0.01$, and $h = 1/64$.}\vspace{-2ex}\label{temporal-convergence-N3-BE}
\end{table}

\begin{table}
	\begin{center}
		\small\begin{tabular}{|c|c|c|c|c|c|c|}\hline
			\multicolumn{2}{|c}{DBDF-2 scheme}&\multicolumn{5}{|c|}{Errors and convergence rates}\\\hline
			$\Delta t$ & $\|e_1\|_{2,1}$ & rate   &$\|e_2\|_{2,1}$ & rate &$\|e_3\|_{2,1}$ & rate  \\ \hline
			$T/4$ & 4.2383e-1 & & 6.4258e-1 &  & 3.2889e-1 &\\ \hline
			$T/8$ & 8.9251e-2    &2.25& 1.5478e-1    &2.05    &6.8971e-2    &2.25\\ \hline
			$T/16$ & 2.1659e-2    &2.04& 3.9555e-2    &1.97    &1.6707e-2    &2.05\\ \hline
			$T/32$ & 5.5089e-3    &1.98& 1.0285e-2    &1.94    &4.2421e-3    &1.98\\ \hline
			$T/64$ & 1.4050e-3    &1.97& 2.6519e-3    &1.96    &1.0800e-3    &1.97\\ \hline
			$T/128$ & 3.5669e-4    &1.98& 6.7677e-4    &1.97    &2.7366e-4    &1.98\\ \hline
		\end{tabular}
	\end{center}
	\caption{\footnotesize Temporal convergence with $\gamma_1=0.001$, $\gamma_2=0.0006$, $\gamma_3=0.01$, and  $h = 1/128$.}\vspace{-4ex} \label{temporal-convergence-N3-BDF2}
\end{table}

In all the tests henceforth, we consider the DBDF-2 scheme with time-step size $\Delta t=0.1$, and initial density ${u}^0_{i,h}=1.6$. In each of the test, we consider either the following non-stationary carrying capacity
\begin{align}
    K(t,\bx)\equiv\big(1.2+2.5\pi^2e^{-(x-0.5)^2-(y-0.5)^2}\big)\big(1.0+0.3\cos t\big),\label{non-stationary-carrying-capacity}
\end{align}
or the stationary carrying capacity
\begin{align}
    K(t,\bx)\equiv1.2+2.5\pi^2e^{-(x-0.5)^2-(y-0.5)^2}.\label{stationary-carrying-capacity}
\end{align}
We use the non-stationary carrying capacity for the simulations, unless otherwise stated. To have a second required initial condition of the DBDF-2 scheme, we run DBE scheme.

\subsection{Two species competition model: Effect of advection on the evolution of population density}
\begin{figure}
	\centering
	\subfloat[]{\includegraphics[width=0.43\textwidth,height=0.27\textwidth]{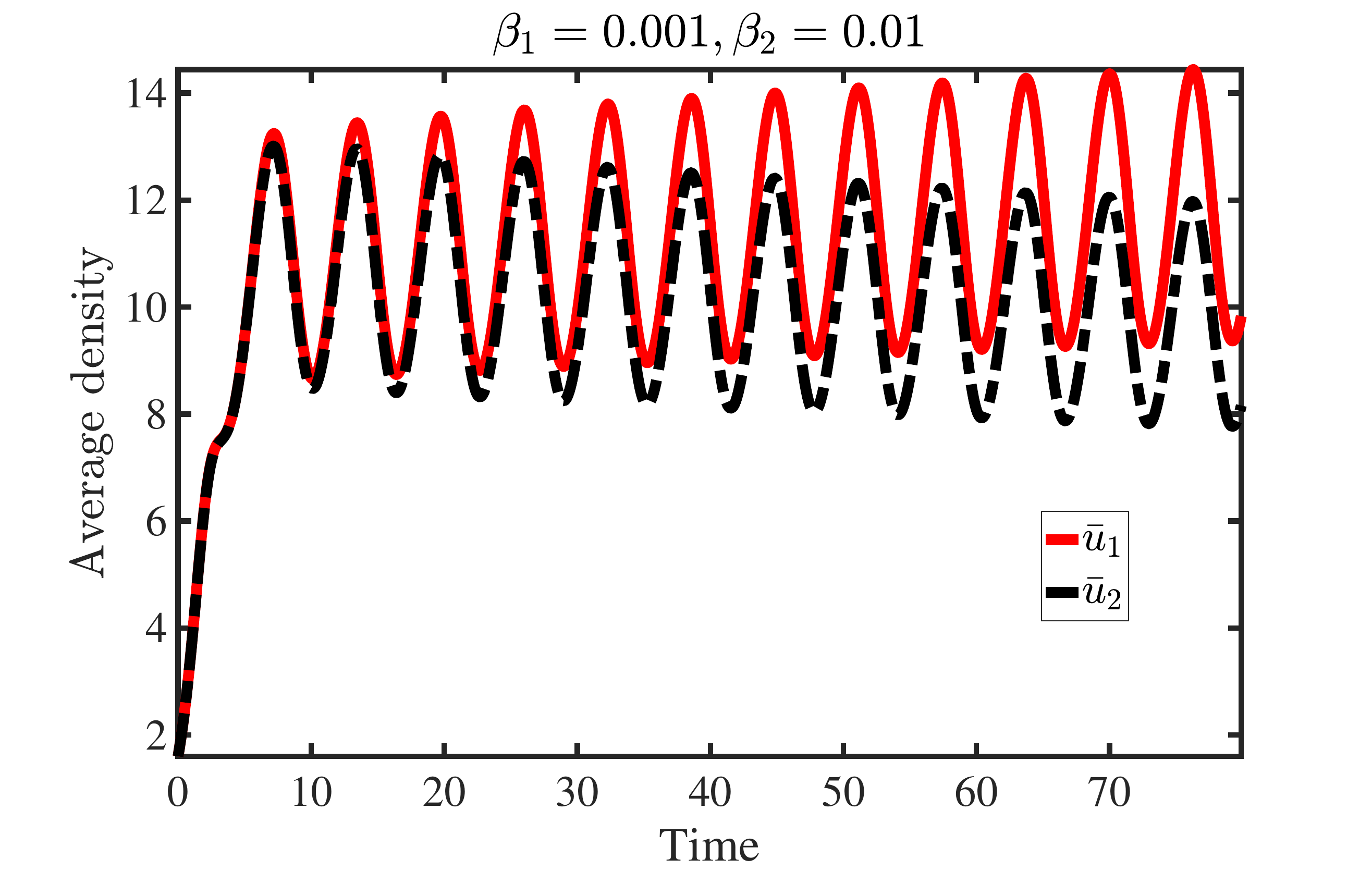}}
	\subfloat[]{\includegraphics[width=0.43\textwidth,height=0.27\textwidth]{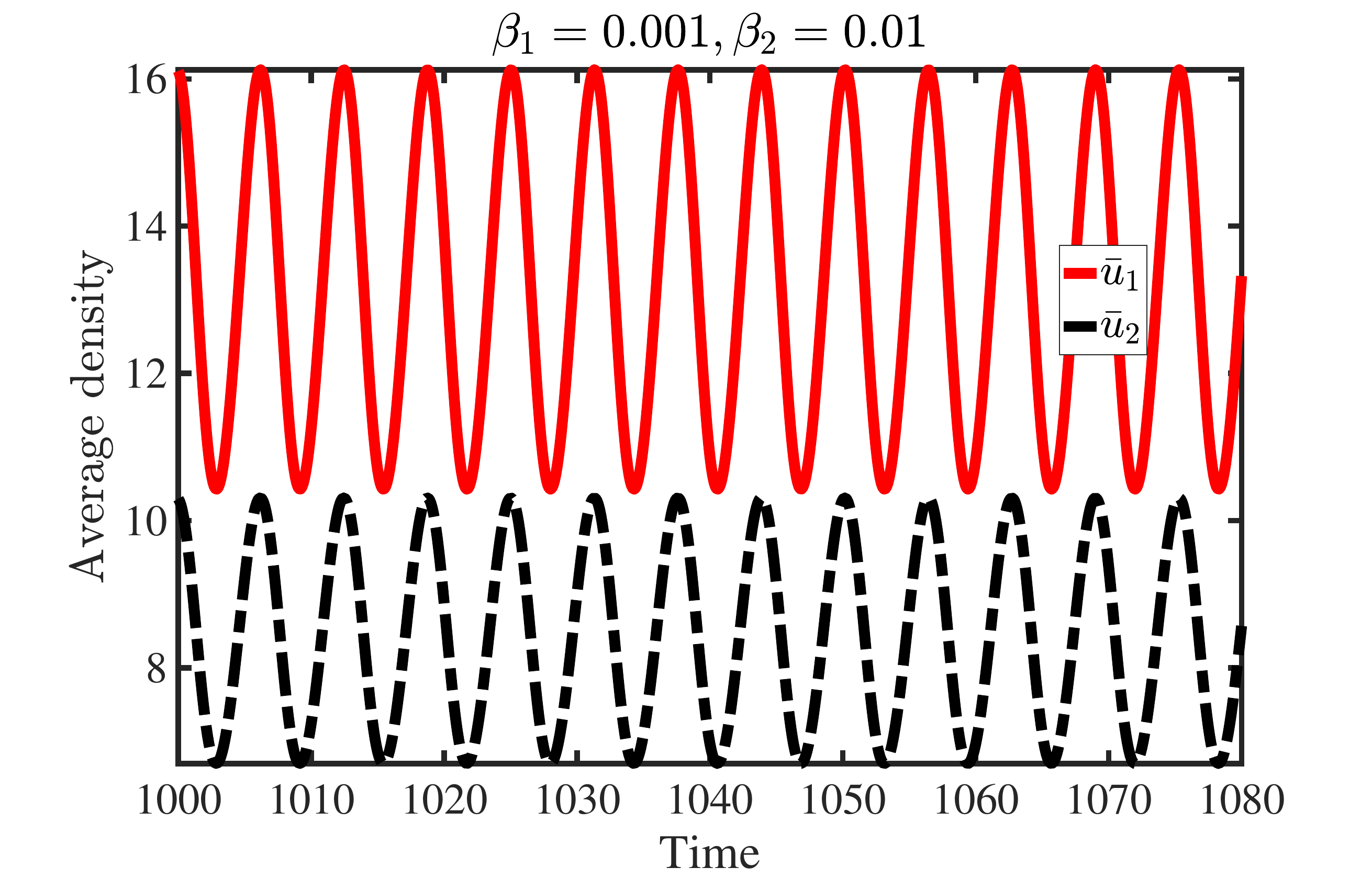}}\vspace{-4ex}
	\caption{Average density of each species: (a) short-, and  (b) long-range with the advection coefficients $\beta_1=0.001$, and $\beta_2=0.01$, intrinsic growth rates $r_i=1$, and harvesting coefficients $\gamma_i=0.0$, for $i=\overline{1,2}$.}\vspace{-4ex}	\label{2S: RDE-energy-b1_0_001_b2_0_01}
\end{figure}

\begin{figure}
	\centering
	\subfloat[]{\includegraphics[width=0.43\textwidth,height=0.27\textwidth]{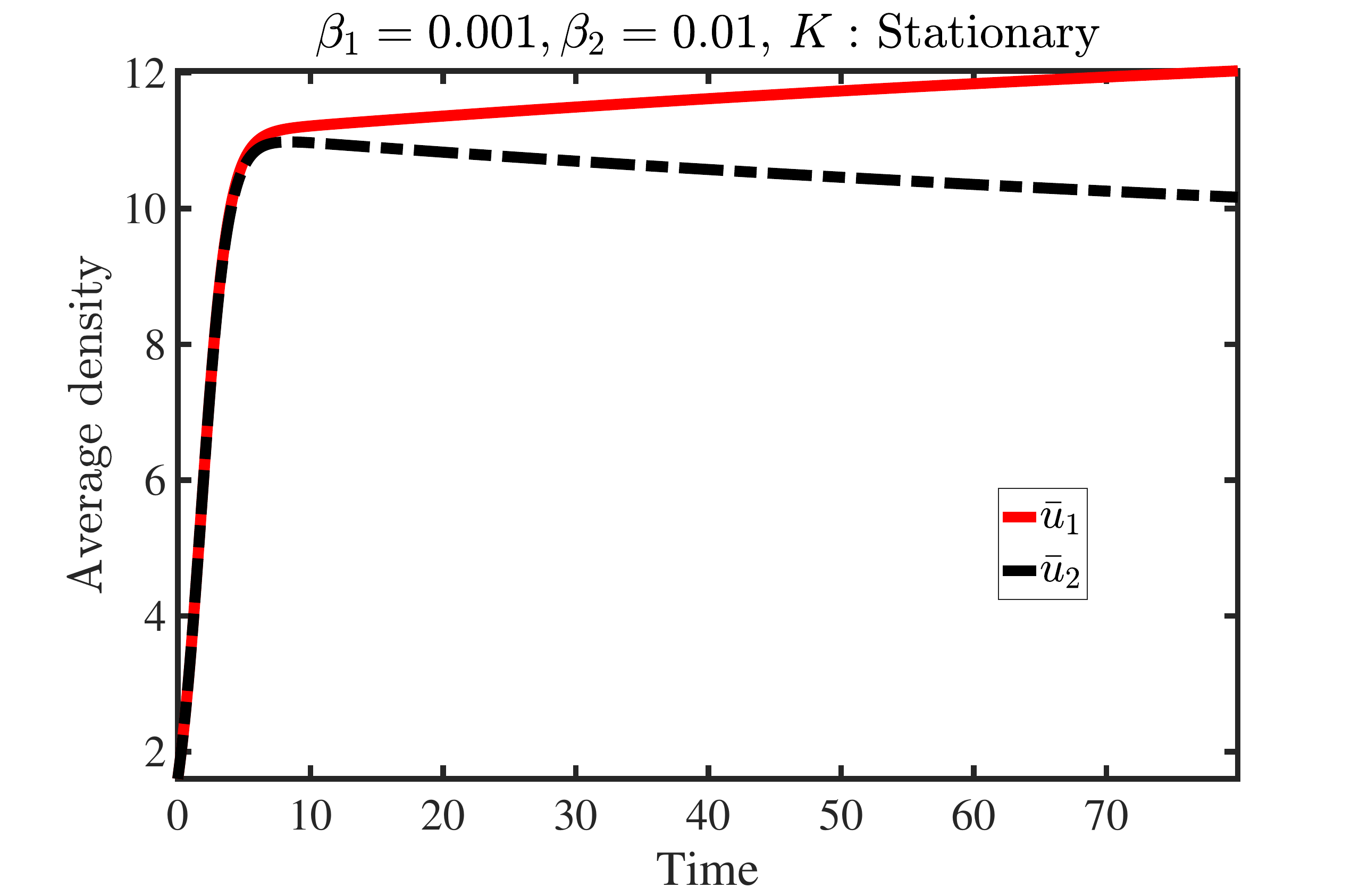}}
	\subfloat[]{\includegraphics[width=0.43\textwidth,height=0.27\textwidth]{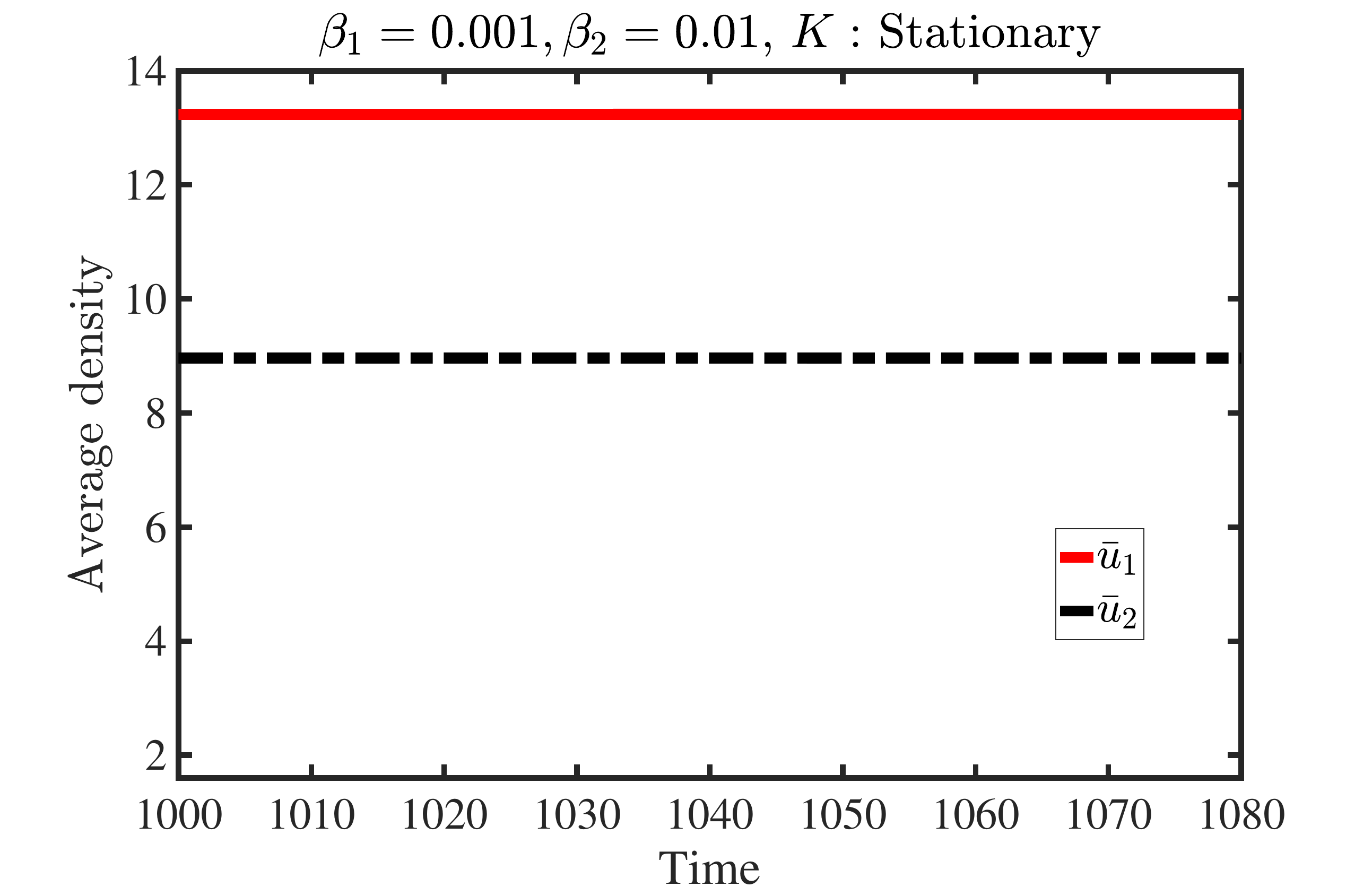}}\vspace{-4ex}
	\caption{Average density of each species: (a) short-, and  (b) long-range with the advection coefficients $\beta_1=0.001$, and $\beta_2=0.01$, intrinsic growth rates $r_i=1$, harvesting coefficients $\gamma_i=0.0$, for $i=\overline{1,2}$, and stationary carrying capacity.}\vspace{-4ex}	\label{2S: RDE-energy-b1_0_001_b2_0_01_K_stationary}
\end{figure}

In this test, we consider the model for two species populations in absence of harvesting (e.g., $\gamma_i=0$) with equal diffusion rate $d_1=d_2=0.1$. To this end, we examine how the average population density of a species changes with the advection parameter. We plot the average density versus time in Figure \ref{2S: RDE-energy-b1_0_001_b2_0_01} for both species varying the advection parameters as $\beta_1=0.001$, and $\beta_2=0.1$. Figure \ref{2S: RDE-energy-b1_0_001_b2_0_01} suggests that the initial value is not important to the final state due to the global convergence of solutions. We observe that as the advection parameter increases, the species density decreases over time. The species with lower advection rate will converge to the stable solution faster and there is a strong co-existence of both the species. Figure \ref{2S: RDE-energy-b1_0_001_b2_0_01}(a), (b)  are plotted for the same data, but for short and long time scenarios. Figure \ref{2S: RDE-energy-b1_0_001_b2_0_01_K_stationary} is plotted for stationary carrying capacity and observed the similar outcomes.


\subsection{Two species competition model: Effect of stocking and harvesting on the evolution of population density}
\begin{figure}
	\centering
	\subfloat[]{\includegraphics[width=0.43\textwidth,height=0.27\textwidth]{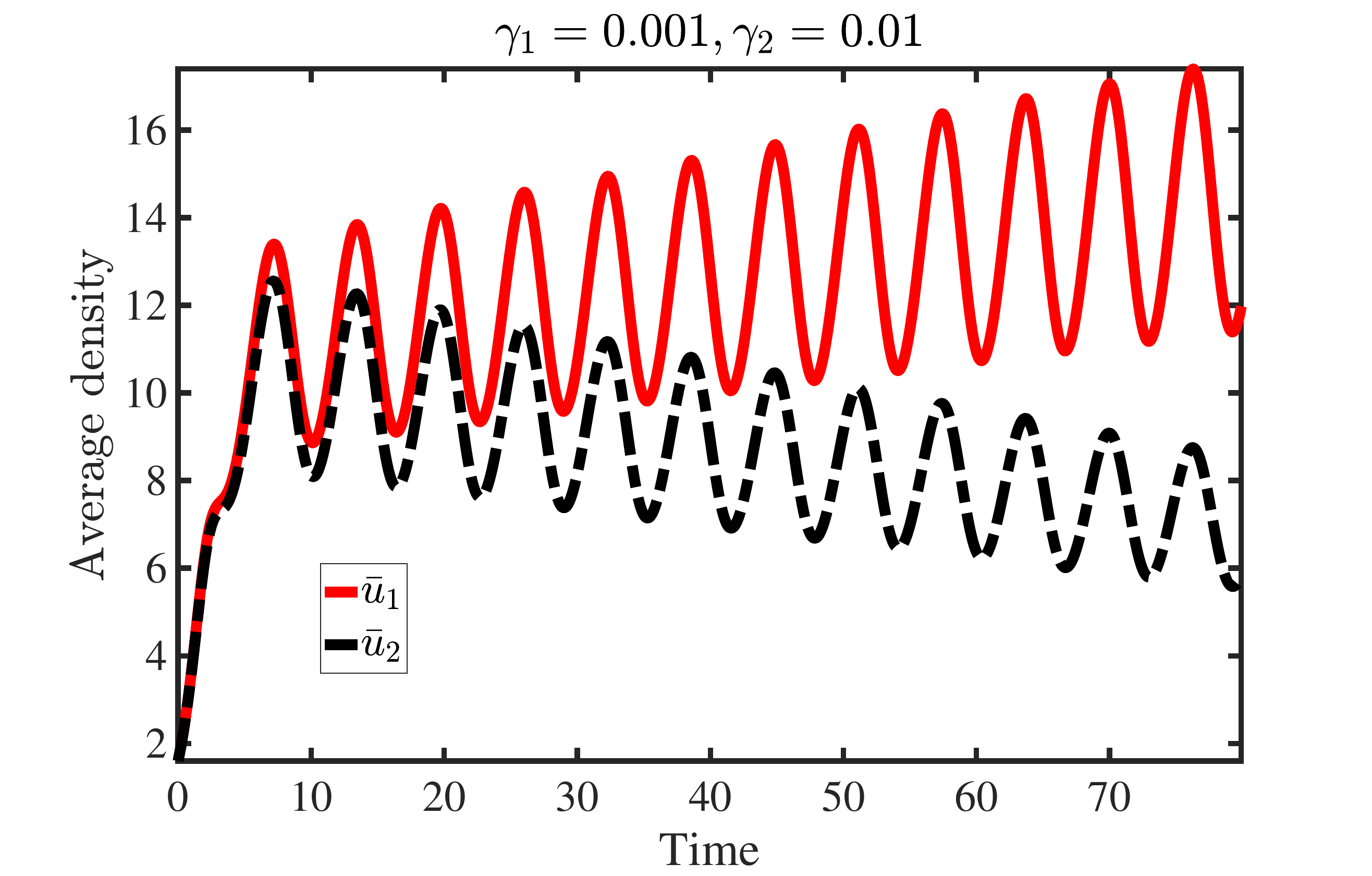}}
	\subfloat[]{\includegraphics[width=0.43\textwidth,height=0.27\textwidth]{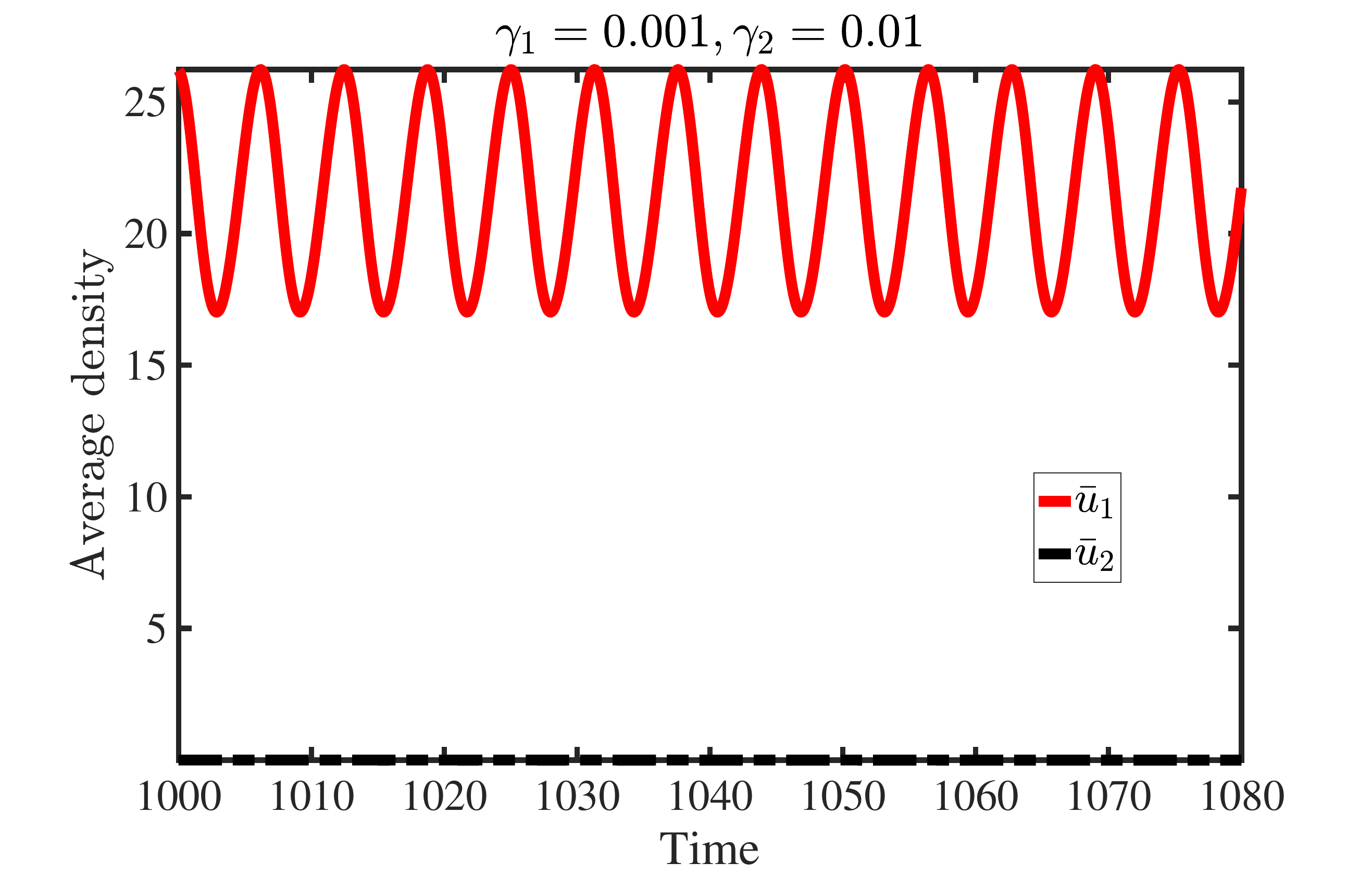}}\vspace{-4ex}
	\caption{Average density of each species: (a) short-, and  (b) long-range with the harvesting coefficients $\gamma_1=0.001$, and $\gamma_2=0.01$ for fixed and equal growth rate and advection coefficients.}	\vspace{-4ex}\label{2S: RDE-energy-g1_0_001_g2_0_01}
\end{figure}

For this experiment, we consider constant intrinsic growth rates $ r_i \equiv 1 $, and zero advection rates ($\beta_i = 0 $,
for i = 1, 2). We consider the carrying capacity given in \eqref{non-stationary-carrying-capacity} and the harvesting coefficients $\gamma_1 = 0.001$, and $\gamma_2=0.01$
in Figures \ref{2S: RDE-energy-g1_0_001_g2_0_01}(a), where the average density of each species versus time is plotted for time $t = 0$ to $ 80 $. From the average density plot, we observe periodic population densities for both species, where the density of $u_2$ is decreasing because of its higher harvesting coefficient (Figure \ref{2S: RDE-energy-g1_0_001_g2_0_01}(a)). It is predicted that the species $ u_2 $ will die out if time is too large (Figure \ref{2S: RDE-energy-g1_0_001_g2_0_01}(b)).

\begin{figure}
	\centering
	\subfloat[]{\includegraphics[width=0.43\textwidth,height=0.27\textwidth]{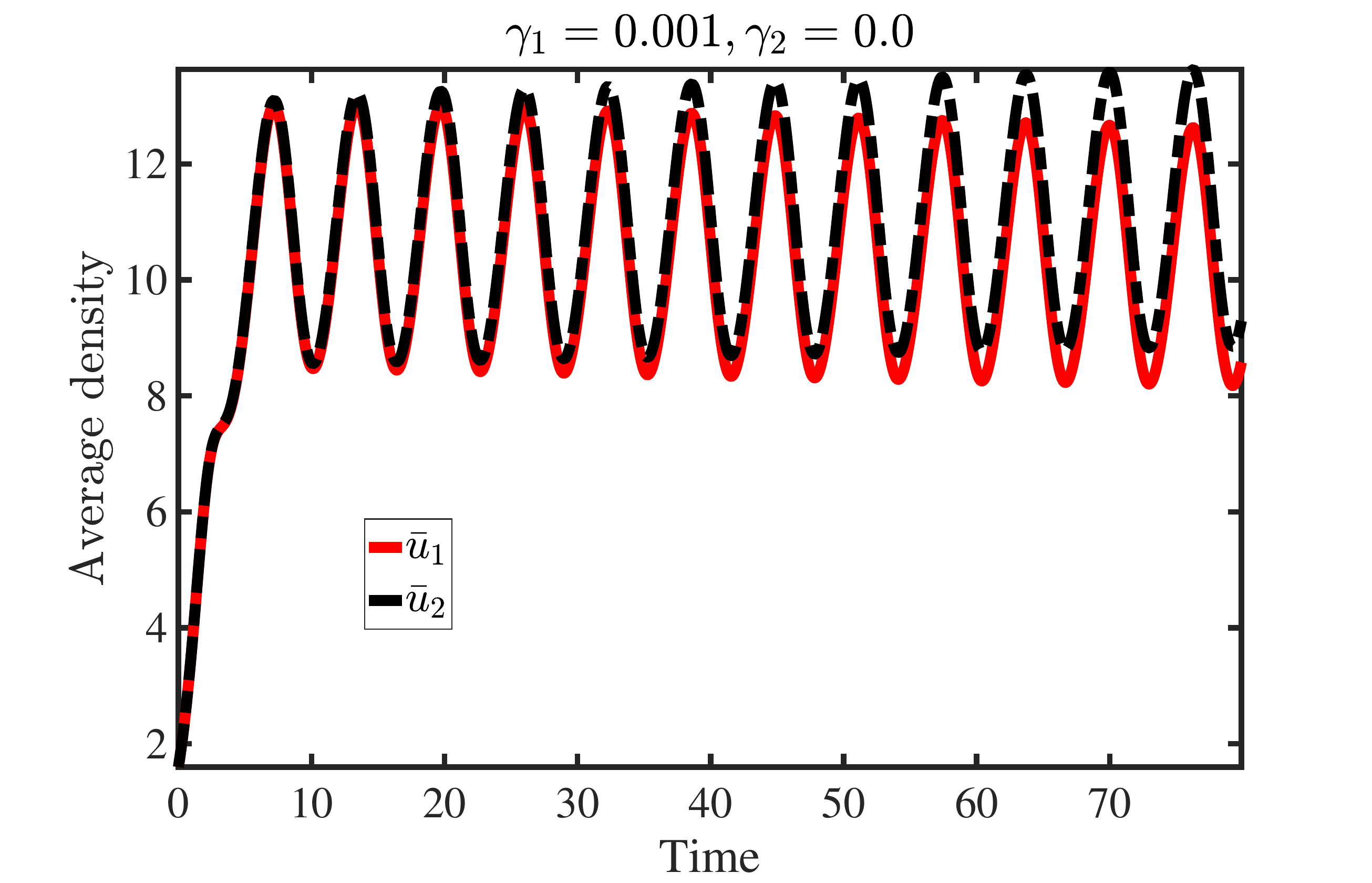}}
	\subfloat[]{\includegraphics[width=0.43\textwidth,height=0.27\textwidth]{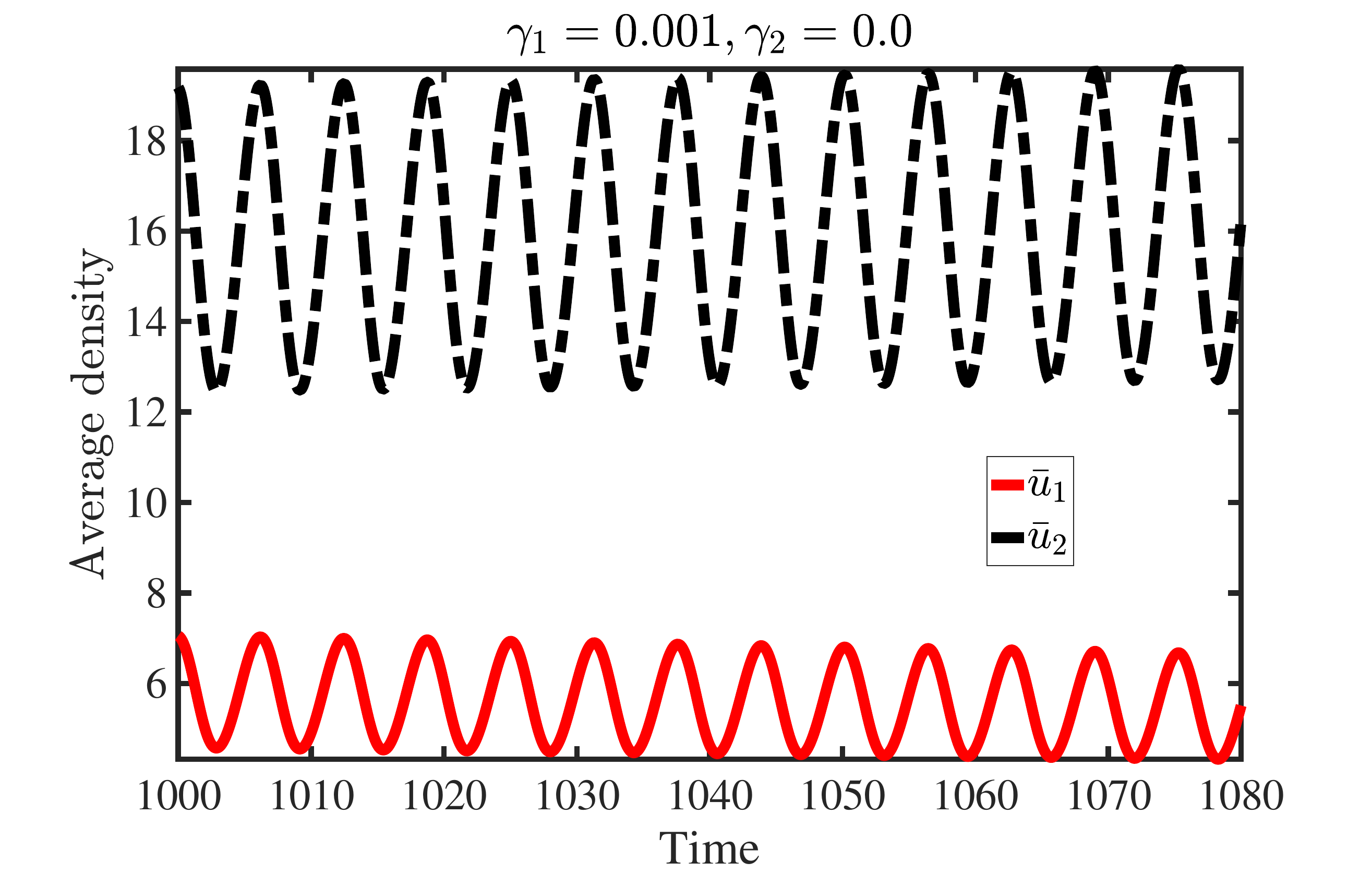}}
	\vspace{-4ex}\caption{Average density of each species: (a) short-, and  (b) long-range with the harvesting coefficients $\gamma_1=0.001$, and $\gamma_2=0.0$ (no harvesting) for fixed and equal growth rate and advection coefficients.}\vspace{-4ex}	\label{2S: RDE-energy-g1_0_001_g2_0}
\end{figure}

\begin{figure}
	\centering
	\subfloat[]{\includegraphics[width=0.43\textwidth,height=0.27\textwidth]{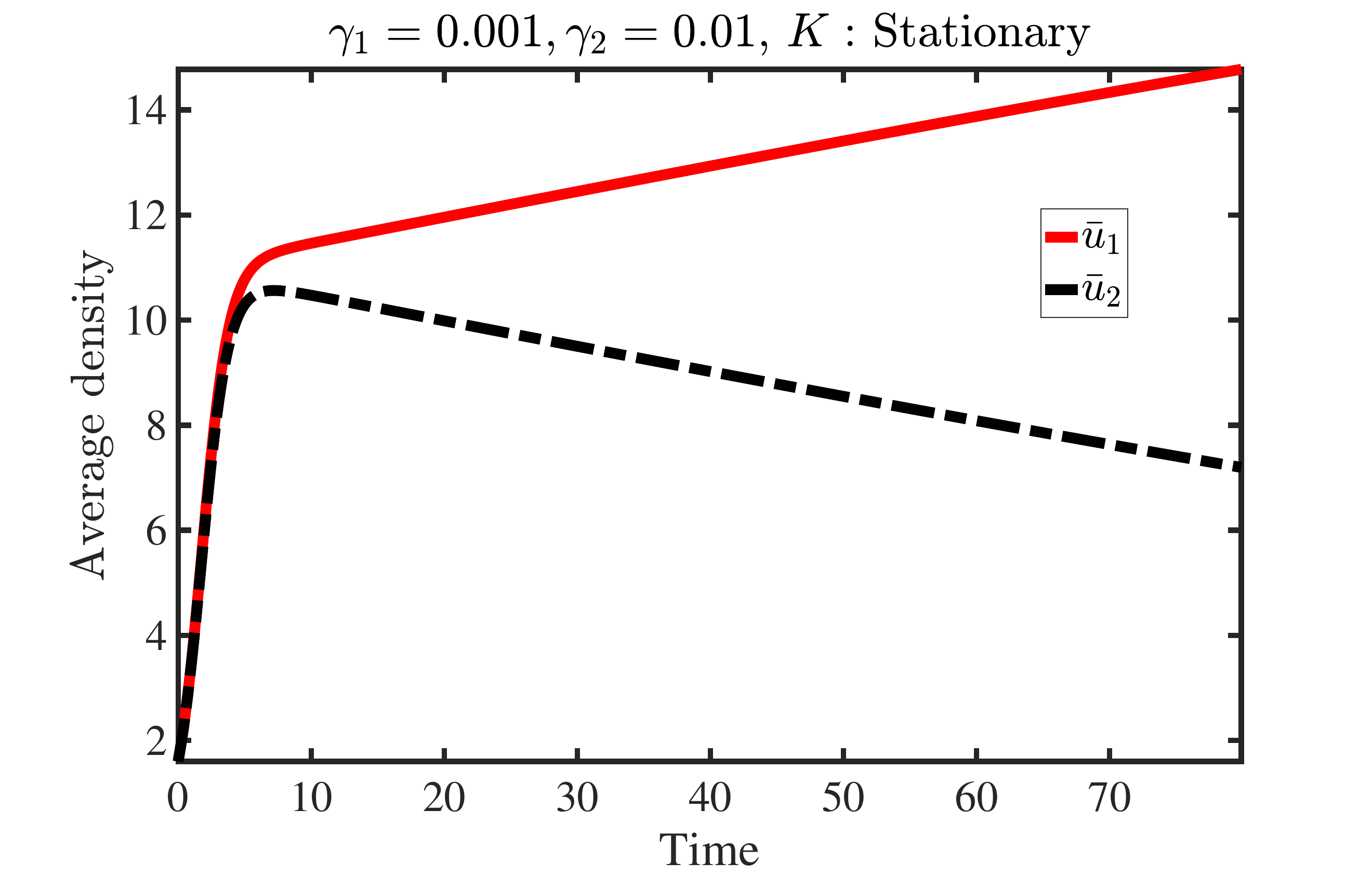}}
	\subfloat[]{\includegraphics[width=0.43\textwidth,height=0.27\textwidth]{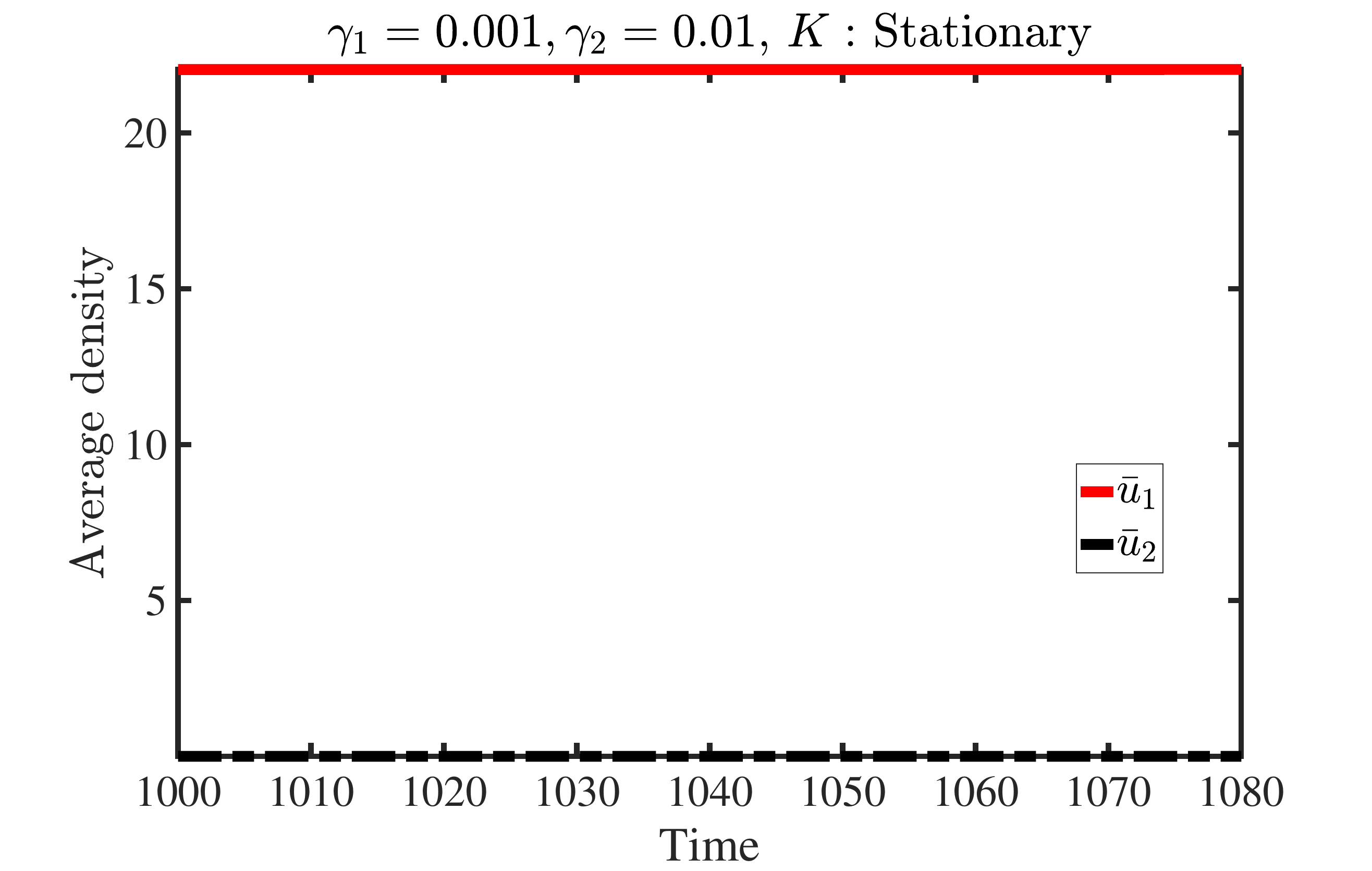}}
	\vspace{-4ex}\caption{Average density of each species: (a) short-, and  (b) long-range with the harvesting coefficients $\gamma_1=0.001$, and $\gamma_2=0.0$, intrinsic growth rates $r_i=1$, advection coefficients $\beta_i=0.0$, for $i=\overline{1,2}$, and stationary carrying capacity.}	\vspace{-4ex}\label{2S: RDE-energy-g1_0_001_g2_0_K_stationary}
\end{figure}

If we further consider the harvesting coefficients as $\gamma_1 = 0.001$, and $\gamma_2=0.0$ (no harvesting) and plot the average density curves versus time for each species in Figure \ref{2S: RDE-energy-g1_0_001_g2_0}, we observe an evolutionary population density feature, especially for the first species. The harvesting in the first species provides an advantage to the second species, and thus an apparent co-existence of both species is visible over the time [0, 80] (Figure \ref{2S: RDE-energy-g1_0_001_g2_0}(a)). Figure \ref{2S: RDE-energy-g1_0_001_g2_0}(b) represents the average density of each of both species on the time interval [1000, 1080], where only the first species is affected by harvesting with coefficient $\gamma_1 = 0.001 $ and it will extinct over time. 
It is observed from Figure \ref{2S: RDE-energy-g1_0_001_g2_0_K_stationary} that for stationary carrying capacity the species with lower harvesting rate will win the competition and for short time period there exists a co-existence between them. On the other hand for long time period the species with higher harvesting rate will extinct.

\begin{figure}
	\centering
	\subfloat[]{\includegraphics[width=0.43\textwidth,height=0.27\textwidth]{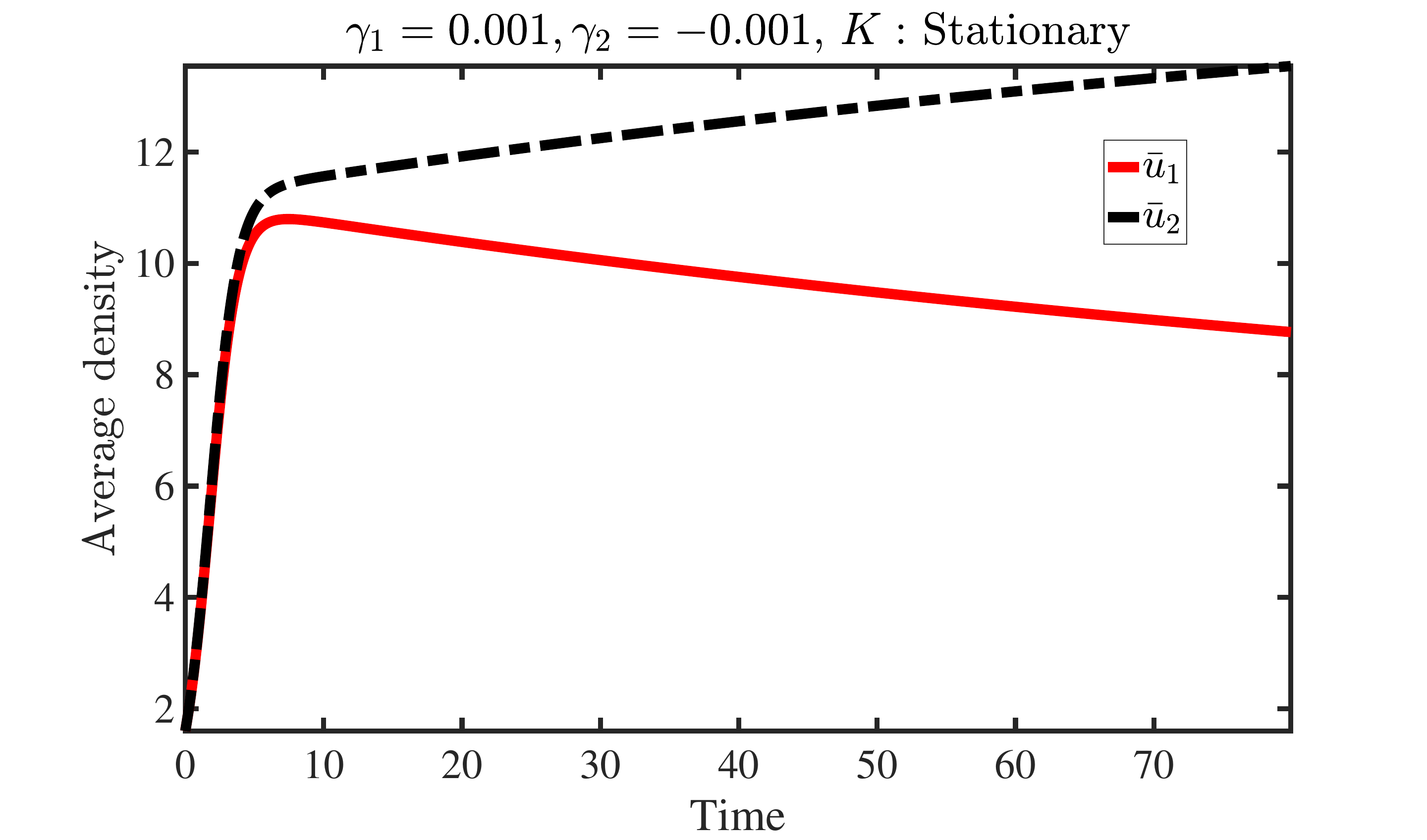}}
	\subfloat[]{\includegraphics[width=0.43\textwidth,height=0.27\textwidth]{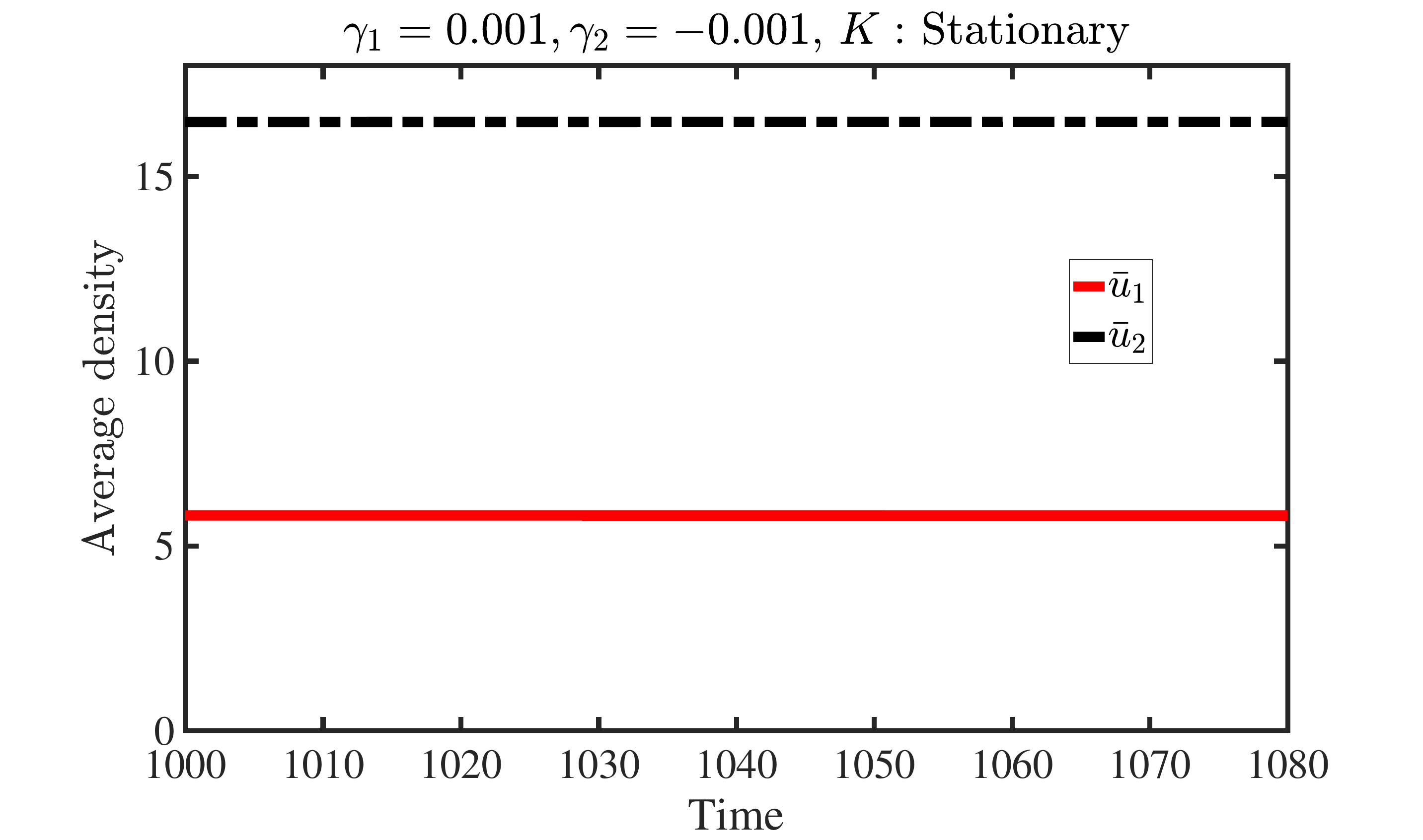}}
	\vspace{-4ex}\caption{Average density of each species: (a) short-, and  (b) long-range with the harvesting coefficient $\gamma_1=0.001$, and stocking coefficient $\gamma_2=-0.001$, intrinsic growth rates $r_i=1$, advection coefficients $\beta_i=0.0$, for $i=\overline{1,2}$, and stationary carrying capacity.}\vspace{-4ex}	\label{2S: RDE-energy-g1_0_001_g2_m0_001_K_stationary}
\end{figure}

Finally, we harvest the first species (with $\gamma_1=0.001)$ and stock the second species (with $\gamma_2=0.001$) and plot the average density versus time curve for each species in Figure \ref{2S: RDE-energy-g1_0_001_g2_m0_001_K_stationary}. This shows an opposite scenario of the case in Figure \ref{2S: RDE-energy-g1_0_001_g2_0_K_stationary}. That is, the second species gets an advantage in the competition, though initially both species coexist, in the long run, the first species will be extinct. 
\subsection{Three species competition model: Evolution of Population Density with Varying Advection Coefficients}\label{exp-2}
\begin{figure}
	\centering
	\subfloat[]{\includegraphics[width=0.43\textwidth,height=0.27\textwidth]{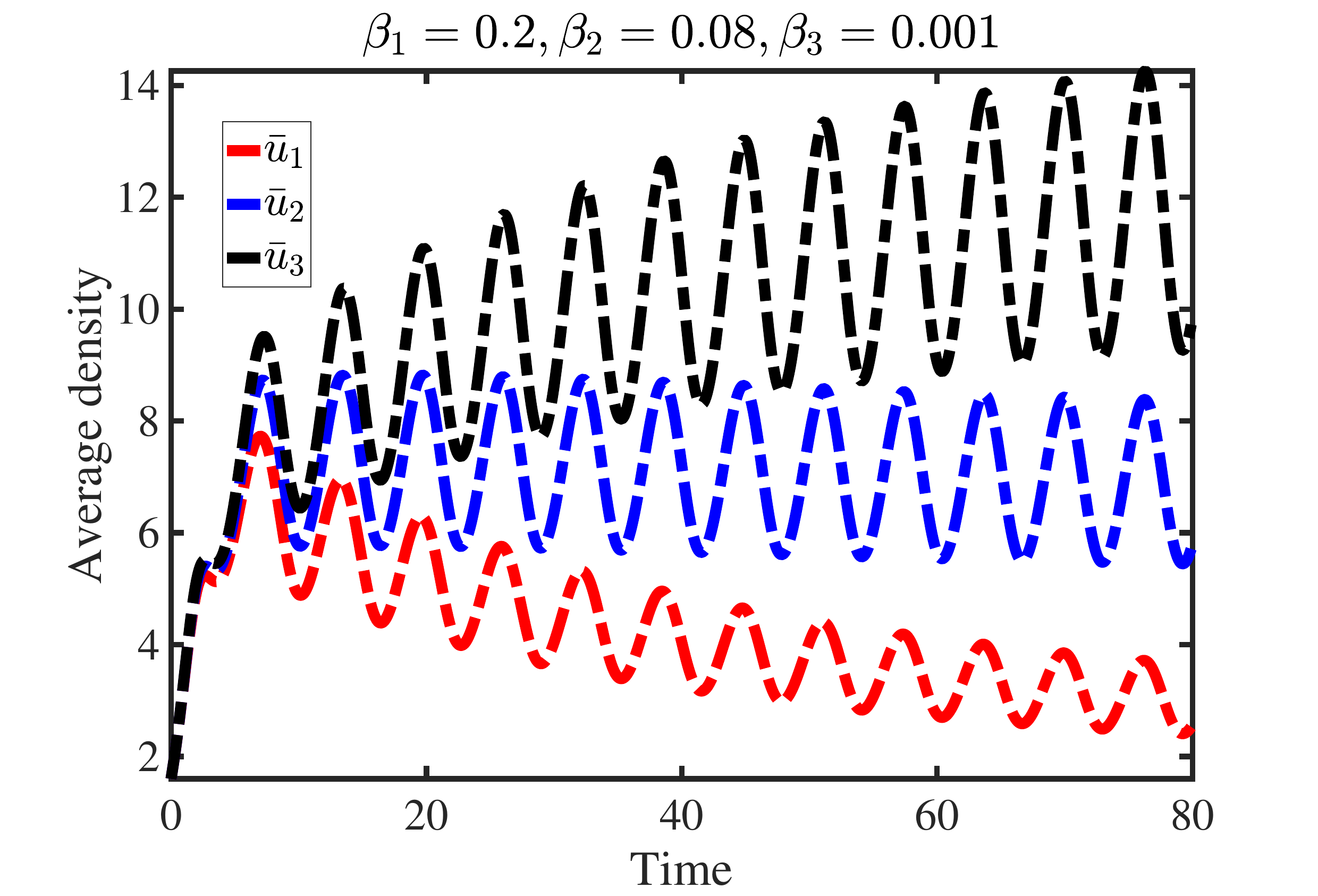}}
	\subfloat[]{\includegraphics[width=0.43\textwidth,height=0.27\textwidth]{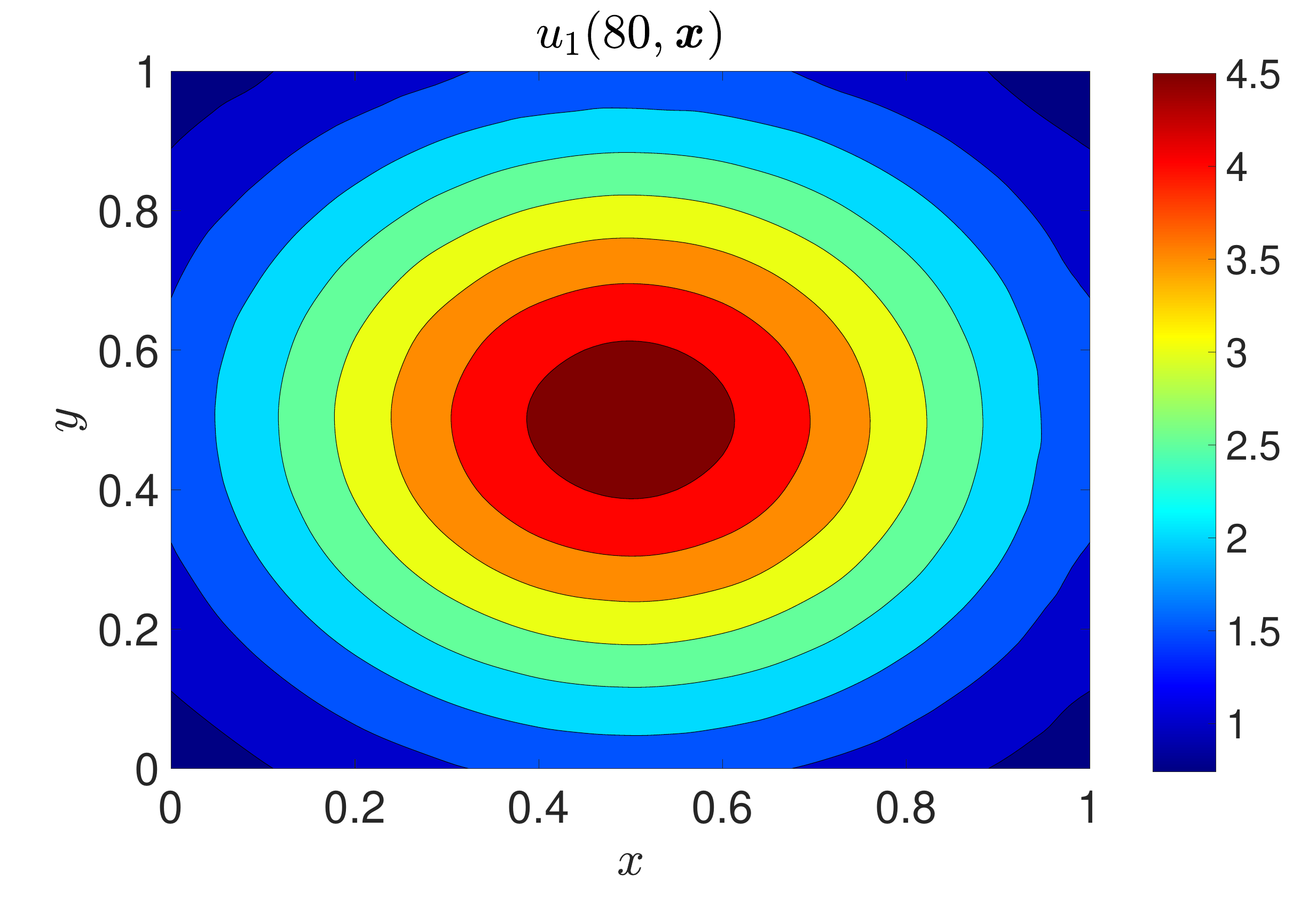}}\\
	\subfloat[]{\includegraphics[width=0.43\textwidth,height=0.27\textwidth]{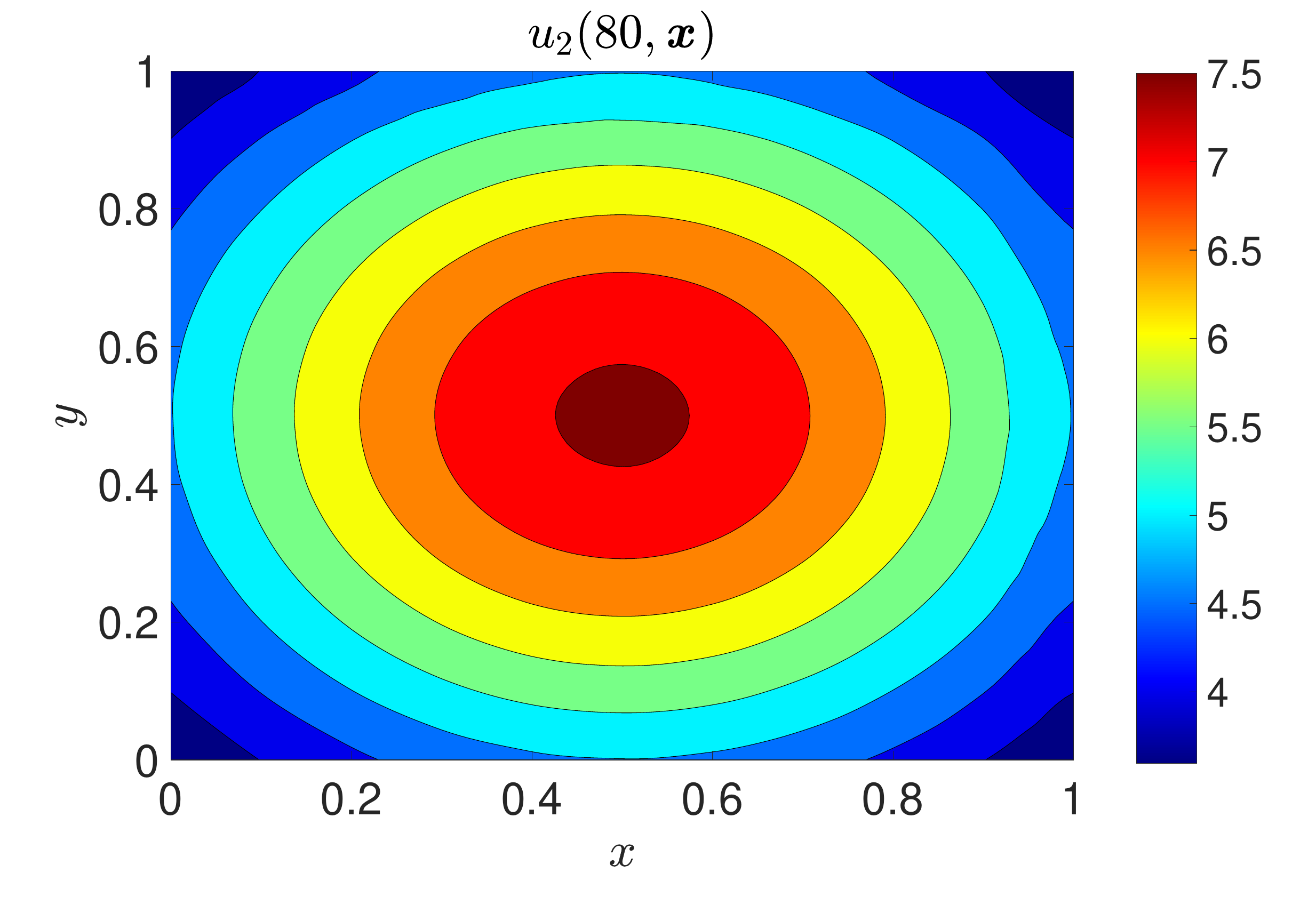}}
	\subfloat[]{\includegraphics[width=0.43\textwidth,height=0.27\textwidth]{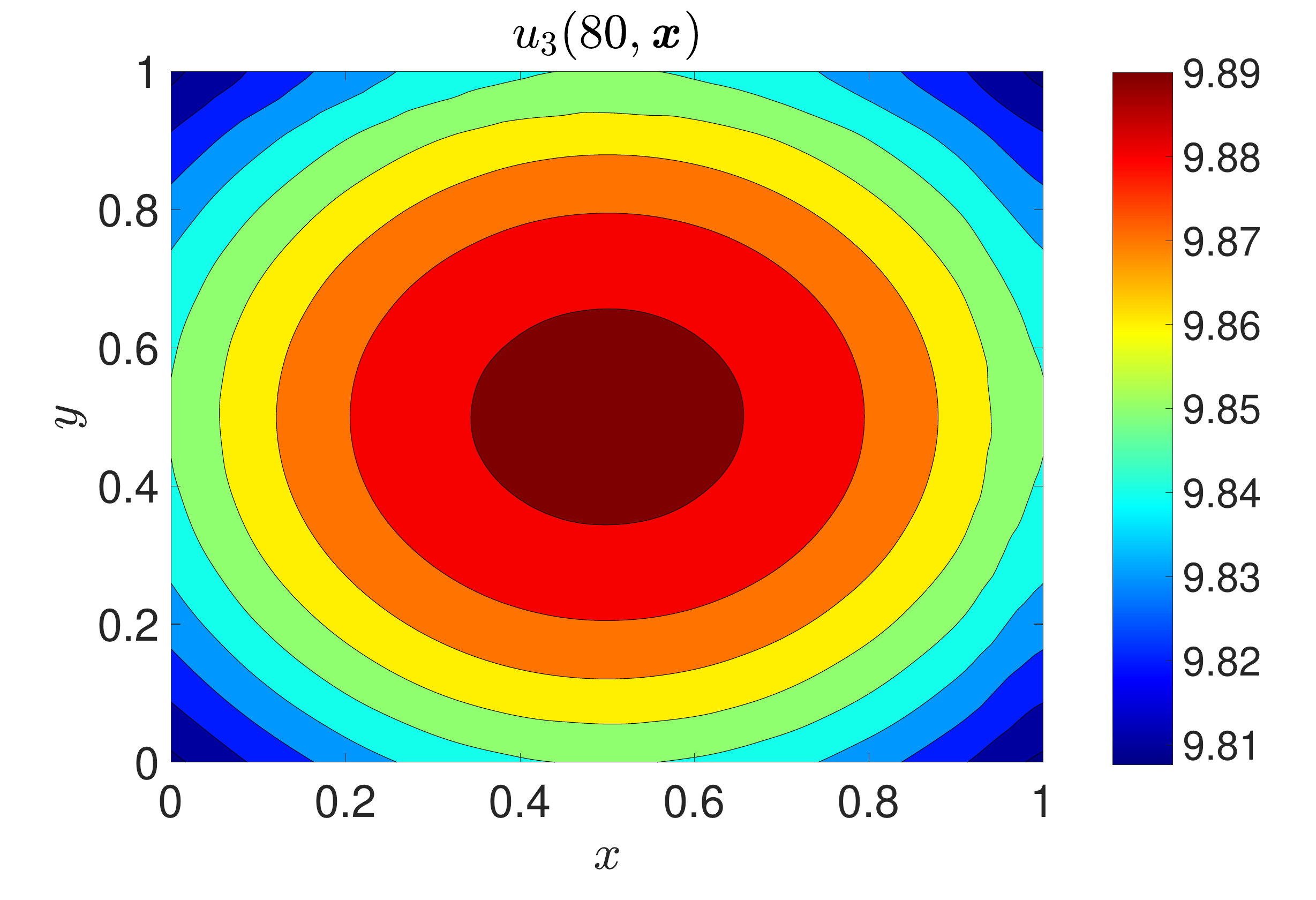}}
	\vspace{-4ex}\caption{ (a) Average density of each species,   (b) contour plot of species density $u_1$, (c) contour plot of species density $u_2$, and (d) contour plot of species density $u_3$ with the advection coefficients $\beta_1=0.2$, $\beta_2=0.08$, and $\beta_3=0.001$, diffusion coefficients $d_i=1$, $r_i=1$, and harvesting coefficients $\gamma_i=0.001$, $i=1,2,3$.}\vspace{-4ex}	\label{RDE-energy-beta_0_2_beta_0_08_beta_0_001-K-time-exp}
\end{figure}

In this test, we consider  constant intrinsic growth rates $r_i=1.0$, diffusion coefficients $d_i=1.0$, harvesting coefficients $\gamma_i=0.001$,  $i=\overline{1,3}$ and varying advection coefficients $\beta_1=0.2$, $\beta_2=0.08$, and $\beta_3=0.001$, and run the simulation until $T=80$. In Figure \ref{RDE-energy-beta_0_2_beta_0_08_beta_0_001-K-time-exp}(a), we represent the average population density over time.
In Figures \ref{RDE-energy-beta_0_2_beta_0_08_beta_0_001-K-time-exp}, the average density of each species versus time is plotted for time $ t = 0 $ to $ 80 $, and the population density contour plot of each of the species at time $ t = 80 $. From the average density plot, we observe periodic population densities for all species, where the density of $ u_1 $ is decreasing because of its higher advection coefficient (Figure \ref{RDE-energy-beta_0_2_beta_0_08_beta_0_001-K-time-exp}(a)). It is predicted that the species $ u_1 $ will die out if time is too large. From the contour plots, it is observed that the highest population density is at the point $ (0.5, 0.5) $ and there is a coexistence of all species, though the population density of the species $ u_3 $ remains bigger than the species $ u_2 $, and $ u_1 $ over the domain (Figure \ref{RDE-energy-beta_0_2_beta_0_08_beta_0_001-K-time-exp}(b)-(d)). This happens because of different advection parameters, and the optimal value of the carrying capacity function is achieved at the point $ (0.5, 0.5) $, which shows the symmetric distribution of the population.

\subsection{Three species competition model: Evolution of Population Density with Varying Harvesting Coefficients}
\begin{figure}
	\centering
	\subfloat[]{\includegraphics[width=0.43\textwidth,height=0.27\textwidth]{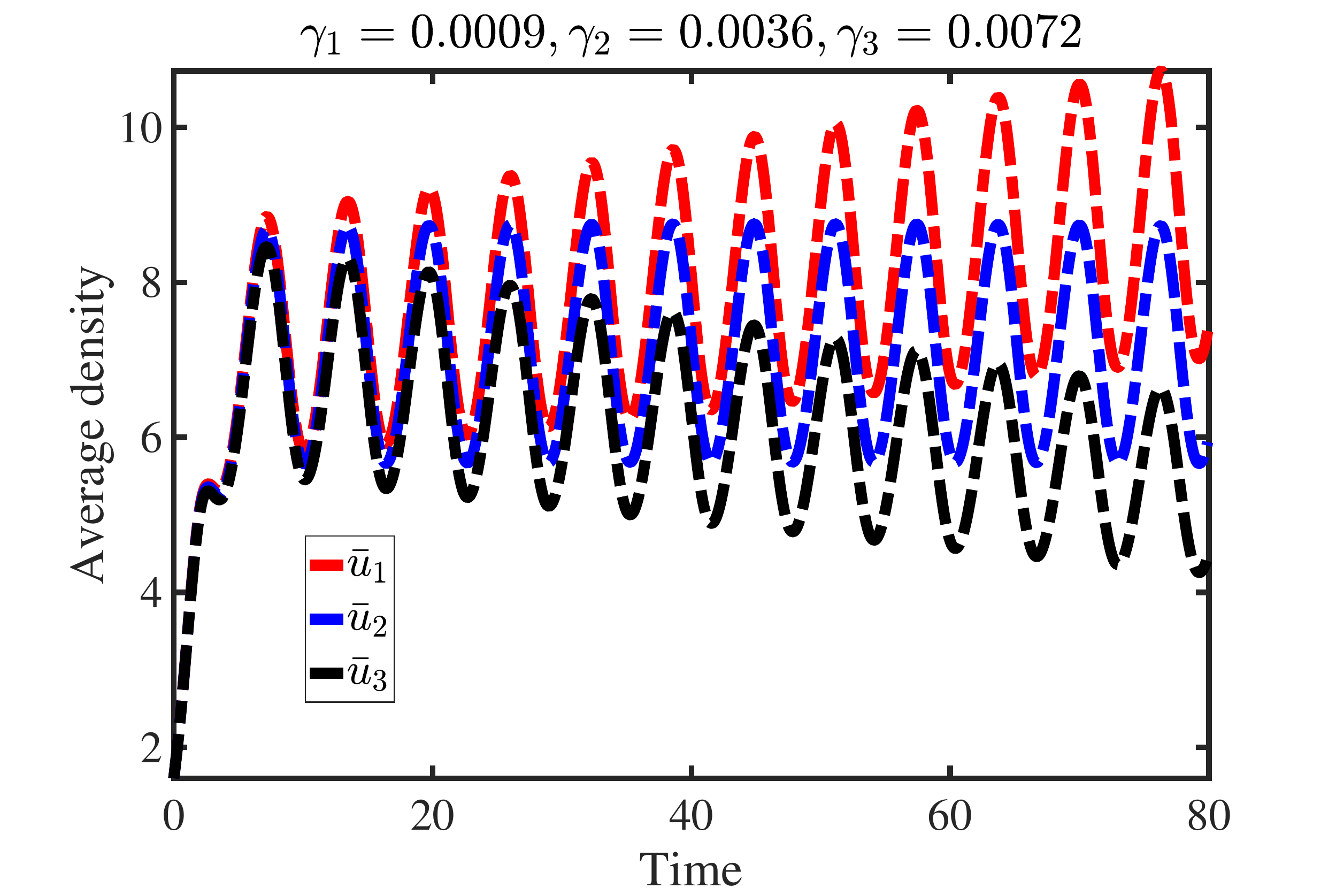}}
	\subfloat[]{\includegraphics[width=0.43\textwidth,height=0.27\textwidth]{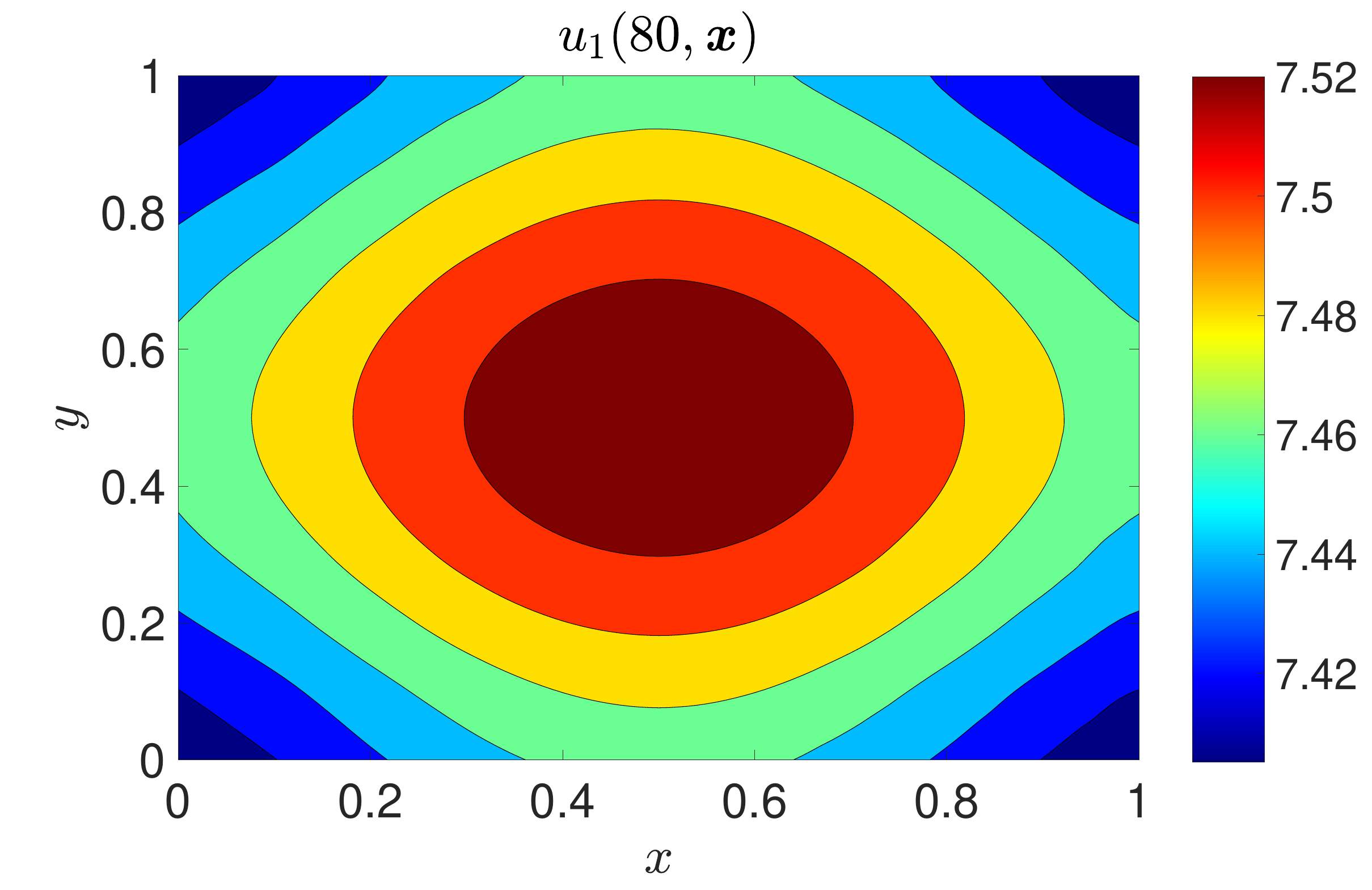}}\\
	\subfloat[]{\includegraphics[width=0.43\textwidth,height=0.27\textwidth]{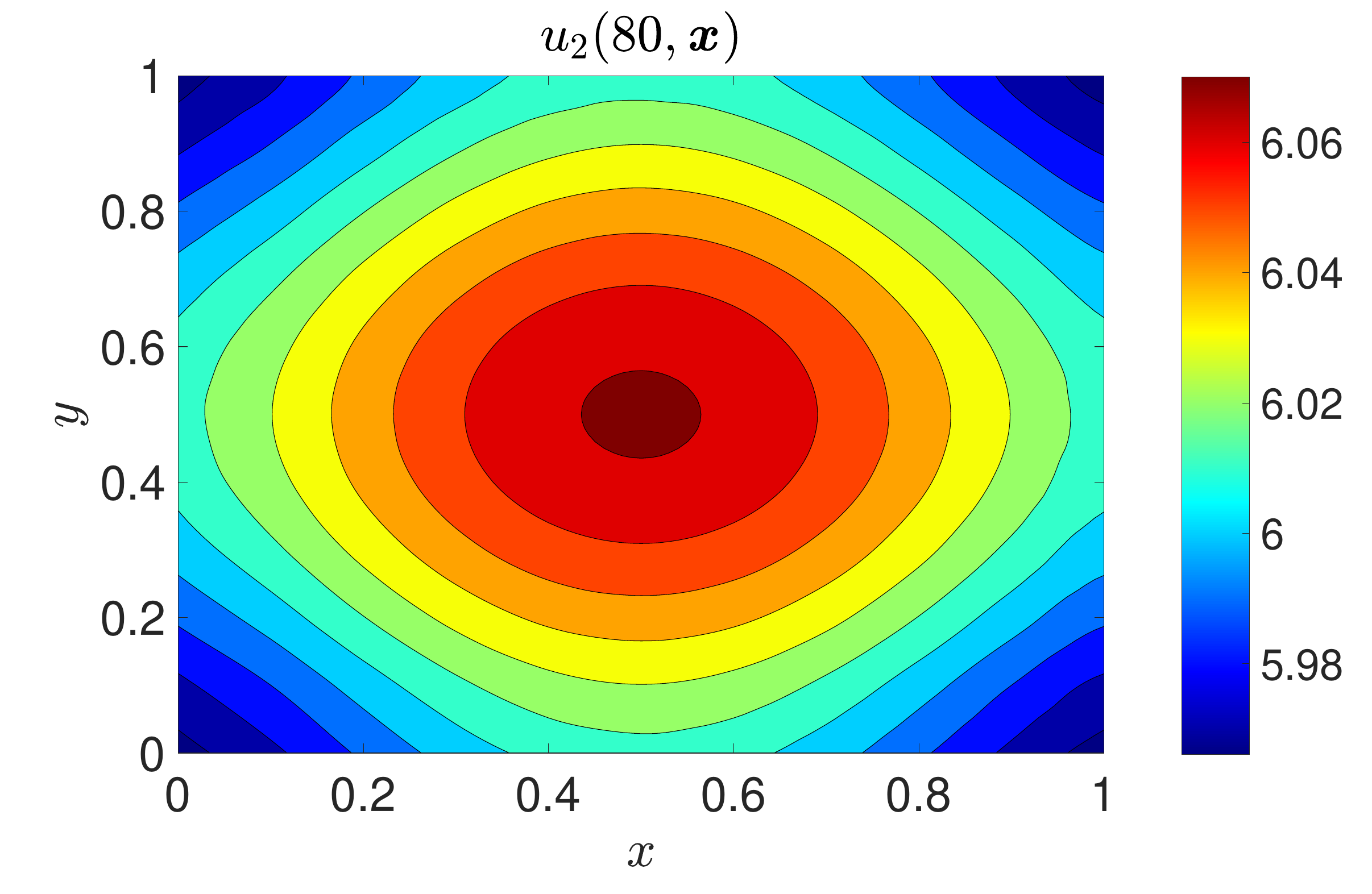}}
	\subfloat[]{\includegraphics[width=0.43\textwidth,height=0.27\textwidth]{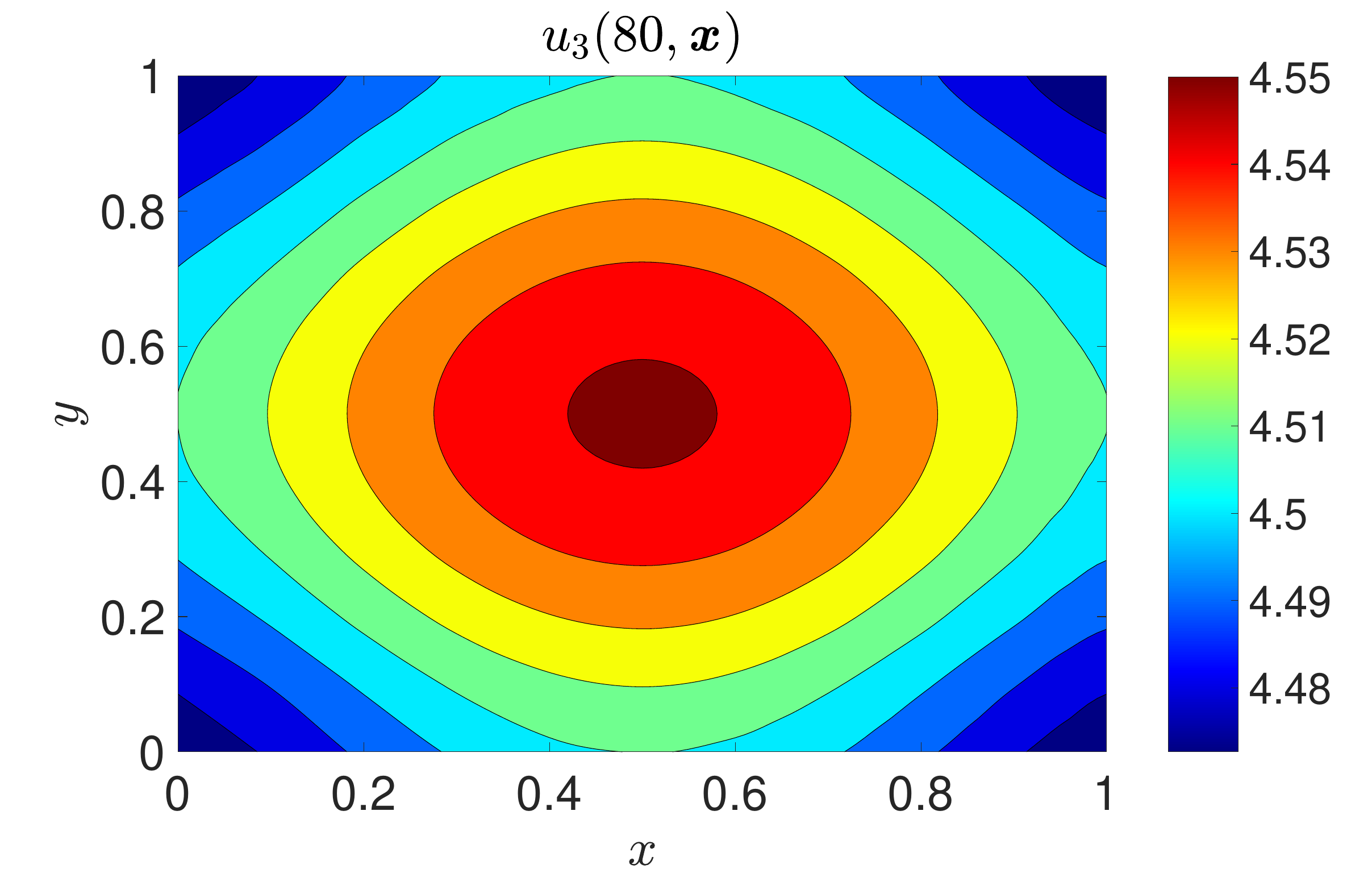}}
	\vspace{-4ex}\caption{ (a) Average density of each species,   (b) contour plot of species density $u_1$, (c) contour plot of species density $u_2$, and (d) contour plot of species density $u_3$ with the harvesting coefficients $\gamma_1=0.0009$, $\gamma_2=0.0036$, and $\gamma_3=0.0072$, diffusion coefficients $d_i=1$, $r_i=1$, and advection coefficients $\beta_i=0.001$, $i=\overline{1,3}$.}\vspace{-4ex}	\label{RDE-energy-g1_0_0009_g2_0_0036_g3_0_0072-K-time-exp}
\end{figure}

\begin{figure}
	\centering
	\subfloat[]{\includegraphics[width=0.43\textwidth,height=0.27\textwidth]{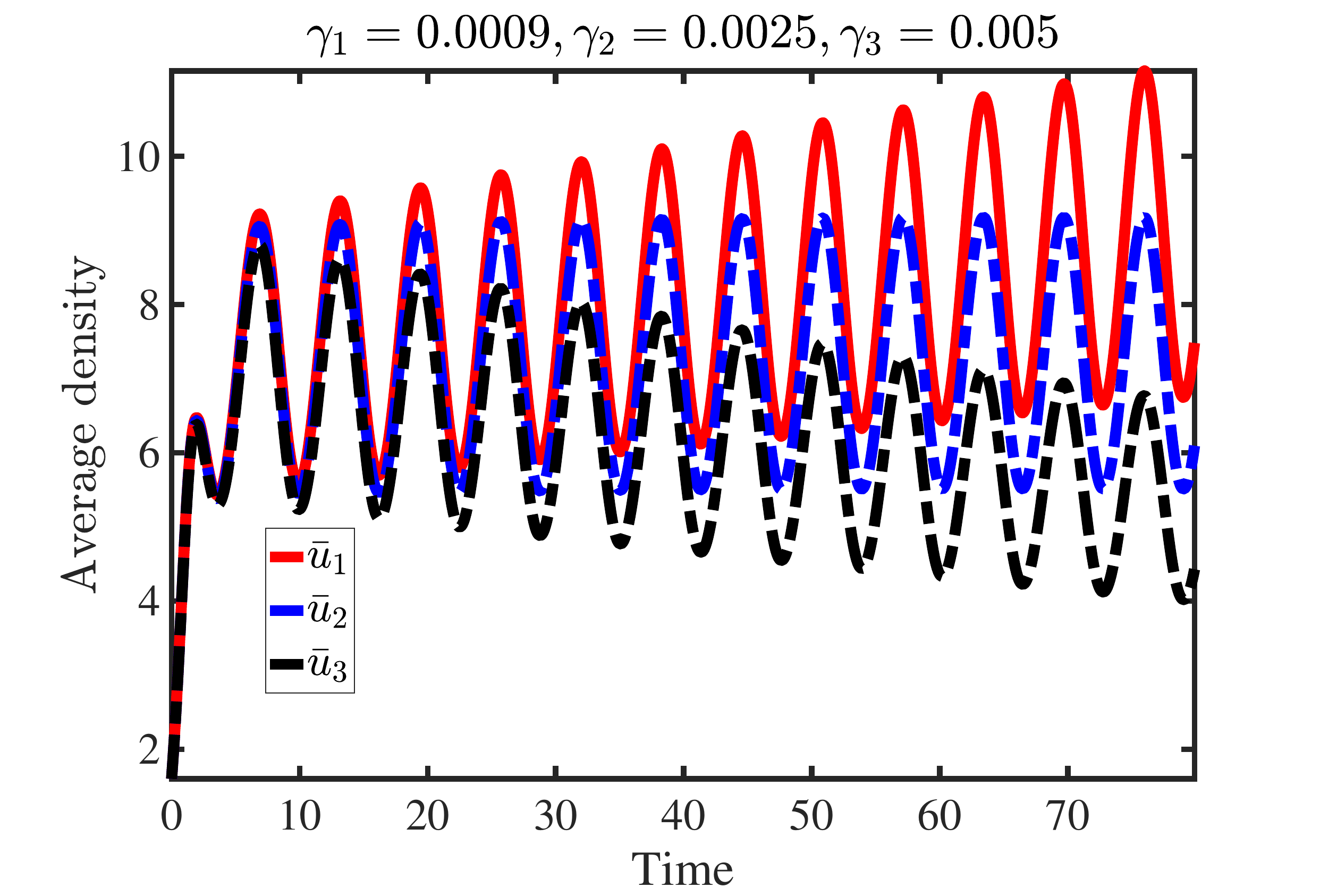}}
	\subfloat[]{\includegraphics[width=0.43\textwidth,height=0.27\textwidth]{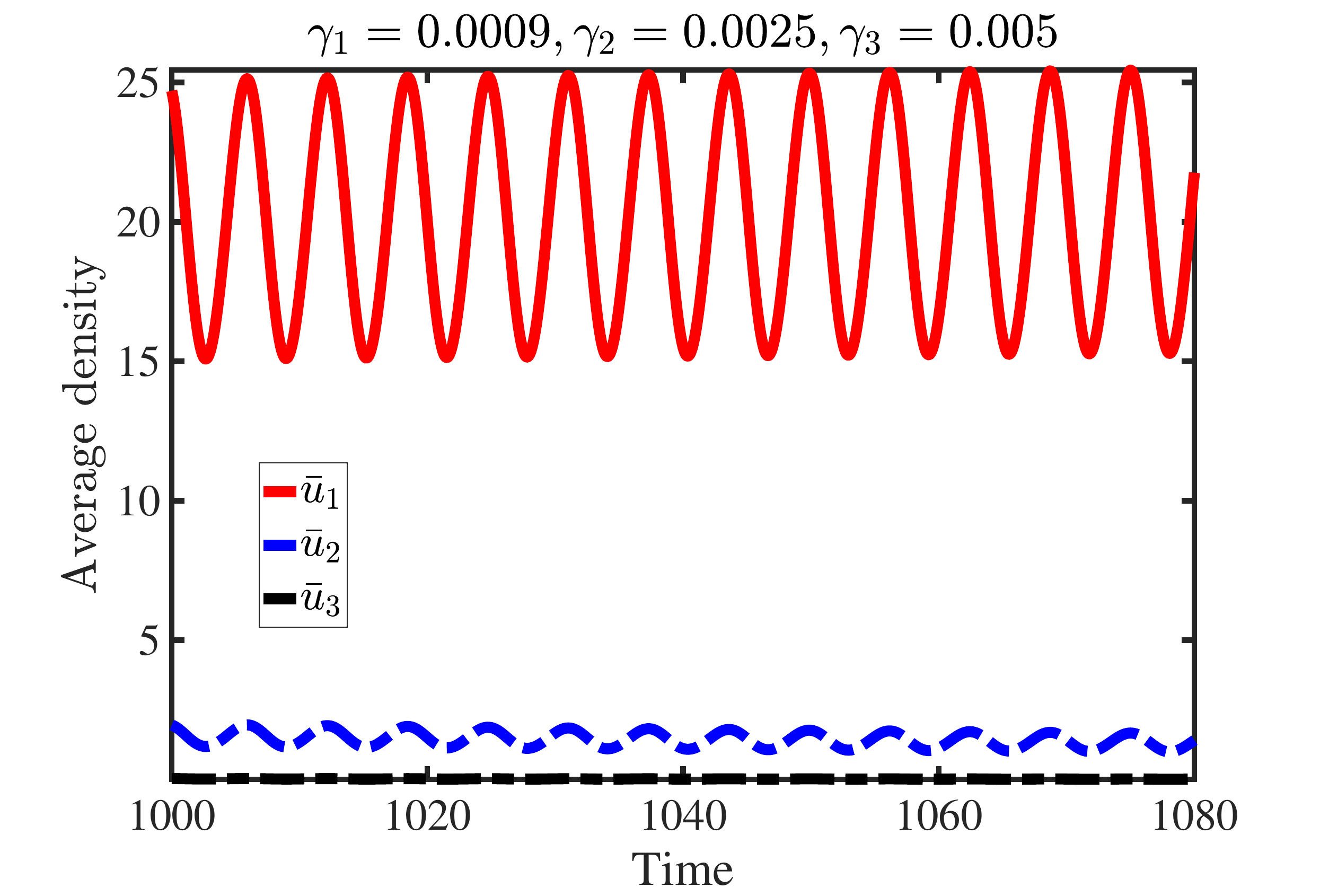}}
	\vspace{-4ex}\caption{ Average density of each species: (a) short-, and  (b) long-range with the harvesting coefficients $\gamma_1=0.0009$, $\gamma_2=0.0025$, and $\gamma_3=0.005$, diffusion coefficients $d_i=0.001$, advection coefficients $\beta_i=0.001$, and intrinsic growth rate $r_i=1.1+0.75\cos x\cos y$, for $i=\overline{1,3}$.}\vspace{-4ex}	\label{RDE-energy-g1_0_0009_g2_0_0025_g3_0_005-K-time-ex}
\end{figure}

In this test, we consider  constant intrinsic growth rates $r_i=1.0$, diffusion coefficients $d_i=1.0$, advection coefficients $\beta_i=0.001$,  $i=\overline{1,3}$ and varying harvesting coefficients $\gamma_1=0.0009$, $\gamma_2=0.0036$, and $\gamma_3=0.0072$, and run the simulation until $T=80$. In Figure \ref{RDE-energy-g1_0_0009_g2_0_0036_g3_0_0072-K-time-exp}(a), we represent the average population density over time.
In Figures \ref{RDE-energy-g1_0_0009_g2_0_0036_g3_0_0072-K-time-exp}, the average density of each species versus time is plotted for time $ t = 0 $ to $ 80 $, and the population density contour plot of each of the species at time $ t = 80 $. From the average density plot, we observe periodic population densities for all species, where the density of $ u_3 $ is decreasing because of its higher harvesting coefficient (Figure \ref{RDE-energy-g1_0_0009_g2_0_0036_g3_0_0072-K-time-exp}(a)). It is predicted that the species $ u_3 $ will die out if time is too large. From the contour plots, it is observed that the highest population density is at the point $ (0.5, 0.5) $ and there is a coexistence of all species, though the population density of the species $ u_1 $ remains bigger than the species $ u_2 $, and $ u_3 $ over the domain (Figure \ref{RDE-energy-g1_0_0009_g2_0_0036_g3_0_0072-K-time-exp}(b)-(d)). This happens because of different harvesting parameters, and the optimal value of the carrying capacity function is achieved at the point $ (0.5, 0.5) $, which shows the symmetric distribution of the population. Because of the reduction in the population density of the third species, the other species get more resources to grow, and a significant boost is observed in the second species’ density and a considerable amount of density increment is observed in the first species.
Because of the periodic resource function, we observe periodic behavior in all the population densities (Figure \ref{RDE-energy-g1_0_0009_g2_0_0025_g3_0_005-K-time-ex}). From both of the plots in Figure \ref{RDE-energy-g1_0_0009_g2_0_0025_g3_0_005-K-time-ex}(a), (b), we observe that as the harvesting parameter increases, the species density decreases over time. The species with lower harvesting rate will converge to the stable solution faster. Figure \ref{RDE-energy-g1_0_0009_g2_0_0025_g3_0_005-K-time-ex}(a), (b) are plotted for the same data, but for short and long time scenarios. We observe that the species with the highest harvesting rate is extinct, whereas the species with the lowest harvesting rate is the winner over the other species. In summary, the species with lower harvesting rate is the sole winner for multiple population competition and is independent of any choice of equal intrinsic growth rate and the initial population size.
\subsection{Three species competition model: Evolution of Population Density with Varying Diffusion Coefficients}\label{exp-4}
\begin{figure} 
	\centering
	\subfloat[]{\includegraphics[width=0.43\textwidth,height=0.27\textwidth]{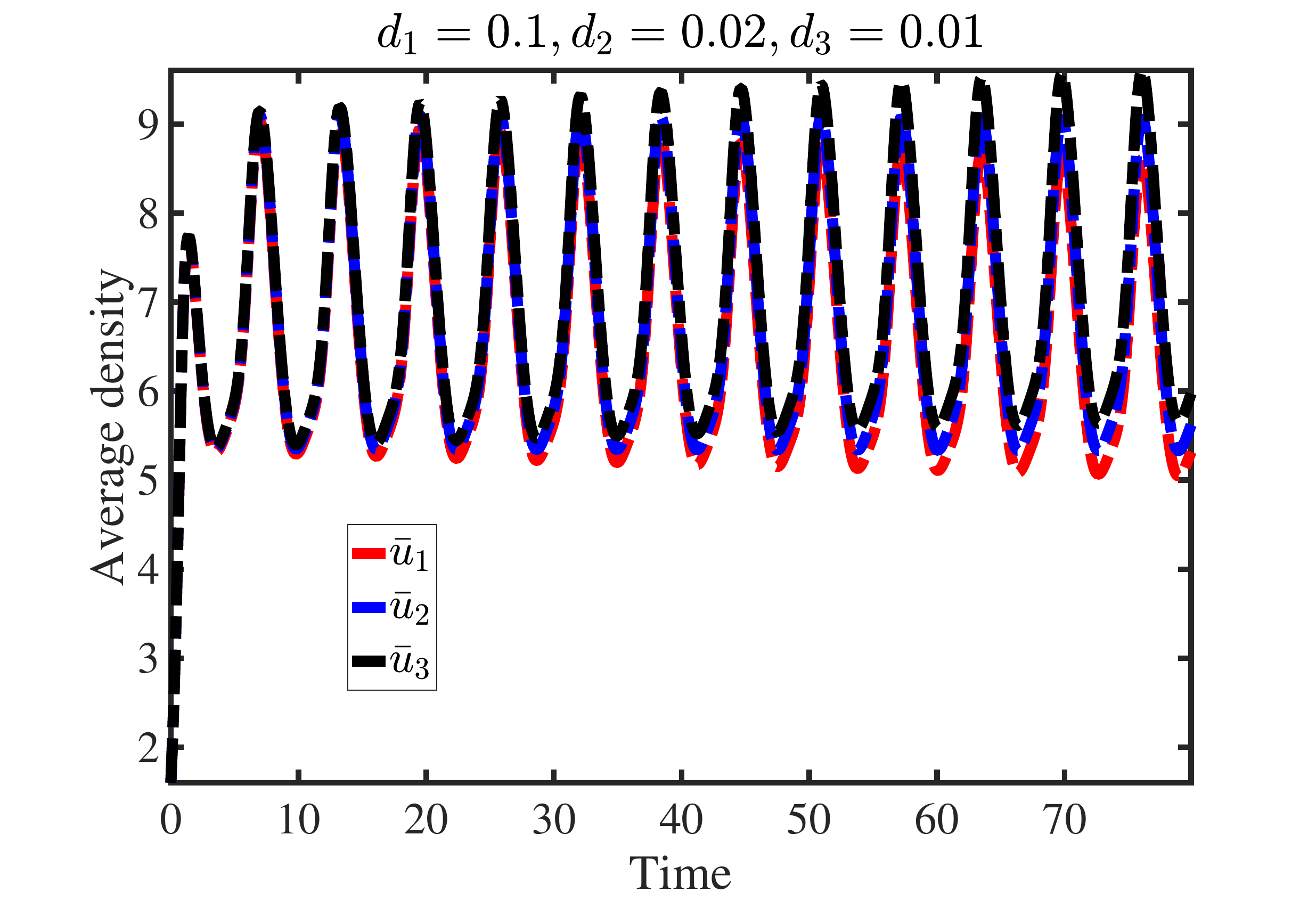}}
	\subfloat[]{\includegraphics[width=0.43\textwidth,height=0.27\textwidth]{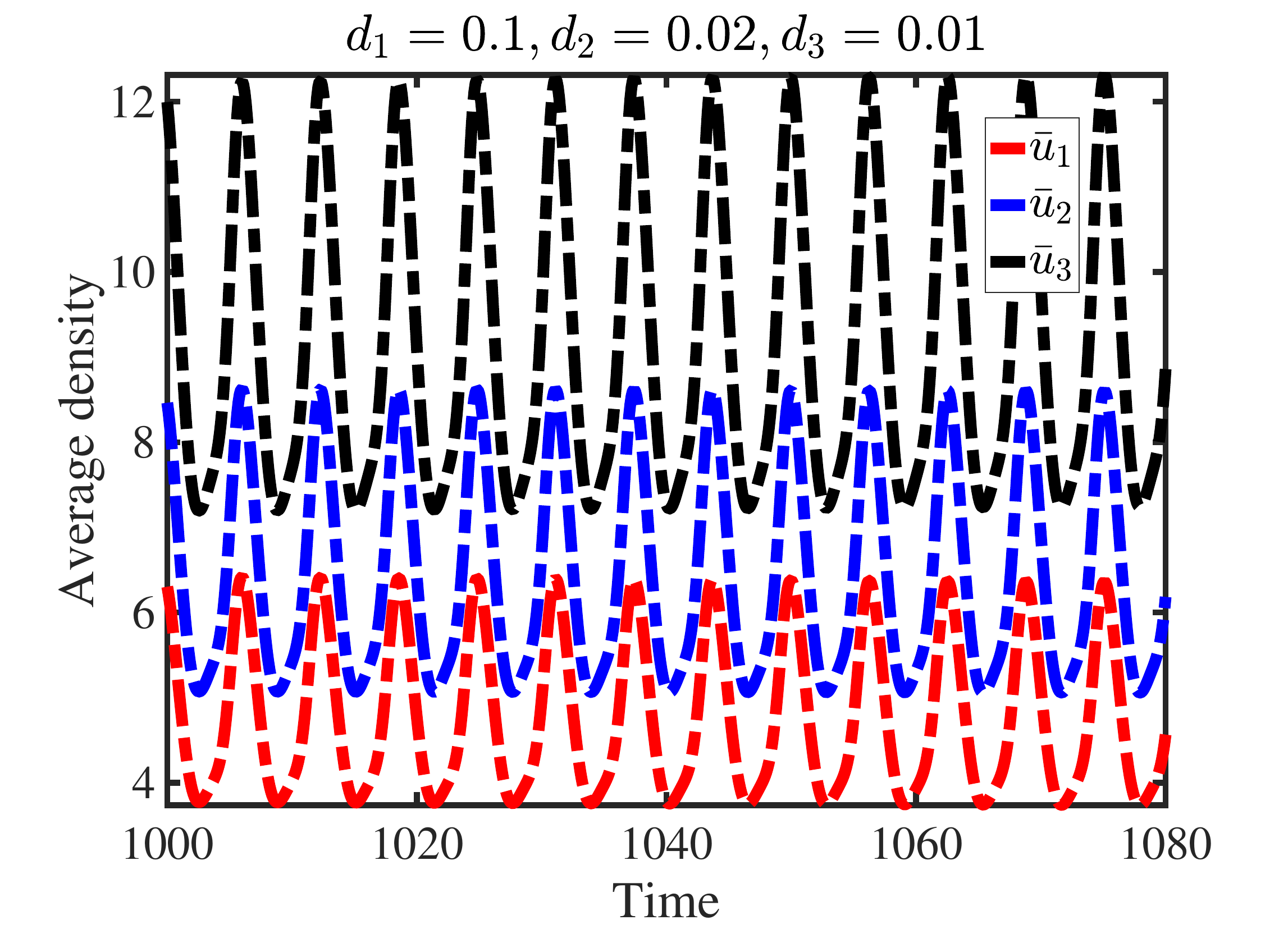}}
	\vspace{-4ex}\caption{Average density of each species: (a) short range, and  (b) long range with the harvesting coefficients $\gamma_1=0.0009$, $\gamma_2=0.0036$, and $\gamma_3=0.0072$, constant intrinsic growth rates $r_i=1$, and advection coefficients $\beta_i=0.001$, for $i=\overline{1,3}$.}\vspace{-4ex}	\label{RDE-energy-d1_0_1_d2_0_002_d3_0_001-K-time-exp}
\end{figure}

\begin{figure}
	\centering
	\subfloat[]{\includegraphics[width=0.43\textwidth,height=0.27\textwidth]{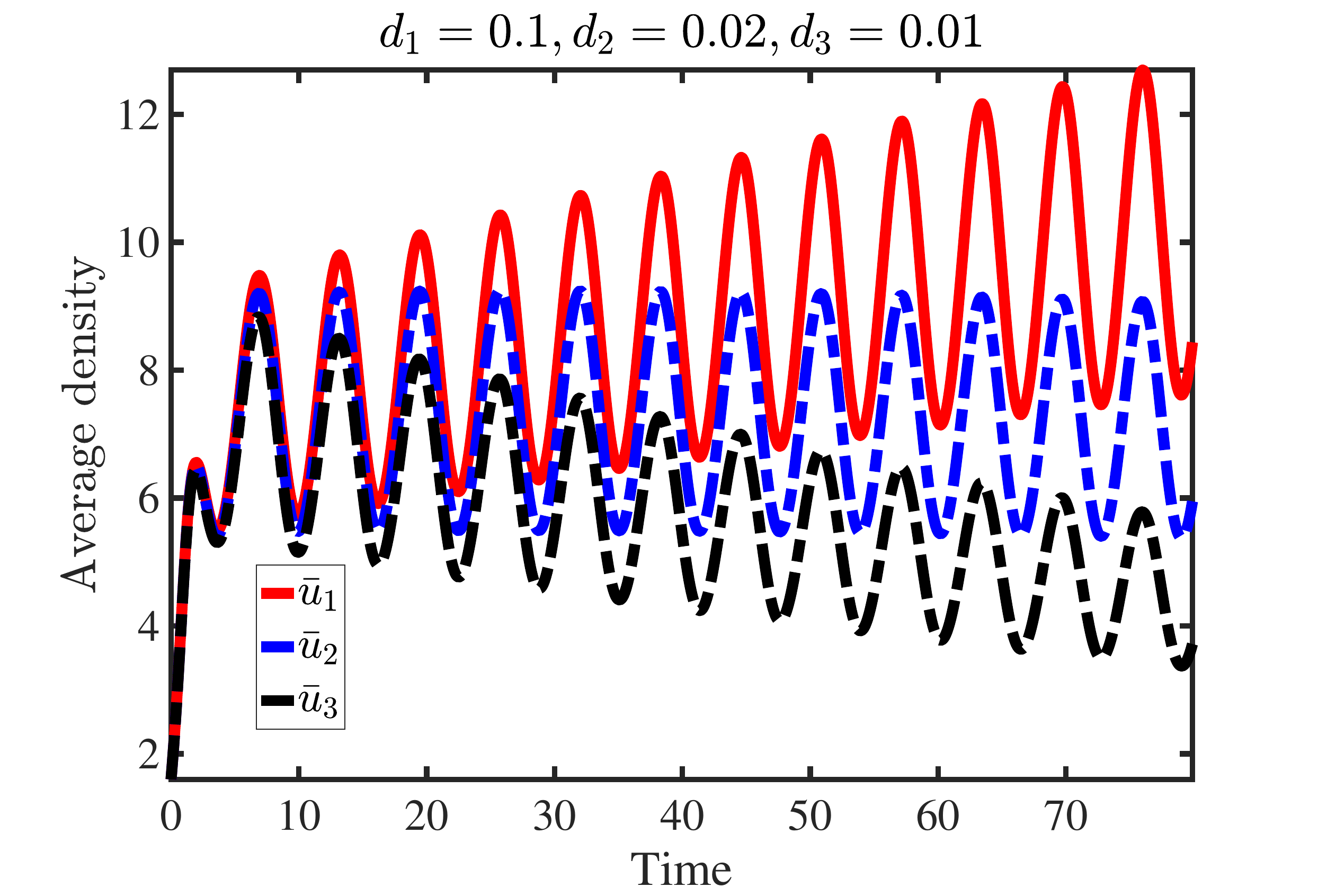}}
	\subfloat[]{\includegraphics[width=0.43\textwidth,height=0.27\textwidth]{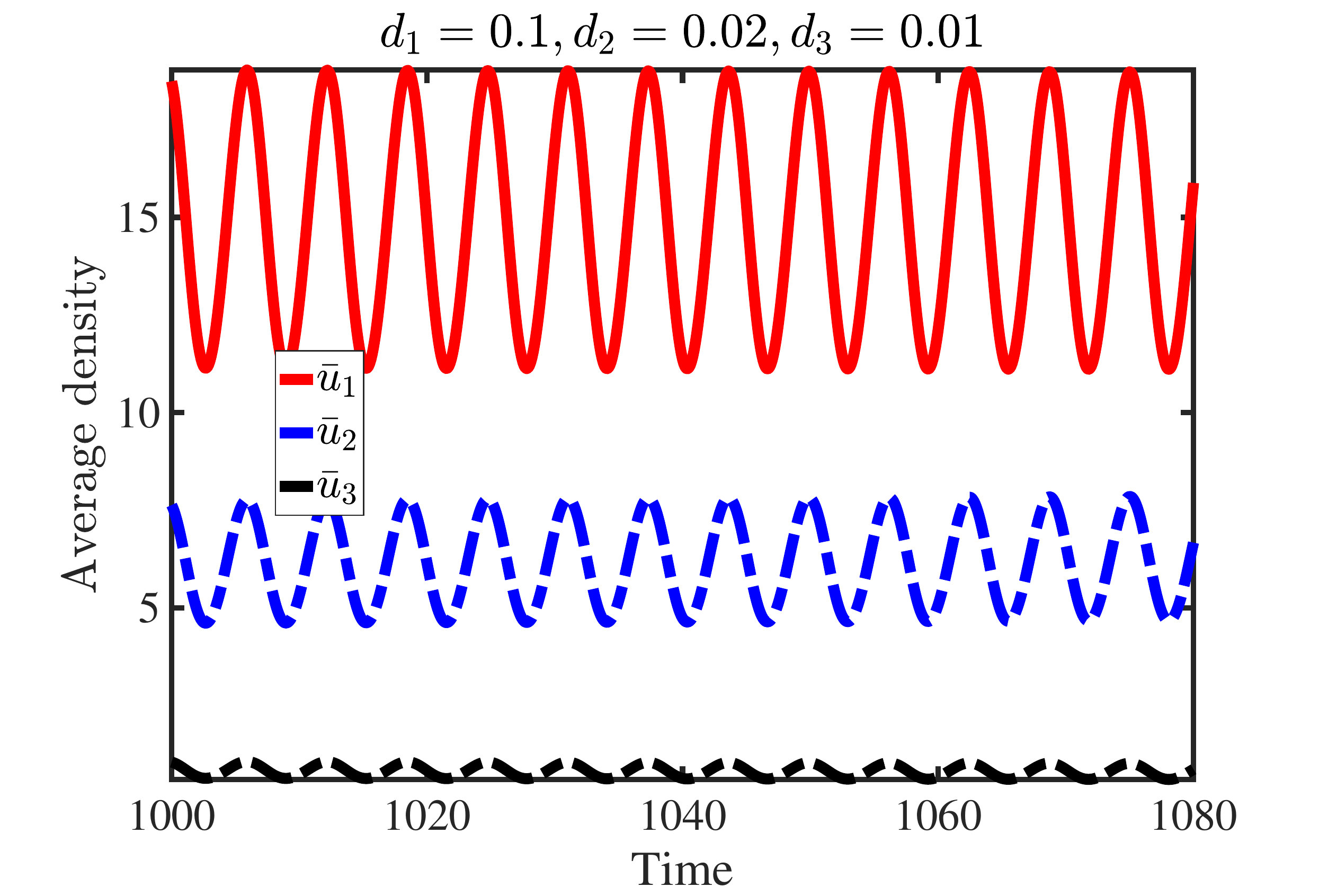}}
	\vspace{-4ex}\caption{Average density of each species: (a) short-, and  (b) long-range with the harvesting coefficients $\gamma_1=0.0009$, $\gamma_2=0.0036$, and $\gamma_3=0.0072$, constant intrinsic growth rates $r_i=1.1+0.75\cos x\cos y$, and advection coefficients $\beta_i=0.001$, for $i=\overline{1,3}$.}\vspace{-4ex}	\label{RDE-energy-d1_0_1_d2_0_002_d3_0_001-r-var}
\end{figure}

In this test, we consider  constant intrinsic growth rates $r_i=1.0$, advection coefficients $\beta_i=0.001$, varying harvesting coefficients $\gamma_i=0.009, 0.0036, 0.0072$,  $i=\overline{1,3}$ and diffusion coefficients $d_1=0.1$, $d_2=0.02$, and $d_3=0.01$, and run the simulation until $T=80$. In Figure \ref{RDE-energy-d1_0_1_d2_0_002_d3_0_001-K-time-exp}, we represent the average population density over time.
In Figures \ref{RDE-energy-d1_0_1_d2_0_002_d3_0_001-K-time-exp}(a), the average density of each species versus time is plotted for time $ t = 0 $ to $ 80 $. From the average density plot, we observe periodic population densities for all species, where the density of $ u_1 $ is decreasing because of its higher diffusion coefficient 
(Figure \ref{RDE-energy-d1_0_1_d2_0_002_d3_0_001-K-time-exp}(a)). It is predicted that the species $ u_1 $ will die out if time is too large.
Because of the periodic resource function, we observe periodic behavior in all the population densities (Figure \ref{RDE-energy-d1_0_1_d2_0_002_d3_0_001-r-var}). From both of the plots in Figure \ref{RDE-energy-d1_0_1_d2_0_002_d3_0_001-r-var}(a), (b), we observe that as the diffusion parameter increases, the species density decreases over time. The species with lower diffusion rate will converge to the stable solution faster. Figure \ref{RDE-energy-d1_0_1_d2_0_002_d3_0_001-r-var}(a), (b) are plotted for the same data, but for short and long time scenarios. We observe that the species with the highest diffusion rate is extinct, whereas the species with the lowest diffusion rate is the winner over the other species. In summary, the species with lower diffusion rate is the sole winner for multiple population competition and is independent of any choice of equal intrinsic growth rate and the initial population size.
\clearpage

	\section{Conclusion}\label{conclusion}
	 Stocking and harvesting strategies can have significant impact on population economics. In this paper, we propose an advection-reaction-diffusion based time evolutionary mathematical model, ARD-PDEs-SH, using a system of nonlinear coupled PDEs. The ARD-PDEs-SH model represents the dynamics of an $N$-species interaction subject to the effect of stocking and harvesting. 
 
 A continuous mathematical analysis of the ARD-PDEs-SH model is given for single species and two competing species in a heterogeneous advective environment with zero Neumann boundary conditions, i.e., no individuals can pass through the upstream end and the downstream end. In our investigation, we assumed that the stocking or harvesting rate was in $[0,1)$. The uniqueness and existence of the model solution have been proven and also, it has been shown that the solution of the model is positive. 
 
 We also proposed two fully discrete decoupled linearized stable methods for numerical approximation of the ARD-PDEs-SH model solutions. The stability and convergence theorems for the both algorithms are proven rigorously. The implicit-explicit schemes are found that the first scheme is first-order accurate and the second scheme is second-order accurate in time and both are optimally accurate in space. Numerical tests using manufactured analytical solution are given that verifies the predicted convergence rates. The efficient feature of these algorithms is that at each time step, each of the discrete linearized PDEs in the system can be solved simultaneously.   
		Numerical illustrations demonstrate:
		\begin{itemize}
			\item [(a)] If the harvesting rate increases, then the population density will decrease and the species with a higher harvesting rate will extinct; An opposite scenario is observed for the stocking parameter. 
			\item[(b)] For space-dependent carrying capacity, there is a periodic behavior in the solution.
			\item[(c)] With the absence of stocking and harvesting effect and a constant advection rate, the species with a higher diffusion rate will extinct. 
			\item[(d)] For different advection rates, there exists a coexistence among three species with the absence of stocking or harvesting effect and constant diffusion rate.
			\item[(e)] A strong coexistence is observed for space-dependent intrinsic growth rates. 
			\item[(f)] the higher the harvesting rate, the lower the density, and 
			\item[(g)] if the harvesting rate of one species surpass the intrinsic growth rate, then it will totally vanish and the other two will coexist.
		\end{itemize}
 A continuous analysis for the two species competition traveling wave solutions of \eqref{RDE1} will be the future research avenue. For the SH model with small diffusion, as future research work, we will propose and analyze FE methods following the work in \cite{john2008finite}. We will design the discrete algorithm elegantly so that at each time-step, the FE system matrix will be the same for all the species but the right-hand-side vector will be different \cite{mohebujjaman2017efficient, Mohebujjaman2022High}. This will reduce a huge computational complexity and will allow us to take advantage of block linear solvers.
\end{document}